\documentclass{aims}

\usepackage{graphicx}
\usepackage{amsmath}
\usepackage{amssymb}
\usepackage{enumerate}
\usepackage{color}
\usepackage{mathrsfs}
\usepackage[colorlinks=true]{hyperref}
\usepackage{longtable,booktabs}
\usepackage{subfigure}
  \usepackage{epsfig} 
\usepackage{graphicx}
\usepackage{epstopdf}
 \usepackage[colorlinks=true]{hyperref}
\hypersetup{urlcolor=blue, citecolor=red}
 \usepackage[pagewise]{lineno}
\def\currentvolume{X}
 \def\currentissue{X}
  \def\currentyear{200X}
   \def\currentmonth{XX}
    \def\ppages{X--XX}
     \def\DOI{10.3934/xx.xx.xx.xx}

\newtheorem{theorem}{Theorem}[section]
\newtheorem{corollary}[theorem]{Corollary}

\newtheorem{lemma}[theorem]{Lemma}
\newtheorem{proposition}[theorem]{Proposition}

\theoremstyle{definition}
\newtheorem{definition}[theorem]{Definition}
\newtheorem{example}{Example}[section]
\newtheorem{remark}{Remark}[section]

\DeclareMathOperator{\diag}{diag}

\newcommand{\ma}{\color{magenta}}

\title[Permanence and classification for discrete-time systems] 
      {Permanence and universal classification of discrete-time competitive systems via the carrying simplex }

\author
[M. Gyllenberg, J. Jiang, L. Niu and P. Yan]{}


 \email{mats.gyllenberg@helsinki.fi}
 \email{jiangjf@shnu.edu.cn}
 \email{lei.niu@helsinki.fi}
 \email{ping.yan@helsinki.fi}

\thanks{This work is supported by the National Natural Science Foundation of China (NSFC) under Grant No. 11371252 and Grant No. 11771295, Shanghai Gaofeng Project for University Academic Program Development, and the Academy of Finland.}

\thanks{$^*$ Corresponding author: Lei Niu}

\thanks{{\ma This article has been accepted for publication in Discrete and Continuous Dynamical Systems - Series A, published by the American Institute of Mathematical Sciences}}

\begin{document}
\maketitle

\centerline{\scshape Mats Gyllenberg}
\smallskip
{\footnotesize
 \centerline{Department of Mathematics and Statistics, University of Helsinki}
   \centerline{Helsinki FI-00014, Finland}
} 

\medskip

\centerline{\scshape Jifa Jiang}
\smallskip
{\footnotesize
 \centerline{Mathematics and Science College, Shanghai Normal University}
   \centerline{Shanghai 200234, China}
}
\medskip

\centerline{\scshape Lei Niu$^*$}
\smallskip
{\footnotesize
 \centerline{Department of Mathematics and Statistics, University of Helsinki}
   \centerline{Helsinki FI-00014, Finland}
}
\medskip
\centerline{\scshape Ping Yan}
\smallskip
{\footnotesize
 \centerline{School of Sciences, Zhejiang A $\&$ F University}
   \centerline{Hangzhou 311300, China}
  \centerline{Department of Mathematics and Statistics, University of Helsinki}
   \centerline{Helsinki FI-00014, Finland}
}

\bigskip

\centerline{(Communicated by the associate editor name)}
\medskip

\begin{abstract}
We study the permanence and impermanence for discrete-time Kolmogorov systems admitting a carrying simplex. Sufficient conditions to guarantee permanence and impermanence are provided based on the existence of a carrying simplex. Particularly, for low-dimensional systems, permanence and impermanence can be determined by boundary fixed points. For a class of competitive systems whose fixed points are determined by linear equations, there always exists a carrying simplex. We provide a universal classification via the equivalence relation relative to local dynamics of boundary fixed points for the three-dimensional systems by the index formula on the carrying simplex.  There are a total of $33$ stable equivalence classes which are described in terms of inequalities on parameters, and we present the phase portraits on their carrying simplices. Moreover, every orbit converges to some fixed point in classes $1-25$ and $33$; there is always a heteroclinic cycle in class $27$; Neimark-Sacker bifurcations may occur in classes $26-31$ but cannot occur in class $32$. Based on our permanence criteria and the equivalence classification, we obtain the specific conditions on parameters for permanence and impermanence. Only systems in classes $29,31,33$ and those in class $27$ with a repelling heteroclinic cycle are permanent. Applications to discrete population models including the Leslie-Gower models, Atkinson-Allen models and Ricker models are given.
\end{abstract}

\section{Introduction}
The theory of the carrying simplex was first developed in \cite{hirsch1988} by Hirsch for continuous-time competitive systems of Kolmogorov ODEs, which states that every strongly competitive and dissipative system for which  the origin is a repeller possesses a globally attracting hypersurface $\Sigma$ of codimension one, called the carrying simplex in \cite{Zeeman1994,Z993}. Furthermore,  $\Sigma$ is homeomorphic to the $(n-1)$-dimensional standard probability simplex $\Delta^{n-1}=\{x\in \mathbb{R}_+^n:\sum_i x_i=1\}$ by radial projection, and has the property that every nontrivial orbit in the nonnegative cone $\mathbb{R}_+^n$ is asymptotic to one in $\Sigma$. It has been proved as a powerful tool to investigate global dynamics of competitive systems, especially for the lower dimensional systems. The reader can consult, for instance, \cite{HofSo,zeeman1998three,xiao2000,zz2,zz1,lu2003three,Gyllenberg06,Baigent2012Global,Ba,jnz,Hou2015GLOBAL,JN-Siam}, for more results on continuous-time competitive systems via the carrying simplex.

The theory of the carrying simplex has been extended to discrete-time competitive systems due to the early work of de Mottoni and Schiaffino \cite{Mottoni1981}, Hale and Somolinos \cite{Hale1983} and Smith \cite{smith1986}. By introducing mild conditions, Wang and Jiang proved in \cite{wang2002uniqueness} the existence of a carrying simplex for competitive mappings, which solved the conjecture on the carrying simplex proposed by Smith in \cite{smith1986} and was improved further by Diekman, Wang and Yan \cite{diekmann2008carrying}. Further, Hirsch announced a theory on the existence of a carrying simplex in \cite{hirsch2008existence} for the continuous Kolmogorov map (not necessarily invertible)
\begin{equation}\label{map-T}
T(x)=(x_1F_1(x),\ldots,x_nF_n(x)),\quad x\in  \mathbb{R}_+^n,
\end{equation}
where $F_i$ are continuous satisfying $F_i(x)>0$ for all $x\in \mathbb{R}_+^n$, $i=1,\ldots,n$. The discrete-time dynamical system induced by such $T$ has been much used to describe the interactions of $n$ species with non-overlapping generations; see \cite{May1974,Hutson1982,hofbauer1987coexistence,Franke1991,Rees1997,Levine2002,Kon2004,hirsch2008existence}. The statement of Hirsch's Theorem in \cite{hirsch2008existence} was rigorously proved by
Ruiz-Herrera \cite{ruiz2011exclusion} under similar assumptions to Hirsch's. Other criteria on the existence of a carrying simplex for Kolmogorov mappings were also established in \cite{Baigent2015,jiang2015,LG}, and we refer the readers to the paper \cite{LG} for a review.

The importance of the existence of a carrying simplex stems from the fact that it captures the relevant long-term dynamics. In particular, it contains all non-trivial fixed points, periodic orbits, invariant circles and heteroclinic cycles (see, for example, \cite{jiang2014,jiang2015, LG, GyllenbergCGAA}).
Therefore, the common approach in the study of these systems is to focus on the dynamics on the carrying simplex.

 In \cite{ruiz2011exclusion}, Ruiz-Herrera provided an exclusion criterion for discrete-time competitive models of two or three species via the carrying simplex. Jiang and Niu derived a fixed point index formula on the carrying simplex for three-dimensional maps in \cite{jiang2014}, which states that the sum of the indices of all fixed points on the carrying simplex is one. Based on this formula, an alternative classification for $3$-dimensional Atkinson-Allen models was provided in \cite{jiang2014} and an alternative classification  for $3$-dimensional Leslie-Gower models was provided in \cite{LG}. Such a classification has also been given for the $3$-dimensional generalized Atkinson-Allen models \cite{GyllenbergCGAA} and Ricker models admitting a carrying simplex \cite{GyllenbergRicker}, respectively. Jiang, Niu and Wang \cite{jiang2015} studied the occurrence and stability of heteroclinic cycles for competitive maps with a carrying simplex. Recently, Niu and Ruiz-Herrera proved in \cite{nonlinearity2018} that every orbit converges to a fixed point for three-dimensional maps with a carrying simplex when there is a unique positive fixed point such that its index is $-1$. For the geometrical properties of the carrying simplex and their impact on the dynamics, we refer the readers to \cite{Baigent2015, Baigent2019, Baigent2017,Mierczynski2018b,Mierczynski2018a}.

In population biology, the question of persistence of interacting species is one of the most important. There has been many papers on permanence for discrete-time Kolmogorov systems; see \cite{Hutson1982,hofbauer1987coexistence,Lu1999,Kon2001,Kon2004}. Here we study the permanence and impermanence for the discrete-time Kolmogorov systems \eqref{map-T} admitting a carrying simplex. The main mathematical tools involve average Liapunov functions (see \cite{Hutson1982,jiang2015}) and the theory of the carrying simplex. A successful case of combining these two approaches is the stability criterion for heteroclinic cycles established by Jiang, Niu and Wang in \cite{jiang2015}. Our project here is to provide the minimal conditions to ensure the permanence for such systems via the carrying simplex. Criteria for permanence and impermanence which are simple and easy to apply are provided for systems \eqref{map-T} admitting a carrying simplex. In particular, for three-dimensional systems, permanence and impermanence can be determined by boundary fixed points. As a special case, when the boundary of the carrying simplex is a heteroclinic cycle, the system is permanent if the heteroclinic cycle is repelling, while the system is impermanent if it is attracting.

Finally, we restrict attention to the class of maps given by
\begin{equation}\label{equ:T1}
    T_i(x)=x_if_i((Ax^\tau)_i,r_i),\quad i=1,\cdots,n,
\end{equation}
where $r_i>0$, $A$ is an ${n\times n}$ matrix with entries $a_{ij}>0$, $f_i:\mathbb{R}_+\times \dot{\mathbb{R}}_+\mapsto \dot{\mathbb{R}}_+$ are $C^1$, and $\tau$ denotes transpose. In addition, in this article, we always assume that $f_i$ satisfies
\begin{equation}\label{cons:f}
\begin{array}{ll}
 {\rm(i)}~\displaystyle f_i(r,r)=1,~\frac{\partial f_i(z,r)}{\partial z}<0,\quad \forall (z,r)\in \mathbb{R}_+\times \dot{\mathbb{R}}_+;& \\
\noalign{\medskip}
{\rm(ii)}~\displaystyle f_i(z,r)+z\frac{\partial f_i(z,r)}{\partial z}>0,\quad \forall (z,r)\in \mathbb{R}_+\times \dot{\mathbb{R}}_+.&
\end{array}
\end{equation}
Note that $f_i$ enjoys the properties: $f_i(z,r)>1$ for $z<r$, $f_i(z,r)=1$ for $z=r$, and $f_i(z,r)<1$ for $z>r$.

The discrete-time system induced by map \eqref{equ:T1} is often used in the modeling of $n$-species in competition; see \cite{Hassell1976,Law1987,Franke1991,Franke1992,Rees1997,Levine2002,roeger2004discrete,cushing2004some,Kon2006Convex,H2014Qualitative}. The variable $x_i$ is the density of species $i$ and $f_i$ is its growth function (or fitness function). Per capita growth rate for species $i$ in the absence of competition is given by $f_i(0,r_i)$.  The parameter $a_{ij}$ is the competition coefficient which measures the effect of species $j$ relative to $i$ on the function $f_i$.

Note that when we consider $T$ restricted to the $i$-th coordinate axis, we have $g_i(x_i)=x_if_i(a_{ii}x_i,r_i)$, which describes the dynamics of the $i$-th species without inter-specific competition. By \eqref{cons:f} (i), $f_i(a_{ii}x_i,r_i)$ is a decreasing function of $x_i$. The biological meaning is that the per capita growth rate of species $i$ is a decreasing function of population density due to negative density dependent mechanism such as intra-specific competition \cite{Fishman1997,Townsend3}.
By \eqref{cons:f} (ii), $g_i(x_i)$ is an increasing function of $x_i$, so the population dynamics of each species is monotone. Biologically, this means that the intra-specific competition is contest due to increasing utilization of available resources, where each successful competitor gets all resources it requires for survival or reproduction (see \cite{Varley1973,Hassell1975,Br2005,W2006Contest}). Furthermore, it follows from \eqref{cons:f} that $f_i(0,r_i)>1$ and
\begin{equation}\label{cons:f-2}
f_i(r_i,r_i)+r_i\frac{\partial f_i}{\partial z}(r_i,r_i)>0,~\mathrm{i.e.}~1+r_i\frac{\partial f_i}{\partial z}(r_i,r_i)>0,
\end{equation}
so $0$ is a repeller (growth of small populations), and $x_i=\frac{r_i}{a_{ii}}$ is an attracting fixed point for $g_i$.

Many functions satisfy \eqref{cons:f}, such as
$$
\begin{array}{l}
  f(z,r)=(\frac{1+r}{1+z})^s\quad (0<s\leq 1), \\
  \noalign{\medskip}
  f(z,r)=\frac{1+r^s}{1+z^s}\quad (0<s\leq 1), \\
  \noalign{\medskip}
  f(z,r)=\frac{(1+r)(1-c)}{1+z}+c\quad (0<c<1), \\
  \noalign{\medskip}
  f(z,r)=\frac{1+\ln(1+r)}{1+\ln(1+z)}.
\end{array}
$$
We refer the readers to \cite{Gyllenberg1997,Geritz-Kisdi2004,Br2005,Eskola2007} for the mechanistic derivation of various discrete-time single-species population models with such growth functions.

By \eqref{cons:f} (i), the fixed points of $T$ are determined by the linear algebraic equations
\begin{equation}
 x_i=0~~\textrm{or}
    ~~(Ax^\tau)_i=r_i,\quad i=1,\ldots,n.
\end{equation}
We call the map \eqref{equ:T1} {\it fixed points linearly determined}. Jiang and Niu proved in \cite{LG} that all maps given by \eqref{equ:T1} with $f_i$ satisfying \eqref{cons:f} admit a carrying simplex  unconditionally. In this article, we focus on studying the parameter conditions that guarantee permanence and impermanence for the three-dimensional map \eqref{equ:T1} with given functions $f_i$. By noticing the above linear structure, we can define an equivalence relation on the parameter space for the three-dimensional map \eqref{equ:T1} as that for the two specific cases: the Atkinson-Allen model \cite{jiang2014} and Leslie-Gower model \cite{LG}. Two maps \eqref{equ:T1} are said to be equivalent relative to the boundary of $\Sigma$ if their boundary fixed points have the same locally dynamical property on $\Sigma$ after a permutation of the indices $\{1,2,3\}$.
Map \eqref{equ:T1} is said to be stable relative to the boundary of $\Sigma$ if
all the fixed points on the boundary are hyperbolic.
Via the index formula on the carrying simplex established in \cite{jiang2014}, we list the equivalence classes for all stable maps in Table \ref{biao0}. There are always a total of $33$ stable equivalence classes which can be described in terms of inequalities on parameters, and the equivalence  classification is independent of the choice of functions $f_i$, which presents a clear picture of the
essence of the dynamics for this class of maps. Moreover, based on this classification, one can easily get the parameter conditions that guarantee permanence and impermanence for such systems. Specifically, applying the permanence criteria to each class, we obtain that systems in classes $29,31,33$ and those in class $27$ with repelling heteroclinic cycles are permanent, while systems in classes $1-26$, $28$, $30$, $32$ and class $27$ with attracting heteroclinic cycles are impermanent. {For systems in class $33$, the permanence can guarantee the global stability of the positive fixed point.} It is emphasized that the fixed points linearly determined systems \eqref{equ:T1} contain many classical systems, such as the Atkinson-Allen model \cite{jiang2014} and the Leslie-Gower model \cite{LG}; see Section \ref{sec:application} for more. Our investigations will stimulate further study on the global behavior of these systems including the higher order bifurcations, multiplicity of closed invariant curves, and so on.

As we mentioned earlier, many of the important ideas used in this work
 are due to Hutson and Moran \cite{Hutson1982}, Hofbauer, Hutson and Jansen \cite{hofbauer1987coexistence}, Jiang and Niu \cite{jiang2014,LG} and Jiang, Niu and Wang \cite{jiang2015}.

The paper is organized as follows. In Section \ref{sec:notation}, we present some notions and recall some known results on average Liapunov functions. In Section \ref{sec:permanence}, we provide the criteria on permanence and impermanence for dissipative Kolmogorov systems and maps admitting a carrying simplex. In Section \ref{sec:extensions}, we define the equivalence relation relative to the boundary dynamics on the parameter space for the three-dimensional map \eqref{equ:T1}, and derive the $33$ stable equivalence classes. Based on this classification, we obtain the parameter conditions that guarantee permanence and impermanence. In Section \ref{sec:application}, we apply our results to some classical discrete population models including the Leslie-Gower models, Atkinson-Allen models and Ricker models. The paper ends with a discussion in Section \ref{sec:Discussion}.

\section{Notation and preliminaries}\label{sec:notation}
Suppose that $X$ is a metric space with metric $d_X(\cdot,\cdot)$ and $T: X\to X$ is a continuous mapping. Let $\mathbf{Z}_+=\{0,1,2,\ldots\}$.
For any $x\in X$, we define the positive orbit through $x$ as $\gamma^+(x):=\{T^kx: k\in \mathbf{Z}_+\}$, and denote the tail from the moment $m\geq 1$ of $\gamma^+(x)$ by $\gamma^+_m(x):=\{T^kx: k\geq m\}$.
A negative orbit through $x$ is a sequence $\{x(-k):k\in \mathbf{Z}_+\}$ such that $x(0)= x$, $Tx(-k-1)= x(-k)$ for all $k\in \mathbf{Z}_+$. The omega limit set $\omega(x):=\bigcap_{k\geq 0}\overline{\bigcup_{m\geq k}T^mx}$ of $x$ is the set of limit points of the positive orbit $\gamma^+(x)$. The alpha limit set $\alpha(x):=\bigcap_{k\geq 0}\overline{\bigcup_{m\geq k}x(-m)}$ associated to a negative orbit $\{x(-k):k\in \mathbf{Z}_+\}$ through $x$ is the set of limit points of this negative orbit.

For a subset $D\subseteq X$, $\overline{D}$ and $D^c$ denote the closure of $D$ in $X$ and the complement of $D$ respectively. Given any set $D$, let
$$
\gamma^+(D)=\bigcup_{x\in D}\gamma^+(x),\qquad \Omega(D)=\overline{\bigcup_{x\in D}\omega(x)}.
$$
Note that $\Omega(D)$ is a subset of the omega limit set $\omega(D)$ of the orbit through $D$, i.e.
$$
\Omega(D)\subseteq  \omega(D):=\bigcap_{k\geq 0}\overline{\bigcup_{m\geq k}T^mD}.
$$
We denote $\mathcal{E}(T)=\{x\in X:Tx=x\}$ to be the set of the fixed points of $T$.

 A set $D\subseteq X$ is called positively invariant (with respect to $T$) if $TD\subseteq D$; negatively invariant if $TD\supseteq D$; and invariant if $TD=D$.

  A set $J$ is said to attract a set $D$ under $T$ if for any $\varepsilon > 0$, there exists a
$k_0\geq 1$ such that $T^kD$ belongs to the $\varepsilon$-neighborhood $O_\varepsilon(J)$ of $J$ for $k \geq k_0$, where $O_\varepsilon(J)=\{y\in X:d_X(y,J)<\varepsilon\}$. A nonempty compact invariant set $J\subseteq X$ is said to be a global attractor of $T$ if $J$ attracts each bounded set $D \subseteq X$. $T$ is said to be dissipative if it admits a global attractor.

We say that a nonempty, compact, positively invariant subset $S\subseteq X$ repels for $T$ if there exists an $\varepsilon$-neighborhood $O_\varepsilon(S)$ of $S$ such that for all $x\in X\setminus S$ there exists a $k_0=k_0(x)>0$ satisfying $T^kx\notin O_\varepsilon(S)$ for all $k \geq k_0$, i.e. $\gamma^+_{k_0}(x)\subseteq X\setminus O_\varepsilon(S)$.

\smallskip

In the following part of this section, we recall two known results on average Liapunov functions.

\begin{lemma}[Lemma 2.1 in \cite{hofbauer1987coexistence}]\label{lemma:absorbing}
Let $T:X\mapsto X$ be a continuous map, where $X$ is a metric space.
Let $W$ be open with compact closure, and suppose that $U$ is open and positively invariant, where $\overline{W}\subseteq U \subseteq X$. If $\gamma^+(x) \cap W\neq \emptyset$ for any $x\in U$, then $\gamma^+(\overline{W})$ is compact and positively invariant such that for any $x\in U$, there exists a $m=m(x)>0$ satisfying $\gamma^+_m(x)\subseteq \gamma^+(\overline{W})$.
\end{lemma}

\begin{lemma}[Repelling Set \cite{Hutson1982}]\label{repellor_T}
Suppose that $M$ is a compact metric space and $T:M\mapsto M$ is a continuous mapping. Assume that $S$ is a compact subset of $M$ with empty interior such that $S$ and $M\setminus S$ are positively invariant under $T$. Suppose that there is a continuous function $V: M\mapsto \mathbb{R}_+$ satisfying that

  {\rm (i)}~$V(x)=0\Leftrightarrow x\in S$,

  {\rm (ii)}~$\displaystyle\psi_M(x):=\sup_{k>0}\vartheta_M(k,x)>1$ for all $x\in S$, where
 \begin{equation}
 \vartheta_M(k,x):=\liminf_{\mbox{\tiny$\begin{array}{c}
y\rightarrow x\\
y\in M\setminus S
\end{array}$}}\!\frac{V(T^k y)}{V(y)},~~k\in \mathbf{Z}_+.
 \end{equation}
Then $S$ repels for $T$.
\end{lemma}
\begin{remark}\label{con:repel}
The function $V$ in Lemma \ref{repellor_T} is called an average Liapunov function. Moreover, by Corollary $2.2$ in \cite{Hutson1982}, the condition {\rm (ii)} in Lemma \ref{repellor_T} is implied by the following condition

  {\rm (ii')}~$\psi_M(x)>1$ for all $x\in \Omega(S)$, and $\psi_M(x)>0$ for all $x\in S$.
\end{remark}

\begin{lemma}[Attracting Set  \cite{jiang2015}]\label{attractor_T}
Suppose that $M$ is a compact metric space and $T:M\mapsto M$ is a continuous mapping. Assume that $S$ is a compact subset of $M$ with empty interior such that $S$ and $M\setminus S$ are positively invariant under $T$. Suppose that there is a continuous function $V: M\mapsto \mathbb{R}_+$ and a constant $C>0$ satisfying that

{\rm (i)}  $V(x)=0\Leftrightarrow x\in S$, and $\frac{V(Tx)}{V(x)}\leq C$ for all $x\in M\setminus S$,

{\rm (ii)}$\displaystyle\varphi_M(x):=\inf_{k>0}\zeta_M(k,x)<1$ for all $x\in \Omega(S)$, where
\begin{equation}
\zeta_M(k,x):=\limsup_{\mbox{\tiny$\begin{array}{c}
y\rightarrow x\\
y\in M\setminus S
\end{array}$}}\!\frac{V(T^k y)}{V(y)},~~k\in \mathbf{Z}_+.
\end{equation}
Then $S$ attracts for $T$, that is, there is an $\varepsilon$-neighborhood $O_\varepsilon(S)$ of $S$ such that $\omega(x)\subseteq S$ for every $x \in O_\varepsilon(S)$.
\end{lemma}

The following proposition shows that the presence of an invariant repelling set $S$ implies the existence of a compact set $K$ contained in $S^c$ that absorbs the orbits contained in $S^c$. See \cite{Smith2010} for more details.
\begin{proposition}\label{prop:repel}
Suppose that $M$ is a compact metric space and $T:M\mapsto M$ is a continuous mapping. Assume that $S\subseteq M$ is compact with empty interior such that $S$ and $M\setminus S$ are positively invariant under $T$. If $S$ repels, then there is a compact positively invariant set $K\subseteq M\setminus S$ such that for every $x\in M\setminus S$, there exists a $m=m(x)>0$ satisfying $\gamma^+_m(x)\subseteq  K$.
\end{proposition}
\begin{proof}
Set $U=M\setminus S$. Since $S$ repels, there exists an $\varepsilon$-neighborhood $O_\varepsilon(S)$ of $S$ such that for any $x\in U$ there exists a $k=k(x)>0$ satisfying $\gamma^+_{k}(x)\subseteq M\setminus O_\varepsilon(S)$, and hence $\gamma^+_{k}(x)\subseteq W$, where $W=M\setminus \overline{O_{\frac{\varepsilon}{2}}(S)}$.  Thus, $\gamma^+(x) \cap W\neq \emptyset$ for any $x\in U$. Note that $W$ is open and $\overline{W}\subseteq U$ is compact, so it follows from Lemma \ref{lemma:absorbing} that $\gamma^+(\overline{W})$ is compact, positively invariant such that for any $x\in U$, there exists a $m=m(x)>0$ satisfying $\gamma^+_m(x)\subseteq \gamma^+(\overline{W})$. Let $K=\gamma^+(\overline{W})$, which is the desired set. Clearly, $K\subseteq U$ because $TU\subseteq U$.
\end{proof}

\section{Permanence criteria}\label{sec:permanence}
From now on we reserve the symbol $n$ for the dimension of the Euclidean space $\mathbb{R}^n$ and the symbol $N$ for the set $\{1, \ldots, n\}$. We will denote by $\{\mathbf{e}_{1}, \ldots, \mathbf{e}_{n}\}$ the usual basis for $\mathbb{R}^n$, and by $d(\cdot,\cdot)$ the usual Euclidean distance. We use $\mathbb{R}^n_+$ to denote the nonnegative cone $\{x\in \mathbb{R}^n: x_i\geq 0,\forall i\in N\}$. The interior of $\mathbb{R}^n_+$ is the open cone $\dot{\mathbb{R}}^n_+:= \{x\in \mathbb{R}^n_+: x_i>0,\forall i\in N \}$ and the boundary of $\mathbb{R}^n_+$ is $\partial \mathbb{R}^n_+:=\mathbb{R}^n_+\setminus \dot{\mathbb{R}}^n_+$.  The symbol $0$  stands for both the origin of $\mathbb{R}^n$ and the real number $0$.

Given two points $x,z$ in $\mathbb{R}^n$, we write $x \leq z$ if $z - x \in \mathbb{R}^n_+$, $x < z$ if $z - x \in \mathbb{R}^n_+\setminus \{0\}$, and $x \ll z $ if $z-x\in \dot{\mathbb{R}}^n_+$.  The reverse relations are denoted by $\geq, >,\gg$, respectively.

Given an $m\times m$ matrix $A$, we write $A\geq 0$ if $A$ is a nonnegative matrix (i.e., all its entries are nonnegative) and $A>0$ if $A$ is a positive matrix (i.e., all its entries are positive). We shall use $I$ to denote the identity matrix.

Consider the map $T:\mathbb{R}_+^n\to \mathbb{R}_+^n$ given by
\begin{equation}\label{map_T_dimn}
T(x)=(x_1F_1(x),\ldots,x_nF_n(x))
\end{equation}
with continuous functions $F_i$ satisfying $F_i(x)>0$ for all $x\in \mathbb{R}_+^n$. Note that this implies that $T_i(x)>0$ if and only if $x_i>0$, and hence $T\dot{\mathbb{R}}^n_+\subset \dot{\mathbb{R}}^n_+$ and $T(\partial \mathbb{R}^n_+)\subset \partial \mathbb{R}^n_+$. In particular, $T^{-1}(\{0\})=\{0\}$.
\begin{definition}[\cite{hofbauer1987coexistence,Hutson1982}]
The map $T$ is said to be permanent if there exists a compact positively invariant set $K\subseteq \dot{\mathbb{R}}^n_+$ such that for every $x\in \dot{\mathbb{R}}^n_+$, there exists a $m=m(x)>0$ such that the tail $\gamma^+_m(x)\subseteq K$.
\end{definition}
\begin{remark}
Since $K\subseteq \dot{\mathbb{R}}^n_+$ is compact, the distance $d(K,\partial \mathbb{R}^n_+)$ of $K$ from the boundary $\partial \mathbb{R}^n_+$ is thus non-zero, and $\omega(x)\subseteq K$ for all $x\in \dot{\mathbb{R}}^n_+$. Equivalently, permanence means that there exist $\delta,D>0$ such that
\begin{equation}
\delta\leq \liminf_{k\to +\infty} T_i^kx\leq \limsup_{k\to +\infty} T_i^kx \leq D,\quad i=1,\ldots,n,
\end{equation}
for all $x\in \dot{\mathbb{R}}^n_+$ {\rm(\cite{hofbauer1988,Kon2004,Smith2010})}. In a permanent system, species can coexist permanently in the sense that when the population densities of all species are positive, after some generations each population density will be bounded away from zero and infinity for all the time. As a consequence extinction and explosion cannot occur. In this paper, the map $T$ is said to be impermanent if it is not permanent.
\end{remark}
\subsection{Permanence for dissipative systems}
Our first criteria on permanence and impermanence for dissipative systems defined on $\mathbb{R}^n_+$ is the following theorem based on the technique of average Liapunov functions. A similar result was given by Garay and Hofbauer for discrete-time replicator dynamics in \cite{Garay2003}. See also \cite{Smith2010} for a detailed discussion on the technique of average Liapunov functions.
\begin{theorem}\label{theorem:permanence}
Suppose that $T$ is dissipative with global attractor $J$. If there are real numbers $\nu_1,\ldots,\nu_n>0$ such that
\begin{equation}\label{con-01}
    g(x)=\sum^n_{i=1}\nu_i\ln F_i(x)>0,\quad\forall x\in \Omega(\partial J),
\end{equation}
where $\partial J=J\cap \partial\mathbb{R}_+^n$, then $T$ is permanent; if instead
\begin{equation}\label{con-02}
    g(x)=\sum^n_{i=1}\nu_i\ln F_i(x)<0,\quad\forall x\in \Omega(\partial J),
\end{equation}
then $T$ is impermanent.
\end{theorem}
\begin{proof}
We first show that if \eqref{con-01} holds then $T$ is permanent. Let
\begin{equation}\label{fun-Vx}
V(x)=x^{\nu_1}_1x^{\nu_2}_2\cdots x^{\nu_n}_n, \quad x\in \mathbb{R}^n_+.
\end{equation}
Note that $V(x)=0$ if and only if $x\in \partial \mathbb{R}^n_+$. For $x\in \dot{\mathbb{R}}^n_+$, we let $\theta(k,x)=V(T^kx)/V(x)$, $k\in \mathbf{Z}_+$.
Then
\begin{equation}\label{theta1}
\theta(1,x)=\frac{V(Tx)}{V(x)}=F_1^{\nu_1}(x)\cdots F_n^{\nu_n}(x), \quad x\in \dot{\mathbb{R}}^n_+.
\end{equation}
Since $F_i:\mathbb{R}^n_+\mapsto \mathbb{R}_+\setminus \{0\}$ are continuous, \eqref{theta1}
provides a continuous extension of $\theta(1,\cdot)$ to $\mathbb{R}^n_+$. Thus, we have $\theta(1,x)=F_1^{\nu_1}(x)\cdots F_n^{\nu_n}(x)>0$, $\forall x \in \mathbb{R}^n_+$. Note that for any $x\in \dot{\mathbb{R}}^n_+$ and $k\geq 2$,
\begin{equation}\label{thetak}
\begin{array}{rl}
\theta(k,x)=&\displaystyle\frac{V(T^kx)}{V(T^{k-1}x)}\frac{V(T^{k-1}x)}{V(T^{k-2}x)}\cdots\frac{V(Tx)}{V(x)}\\
\noalign{\medskip}
=&\displaystyle \theta(1,T^{k-1}x)\theta(1,T^{k-2}x)\cdots \theta(1,x).
\end{array}
\end{equation}
So \eqref{thetak}
provides a continuous extension of $\theta(k,\cdot)$ to $\mathbb{R}^n_+$, $\forall k\geq 2$.

Consider the maps $T|_J:J\mapsto J$ and $V|_J:J\mapsto \mathbb{R}_+$ given by \eqref{fun-Vx}, where $T|_J$ denotes the restriction of $T$ on $J$, and similarly for $V|_J$. Note that $V|_J(x)=0$ if and only if $x\in \partial J$. For $x\in \partial J$, one has
\begin{equation}
\vartheta_J(k,x):=\displaystyle\liminf_{\mbox{\tiny$\begin{array}{c}
y\rightarrow x\\
y\in J\setminus \partial J
\end{array}$}}\frac{V|_J(T^ky)}{V|_J(y)}
=\displaystyle\liminf_{\mbox{\tiny$\begin{array}{c}
y\rightarrow x\\
y\in J\setminus \partial J
\end{array}$}} \theta(k,y)=\theta(k,x).
\end{equation}
Let $\psi_J(x):=\sup_{k>0}\vartheta_J(k,x)$ for $x\in \partial J$.
By condition \eqref{con-01}, $\exp\{g(x)\}>1$ for all $x\in \Omega(\partial J)$, that is
\begin{equation}
\vartheta_J(1,x)=F^{\nu_1}_1(x)\cdots F^{\nu_n}_n(x)>1, \quad\forall x\in\Omega(\partial J).
\end{equation}
Therefore, $\psi_J(x)>1$ for all $x\in \Omega(\partial J)$, and obviously $\psi_J(x)\geq \vartheta_J(1,x) >0$ for all $x\in \partial J$. In Lemma \ref{repellor_T} take $M=J$, $S=\partial J$ and $T=T|_J$. Then by Remark \ref{con:repel}, we have
$\psi_J(x)>1$ for all $x\in \partial J$.

Let $W=O_\varepsilon(J)$ be an $\varepsilon$-neighborhood of $J$ in $\mathbb{R}^n_+$. Since $J$ is the global attractor, for all $x\in \mathbb{R}^n_+$, $\gamma^+(x)\cap W\neq \emptyset$. It then follows from Lemma \ref{lemma:absorbing} that $\gamma^+(\overline{W})$ is compact and positively invariant such that for any $x\in \mathbb{R}^n_+$, there exists a $k=k(x)>0$ satisfying $\gamma^+_k(x)\subseteq \gamma^+(\overline{W})$.

Now take $M=\gamma^+(\overline{W})$ and $S=M\cap \partial \mathbb{R}^n_+$, which are compact. Consider the maps $T|_M:M\mapsto M$ and $V|_M:M\mapsto \mathbb{R}_+$ given by \eqref{fun-Vx}. Obviously, $S$ and $M\setminus S$ are positively invariant under $T|_M$.  We show that $S$ repels for $T|_M$.

Note that $J\subseteq M, \partial J\subseteq S$ and $V|_M(x)=0$ if and only if $x\in S$.  For $x\in S$, one has
\begin{equation}
\vartheta_M(k,x):=\displaystyle\liminf_{\mbox{\tiny$\begin{array}{c}
y\rightarrow x\\
y\in M\setminus S
\end{array}$}}\frac{V|_M(T^ky)}{V|_M(y)}
=\displaystyle\liminf_{\mbox{\tiny$\begin{array}{c}
y\rightarrow x\\
y\in M\setminus S
\end{array}$}} \theta(k,y)=\theta(k,x).
\end{equation}
Let $\psi_M(x):=\sup_{k>0}\vartheta_M(k,x)$ for $x\in S$. Note that for all $x\in \partial J$,
$$
\vartheta_M(k,x)=\theta(k,x)=\vartheta_J(k,x),
$$
so we get $\psi_M(x)=\psi_J(x)>1$ for all $x\in \partial J$. Since  $\omega(x)\subseteq J$ for all $x\in \mathbb{R}^n_+$, one has
$$\omega(x)\subseteq J\cap \partial \mathbb{R}^n_+=\partial J,\quad \forall x\in S.$$
Thus, $\Omega(S)\subseteq \partial J$, which implies that $\psi_M(x)>1$ for all $x\in \Omega(S)$. Clearly, $\psi_M(x)\geq \vartheta_M(1,x) >0$ for all $x\in S$. It then follows from Remark \ref{con:repel} that
$\psi_M(x)>1$ for all $x\in S$, and hence $S$ repels for $T|_M$ by Lemma \ref{repellor_T}. Therefore, there exists a compact set $K\subseteq M\setminus S$ of the subspace $M$ which is positively invariant under $T|_M$, such that for every $x\in M\setminus S$, there exists a $m=m(x)>0$ such that $\gamma^+_m(x)\subseteq K$ by Proposition \ref{prop:repel}. Of course, $K$ is also a compact subset of $\mathbb{R}^n_+$ and positively invariant under $T$. Note that $K\subseteq M\setminus S \subseteq \dot{\mathbb{R}}^n_+$.

Recall that for any $x\in \dot{\mathbb{R}}^n_+$, there exists a $k=k(x)>0$ such that $T^k(x) \in M$, and hence $T^k(x) \in M\setminus S$ because $T \dot{\mathbb{R}}^n_+ \subseteq \dot{\mathbb{R}}^n_+$. Therefore, there exists a $m=m(x)\geq k$ such that the tail $\gamma^+_m(x)\subseteq K$, that is $T$ is permanent.

Now suppose that \eqref{con-02} holds. Consider the maps $T|_J:J\mapsto J$ and $V|_J:J\mapsto \mathbb{R}_+$ given by \eqref{fun-Vx}. In Lemma \ref{attractor_T} take $M=J,S=\partial J$ and $T=T_J$. Since $\theta(1,\cdot):\mathbb{R}^n_+ \mapsto \mathbb{R}_+$ is continuous, there exist a constant $C>0$ such that $\frac{V|_J(Tx)}{V|_J(x)}\leq C$ for all $x\in J\setminus \partial J$.  It follows from \eqref{con-02} that $\exp\{g(x)\}<1$ and hence $\theta(1,x)<1$, for all $x\in \Omega(\partial J)$. Then for all $x\in \Omega(\partial J)$, one has
$\varphi_J(x)\leq \zeta_J(1,x)=\theta(1,x)<1$,
where $\varphi_J(x)$ and $\zeta_J(1,x)$ are defined in Lemma \ref{attractor_T}. Thus, $\partial J$ attracts for $T_J$. Therefore, there exists some $x\in \dot{\mathbb{R}}^n_+$ such that $\omega(x)\subseteq \partial J\subseteq \partial \mathbb{R}^n_+$, and hence $T$ is impermanent.
\end{proof}
\subsection{Permanence via carrying simplex}\label{subsection:permanence}
 Before presenting the permanence and impermanence criteria for the map $T$ given by \eqref{map_T_dimn} admitting a carrying simplex $\Sigma$, we first recall the properties of carrying simplex.

A {\it carrying simplex} for the map $T$ is a subset $\Sigma$ of $\mathbb{R}_+^n\setminus \{0\}$ with the following properties:
\begin{enumerate}[(P1)]
\item $\Sigma$ is compact and invariant under $T$;
\item for any $x\in \mathbb{R}_+^n\setminus \{0\}$, there exists some $z\in \Sigma$ such that
$\displaystyle\lim_{k\to \infty}|T^kx-T^kz|=0$;
\item $\Sigma$ is unordered (i.e. if $x,z\in \Sigma$ such that $x_i\geq z_i$ for all $i\in N$, then $x = z$), and homeomorphic to the probability simplex $\Delta^{n-1}$ via radial projection;
\item $T:\Sigma\mapsto \Sigma$ is a homeomorphism.
\end{enumerate}

(P1) and (P2) imply that the long-term dynamics of $T$ is accurately reflected by that in $\Sigma$, and (P3) means that $\Sigma$ is topologically simple. We denote the boundary of $\Sigma$, i.e. $\Sigma\cap \partial\mathbb{R}^n_+$ by $\partial \Sigma$, and the interior of $\Sigma$, i.e. $\Sigma \setminus \partial \Sigma$ by $\dot{\Sigma}$.

We denote by $\mathbb{H}^+_{i}$ the $i$-th positive coordinate axis and by $\pi_i=\{x\in \mathbb{R}^n_+:x_i=0\}$ the $i$-th coordinate plane. Note that each $\pi_i$ is positively invariant under $T$ and $\partial\Sigma\cap \pi_i$ is the carrying simplex of $T|_{\pi_i}$, that is $\partial\Sigma$ is composed of the carrying simplices of $T|_{\pi_i}$, $i=1,\ldots,n$. $\Sigma$ contains all non-trivial fixed points, periodic orbits and heteroclinic cycles, etc. Every vertex of $\Sigma$ is a fixed point of $T$, where $\Sigma$ and some positive coordinate axis meet, and denote by $q_{\{i\}}=q_i\mathbf{e}_{i}$ the fixed point at the vertex where $\Sigma$ and $\mathbb{H}_{i}^+$ meet. For one-dimensional case, $T$ admits a carrying simplex if and only if it has a globally attracting positive fixed point in $\mathbb{R}_+\setminus \{0\}$.

A map $T:\mathbb{R}^n_+\to \mathbb{R}^n_+$ is competitive (or retrotone) in a subset $W\subset \mathbb{R}^n_+$ if for all $x,z\in W$ with $Tx < Tz$ one has that $x_i<z_i$ provided $z_i>0$.

We first recall a readily checked criterion provided by Jiang and Niu \cite{LG} on the existence of a carrying simplex for the competitive map $T$ of type \eqref{map_T_dimn}.

\begin{lemma}[Existence Criterion of Carrying Simplex  \cite{LG}]\label{simplex}
Suppose $F_i$ are $C^1$, $i=1,\ldots,n$. Assume that
\begin{itemize}
  \item[{\rm $\Upsilon$1)}] $\partial F_i(x)/\partial x_j<0$ holds for any $x\in \mathbb{R}^n_+$ and $i,j\in N;$
  \item[{\rm $\Upsilon$2)}] $\forall i \in N$, $T|_{\mathbb{H}_{i}^+}: \mathbb{H}_{i}^+\to \mathbb{H}_{i}^+$ has a fixed point $q_{\{i\}}=q_i\mathbf{e}_{i}$ with $q_i>0;$
  \item[{\rm $\Upsilon$3)}] $\forall x\in [0,q]\setminus \{0\}$, $F_i(x)+\sum_{j\in \kappa(x)}x_j\frac{\partial F_i(x)}{\partial x_j}>0$ holds for any $i\in \kappa(x)$
  $($or $F_i(x)+\sum_{j\in \kappa(x)}x_i\frac{\partial F_i(x)}{\partial x_j}>0$ holds for any $i\in \kappa(x))$,
where $q=(q_1,\ldots,q_n)$ and $\kappa(x)=\{i:x_i>0\}$ is the support of $x$.
\end{itemize}
Then $T$ possesses a carrying simplex $\Sigma\subset [0,q]$.
\end{lemma}

Condition $\Upsilon$1) means that
$F_i(y)<F_i(x)$ for all $i\in N$ provided $x<y$. This follows from
$$
   F_i(y)-F_i(x)=\int_0^1DF_i(x_s)(y-x)ds,
$$
where $x_s=x+s(y-x)$ with $s\in[0,1]$. Together with $\Upsilon$2), $\Upsilon$1) implies $F_i(0)>F_i(q_{\{i\}})=1$ for all $i\in N$, i.e. $0$ is a hyperbolic repeller for $T$. $\Upsilon$3) implies that $\det DT(x)>0$ for all $x\in [0,q]$, and together with $\Upsilon$1) it guarantees $(DT(x)_{\kappa(x)})^{-1}> 0$ for all $x\in [0,q]\setminus \{0\}$ (see \cite[Theorem 3.1]{LG}), so $T$ is competitive and one-to-one in $[0,q]$ by \cite[Proposition 4.1]{ruiz2011exclusion}.

\begin{theorem}\label{Sigma:permanence}
Assume that $T$ admits a carrying simplex $\Sigma$. If there are real numbers $\nu_1,\ldots,\nu_n>0$ such that
\begin{equation}\label{Sigma:con-01}
    g(x)=\sum^n_{i=1}\nu_i\ln F_i(x)>0\quad\forall x\in \Omega(\partial \Sigma),
\end{equation}
then $T$ is permanent; if instead
\begin{equation}\label{Sigma:con-02}
    g(x)=\sum^n_{i=1}\nu_i\ln F_i(x)<0\quad\forall x\in \Omega(\partial \Sigma),
\end{equation}
then $T$ is impermanent.
\end{theorem}
\begin{proof}
We first show that if \eqref{Sigma:con-01} holds then $T$ is permanent. Since $\Sigma$ is invariant and compact such that $\omega(x)\subseteq \Sigma$ for all $x\in \mathbb{R}^n_+\setminus \{0\}$, there exists an $\varepsilon$-neighborhood $O_\varepsilon(\Sigma)\subset \mathbb{R}^n_+\setminus \{0\}$ of $\Sigma$ such that  $\overline{O_\varepsilon(\Sigma)}\subseteq \mathbb{R}^n_+\setminus \{0\}$ and $\gamma^+(x)\cap O_\varepsilon(\Sigma)\neq \emptyset$ for all $x\in \mathbb{R}^n_+\setminus \{0\}$. Then it follows from Lemma \ref{lemma:absorbing} that $\gamma^+(\overline{O_\varepsilon(\Sigma)})$ is a compact positively invariant set. Clearly, $M=\gamma^+(\overline{O_\varepsilon(\Sigma)})$ is a compact neighborhood of $\Sigma$. Set $S=M\cap \partial \mathbb{R}^n_+$. By the property (P2) of $\Sigma$, one has $\omega(x)\subseteq \Omega(\partial \Sigma)$ for any $x\in S$, and hence $\Omega(S)\subseteq \Omega(\partial \Sigma)$. So, if \eqref{Sigma:con-01} holds, then $g(x)>0$ for all $x\in \Omega(S)$.

Now consider the map $T|_M:M\mapsto M$. Let $V(x)=x^{\nu_1}_1\cdots x^{\nu_n}_n$, $x\in M$. Note that $V(x)=0$ if and only if $x\in S$.
By the above analysis, we know $\exp\{g(x)\}>1$ for all $x\in \Omega(S)$, that is, $F^{\nu_1}_1(x)\cdots F^{\nu_n}_n(x)>1$ for all $x\in \Omega(S)$. On the other hand, for $x\in S$, one has
$$
\begin{array}{rl}
 \vartheta_M(1,x)=&\displaystyle\liminf_{\mbox{\tiny$\begin{array}{c}
y\rightarrow x\\
y\in M\setminus S
\end{array}$}}\frac{V(Ty)}{V(y)} \\
\noalign{\medskip}
 =&\displaystyle\liminf_{\mbox{\tiny$\begin{array}{c}
y\rightarrow x\\
y\in M\setminus S
\end{array}$}} F^{\nu_1}_1(y)\cdots F^{\nu_n}_n(y)\\
\noalign{\medskip}
 =&F^{\nu_1}_1(x)\cdots F^{\nu_n}_n(x).
\end{array}
$$
Thus, $\vartheta_M(1,x)>1$ for all $x\in \Omega(S)$, and $\vartheta_M(1,x)>0$ for all $x\in S$. An application of Lemma \ref{repellor_T} and Remark \ref{con:repel} to such $M$, $S$ and $T=T|_M$ shows that $S$ repels for $T_M$. Then the rest of the proof can be completed by repeating the same arguments as in Theorem \ref{theorem:permanence}.

If \eqref{Sigma:con-02} holds, then $\partial \Sigma$ attracts for $T|_{\Sigma}$ by \cite[Theorem 2]{jiang2015}, and hence there exists some $x\in \dot{\mathbb{R}}^n_+$ such that $\omega(x)\subseteq \partial \Sigma \subseteq \partial \mathbb{R}^n_+$, which implies that $T$ is impermanent.
\end{proof}
\begin{remark}
Note that in the proof of Theorem \ref{Sigma:permanence}, we do not need the properties {\rm(P3)} and {\rm(P4)} of the carrying simplex. In fact, the results in Theorem \ref{Sigma:permanence} hold for other kinds of maps which have an attracting and invariant manifold $\mathcal{S}$, that is $\omega(x)\subseteq \mathcal{S}$ for all $x\in \mathbb{R}^n_+\setminus \{0\}$, although we mainly focus on the maps with a carrying simplex here.
\end{remark}

For the two-dimensinal (i.e. $n=2$) map $T$ given by \eqref{map_T_dimn} with a carrying simplex $\Sigma$, we know that $T|_\Sigma$ is topologically conjugate to a strictly increasing 
homeomorphism $h$ taking $[0,1]$ onto $[0,1]$ (similarly for $(T|_\Sigma)^{-1}$) because $\Sigma$ is homeomorphic to $[0,1]$ and $T|_\Sigma$ is a homeomorphism taking $\Sigma$ onto $\Sigma$, and hence every nontrivial orbit of $T$ converges to some fixed point on $\Sigma$ (see also \cite{smith1986,smith1998planar}). In particular, it follows from Theorem \ref{Sigma:permanence} that $T$ is permanent if $F_i(q_{\{j\}})>1$ for $i\neq j$, $i,j=1,2$, where $q_{\{j\}}$ is the fixed point on $\mathbb{H}^+_{j}$.

%
%

By the above arguments, we know that for the two-dimensional map $T$ which admits a carrying simplex, the dynamics is relatively simple, that is $\omega(x)\subseteq \mathcal{E}(T)$ for all $x\in \mathbb{R}^2_+$, and hence $\Omega(\Sigma)\subseteq \mathcal{E}(T)$. For the three-dimensional case, since $\partial\Sigma\cap \pi_i$ is the carrying simplex of $T|_{\pi_i}$, which is a two-dimensional map, the boundary dynamics for $T$ is simple, i.e. $\Omega(\partial\Sigma)\subset \mathcal{E}(T)$.

%

\begin{corollary}\label{coro:3-permanence}
Let $n=3$. Suppose that $T(x)=(x_1F_1(x),x_2F_2(x),x_3F_3(x))$ taking $\mathbb{R}^3_+$ into $\mathbb{R}^3_+$ admits a carrying simplex $\Sigma$. If there are real numbers $\nu_1,\nu_2,\nu_3>0$ such that
\begin{equation}\label{Sigma-con-02}
    g(\hat{x})=\sum^3_{i=1}\nu_i\ln F_i(\hat{x})>0\,(resp. <0),\quad\forall \hat{x}\in \mathcal{E}(T)\cap \partial \Sigma,
\end{equation}
then $T$ is permanent \,(resp. impermanent).
\end{corollary}
\begin{proof}
The conclusion follows from the above analysis and Theorem \ref{Sigma:permanence} immediately.
\end{proof}

For low-dimensional systems, one remarkable phenomenon is the occurrence of heteroclinic cycles, i.e., the cyclic arrangements of saddle fixed points and heteroclinic connections; see \cite{roeger2004discrete,Lih2005,Diekmann2005,cushing2015,jiang2015}.

\begin{figure}[h]
 \begin{center}
 \includegraphics[width=0.62\textwidth]{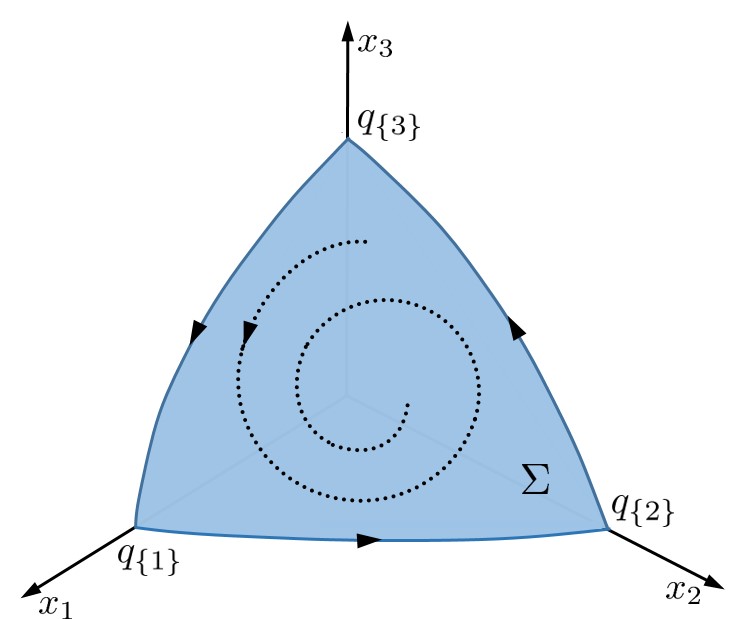}
\caption{A carrying simplex $\Sigma$ with a repelling heteroclinic cycle $\partial\Sigma$.} \label{fig:cs}
  \end{center}
\end{figure}

Let $n=3$. Suppose that $F_i$ are $C^1$ and $T$ admits a carrying simplex $\Sigma$ (homeomorphic to $\Delta^2$) with three axial fixed points $q_{\{1\}}=(q_1,0,0)$, $q_{\{2\}}=(0,q_2,0)$ and $q_{\{3\}}=(0,0,q_3)$, which lie at the vertices of $\Sigma$. Assume that $q_{\{1\}},q_{\{2\}},q_{\{3\}}$ are saddles on $\Sigma$, and $\partial \Sigma\cap \pi_i$ is the heteroclinic connection between $q_{\{j\}}$ and $q_{\{k\}}$. In this case, there are no other fixed points on $\partial \Sigma$ which is a heteroclinic cycle of May-Leonard type: $q_{\{1\}}\to q_{\{2\}} \to q_{\{3\}}\to q_{\{1\}}$ (or the arrows reversed); see Fig. \ref{fig:cs}. For more details, see \cite{Hofbauer1994,hofbauer1988,jiang2015}.

\begin{lemma}[Theorem 3 in \cite{jiang2015}]\label{heteroclinic-cycle}
Suppose $F_i$ are $C^1$, $i=1,2,3$. Assume that $T$ admits a carrying simplex $\Sigma$, and $\partial \Sigma$ is a heteroclinic cycle above. Then the heteroclinic cycle $\partial \Sigma$ repels (resp. attracts), if
\begin{equation}\label{heteroclinic-con}
    \varrho:=\prod_{i=1}^3 \ln F_i(q_{\{i-1\}})+\prod_{i=1}^3\ln F_i(q_{\{{i+1}\}})>0~(resp.~<0),
\end{equation}
 where $i\in \{1,2,3\}$ is considered cyclic.
\end{lemma}
\begin{corollary}\label{stability-permance}
Suppose $F_i$ are $C^1$, $i=1,2,3$. Assume that $T$ admits a carrying simplex $\Sigma$, and $\partial \Sigma$ is a heteroclinic cycle above. If $\varrho>0$ (resp. $<0$), i.e. $\partial \Sigma$ repels (resp. attracts), where $\varrho$ is defined by \eqref{heteroclinic-con}, then $T$ is permanent (resp. impermanent).
\end{corollary}
\begin{proof}
Under the assumption, one has $\Omega(\partial \Sigma)=\{q_{\{1\}},q_{\{2\}},q_{\{3\}}\}$. It follows from the proof of Theorem 3 in \cite{jiang2015} that there are real numbers $\nu_1,\nu_2,\nu_3>0$ such that \eqref{Sigma:con-01} holds if $\varrho>0$,  so $T$ is permanent by Theorem \ref{Sigma:permanence}; see Fig. \ref{fig:cs}. If $\varrho<0$, then $\partial \Sigma$ attracts, and hence $T$ is impermanent.
\end{proof}

\section{Extensions to competitive systems}
\label{sec:extensions}
In this section, we study the Kolmogorov map $T$ given by \eqref{equ:T1}.

For the convenience of the study, we set $a_{ij}=r_i\mu_{ij}$. Let $R=\diag[r_i]$, the ${n\times n}$ diagonal matrix with diagonal entries $r_i$, $i=1,\ldots,n$, and $U$ be the ${n\times n}$ matrix with entries $\mu_{ij}>0$. Then $A=RU$, and \eqref{equ:T1} is written as
\begin{equation}\label{equ:T2}
    T_i(x)=x_iF_i(x)=x_if_i((RUx^\tau)_i,r_i)=x_if_i(r_i\sum_{j=1}^n\mu_{ij}x_j,r_i),\quad i=1,\cdots,n.
\end{equation}
In this form, the fixed points of $T$ are determined by the linear equations
$$
   x_i=0~~\textrm{or}
    ~~\sum_{j=1}^n\mu_{ij}x_j=1,\quad i=1,\ldots,n,
$$
which depend only on the parameters $\mu_{ij}$.

Denote by $\mathscr{F}$ the collection of all $C^1$ functions $f: \mathbb{R}_+\times \dot{\mathbb{R}}_+\to \dot{\mathbb{R}}_+$ satisfying \eqref{cons:f}.  Let $\mathcal{F}_n=\{f=(f_1,\ldots,f_n): f_i\in \mathscr{F},\ i=1,\ldots,n\}$. Given $f\in \mathcal{F}_n$, denote by
$${\rm DCS}(n,f)=\{T\in \mathcal{T}(\mathbb{R}_+^n):T_i(x)=x_i\displaystyle f_i(r_i\sum_{j=1}^n\mu_{ij}x_j,r_i),\ \mu_{ij}>0,r_i>0\}$$
the set of all maps on $\mathbb{R}_+^n$ of the from \eqref{equ:T2} with the given function $f\in \mathcal{F}_n$, where $\mathcal{T}(\mathbb{R}_+^n)$ denotes the set of all maps taking $\mathbb{R}_+^n$ into itself. $f\in \mathcal{F}_n$ is called a {\it generating function} for the map \eqref{equ:T2}. For $T\in\mathrm{DCS}(n,f)$, we always let $F_i(x)=f_i((RUx^\tau)_i,r_i)$ such that $T_i(x)=x_iF_i(x)$, $i=1\ldots,n$.


Let $T \in \mathrm{DCS}(n,f)$. The entries of the
Jacobian matrix $DT(x)$ at $x$ are given by
\begin{equation}\label{jacobian-x}
\begin{array}{rl}
(DT(x))_{ij}=&\left\{
\begin{array}{ll}
\delta_{ij}F_i(x)+x_i\frac{\partial f_i}{\partial z}((RUx^\tau)_i,r_i)r_i\mu_{ij},&\quad i\in \kappa(x),\\
\noalign{\medskip}
F_j(x),&\quad i\notin \kappa(x),\ j=i,\\
\noalign{\medskip}
0, &\quad i\notin \kappa(x),\ j\neq i,
\end{array}\right.
\end{array}
\end{equation}
where $\delta_{ij}=1$ for $i=j$ and $\delta_{ij}=0$ for $i\neq j$, $i,j=1,\ldots,n$.
$DT(x)$ ``splits" into two blocks: the square matrix $DT(x)_{\kappa(x)}$ defines the ``internal" block which corresponds to the Jacobian matrix of the restriction of $T$ to the subspace $\mathbf{R}_+^{\kappa(x)}=\{x\in \mathbb{R}_+^n:x_i=0~\mathrm{for}~i\notin \kappa(x)\}$, where $DT(x)_{\kappa(x)}$ is the submatrix of $DT(x)$ with rows and columns from $\kappa(x)$; the square matrix $DT(x)_{N\setminus \kappa(x)}$ is the ``external" block  which is a diagonal matrix with diagonal entries $F_j(x)>0$, where $j\notin \kappa(x)$ and $DT(x)_{N\setminus \kappa(x)}$ is the submatrix of $DT(x)$ with rows and columns from $N\setminus \kappa(x)$.

Let $\hat{x}=(\hat{x}_1,\ldots,\hat{x}_n)$ be a fixed point of $T$. Then $F_i(\hat{x})=1$ for any $i\in \kappa(\hat{x})$, and
\begin{equation}\label{jacobian-T}
\begin{array}{rl}
(DT(\hat{x}))_{ij}=&\left\{
\begin{array}{ll}
\delta_{ij}+\hat{x}_i\frac{\partial f_i}{\partial z}(r_i,r_i)r_i\mu_{ij},&\quad i\in \kappa(\hat{x}),\\
\noalign{\medskip}
F_j(\hat{x})=f_j((RU\hat{x}^\tau)_j,r_j),&\quad i\notin \kappa(\hat{x}),\ j=i,\\
\noalign{\medskip}
0, &\quad i\notin \kappa(\hat{x}),\ j\neq i.
\end{array}\right.
\end{array}
\end{equation}
The external block $DT(\hat{x})_{N\setminus \kappa(x)}$ of $DT(\hat{x})$ is a diagonal matrix whose entries are the external eigenvalues $F_j(\hat{x})>0$ (we call it the external eigenvalue in direction $j$), where $j\notin \kappa(\hat{x})$.

The eigenvalues of $DT(0)$ are $f_i(0,r_i)>1$, i.e. $F_i(0)>1$,  $i=1,\ldots,n$, so the trivial fixed point $0$ is a hyperbolic repeller.

\begin{lemma}[Gerschgorin Circle Theorem \cite{Meyer}]\label{circle-theorem}
Let $\mathcal{B}$ be an $n\times n$ matrix with entries $b_{ij}$. Define the $i$-{\rm th} Gerschgorin disc $\mathcal{D}_i$ in the complex plane to be the closed disc centered at $b_{ii}$ with radii $\sum_{j\neq i}|b_{ij}|$. Each $\mathcal{D}_i$ contains an eigenvalue of $\mathcal{B}$ and, moreover, for any distinct $i_1,\cdots,i_m$, there are at least $m$ eigenvalues $($counting multiplicities$)$ of $\mathcal{B}$ in $\bigcup_{k=1}^m \mathcal{D}_{i_k}$.
\end{lemma}
\begin{lemma}[Lemma 2.3.4 in \cite{C82}]\label{one-to-one}
Suppose that $M\subset \mathbb{R}^n$ is a connected compact set and the continuous function $g:M\mapsto g(M)$ is a local homeomorphism. Then the cardinal number of $g^{-1}(\{z\}) $ is finite and constant for all $z\in g(M)$.
\end{lemma}
\begin{proposition}\label{prop-competitive}
Every $T\in \mathrm{DCS}(n,f)$ is a diffeomorphism from $\mathbb{R}_+^n$ to its image and also a competitive map on $\mathbb{R}_+^n$.
\end{proposition}
\begin{proof}
We first show that $T$ is a local diffeomorphism. According to the inverse function theorem, it suffices to prove that $\det DT(x)>0$ for all $x\in \mathbb{R}_+^n$. Recall that $DT(x)$ splits into two blocks:  the internal block $DT(x)_{\kappa(x)}$ and the external block $DT(x)_{N\setminus \kappa(x)}$ which is a diagonal matrix with positive diagonal entries, so we only need to show that $\det DT(x)_{\kappa(x)}>0$. Therefore, without loss of generality, we assume that $x\in \dot{\mathbb{R}}_+^n$, and show $\det DT(x)>0$. By \eqref{jacobian-x}, $DT(x)$ can be written as $DT(x)=\diag[F_k(x)]+\diag[x_k] \mathcal{W}$, where $\mathcal{W}$ is the matrix whose $(i,j)$-th entry is given by $\frac{\partial F_i(x)}{\partial x_j}<0$. Note that $\diag[x_k]$ is invertible because $x_k>0$. Therefore, $DT(x)$ is similar to
$$
\mathcal{B}:=\diag[x_k]^{-1}DT(x)\diag[x_k]=\diag[F_k(x)]+\mathcal{W}\diag[x_k].
$$
Note that the $(i,j)$-th entry of $\mathcal{W}\diag[x_k]$ is given by
\begin{equation}\label{partial-xj}
x_j\frac{\partial F_i(x)}{\partial x_j}=r_i\mu_{ij}x_j\frac{\partial f_i}{\partial z}((RUx^\tau)_i,r_i)<0.
\end{equation}
Then by \eqref{cons:f} (ii) one has
\begin{equation}\label{condition-gamma-3}
\begin{array}{rl}
  &\displaystyle F_i(x)+\sum_{j=1}^nx_j\frac{\partial F_i(x)}{\partial x_j} \\
  \noalign{\medskip}
 =&\displaystyle F_i(x)+\sum_{j=1}^nr_i\mu_{ij}x_j\frac{\partial f_i}{\partial z}((RUx^\tau)_i,r_i) \\
 \noalign{\medskip}
 =&\displaystyle f_i((RUx^\tau)_i,r_i)+(RUx^\tau)_i\frac{\partial f_i}{\partial z}((RUx^\tau)_i,r_i)>0.
\end{array}
\end{equation}
It follows from \eqref{partial-xj} and \eqref{condition-gamma-3} that the diagonal entries $F_i(x)+x_i\frac{\partial F_i(x)}{\partial x_i}$ of $\mathcal{B}$ are positive and, moreover, each Gerschgorin disc $\mathcal{D}_i$ of $\mathcal{B}$ which is centered at $F_i(x)+x_i\frac{\partial F_i(x)}{\partial x_i}$ with radii $-\sum_{j\neq i}x_j\frac{\partial F_i(x)}{\partial x_j}$ lies in the right half-plane. Then by Lemma \ref{circle-theorem} all the eigenvalues of $\mathcal{B}$ have positive real parts, and hence
$\det DT(x)= \det \mathcal{B}>0$. At this moment we have proved that $T$ is a local diffeomorphism.

Now we show that $T$ is one-to-one. By a contradiction argument assume that there exist $x\neq y$ such that $Tx=Ty$. Then one can choose some $l>0$ such that $0,x,y\in \overline{B}_{l}$, where
$$
\overline{B}_{l}=\{z\in \mathbb{R}_+^n:|z|\leq l\}.
$$
Consider the restriction $$T|_{\overline{B}_{l}}: \overline{B}_{l} \mapsto T\overline{B}_{l}. $$
It follows from Lemma \ref{one-to-one} that $T|_{\overline{B}_{l}}^{-1}(\{z\})$ is finite and constant for all $z\in T\overline{B}_{l}$.
Since, as noticed above, $T|_{\overline{B}_{l}}^{-1}(\{0\})=\{0\}$, this constant is one and hence $T|_{\overline{B}_{l}}$ is one-to-one, contradicting that $T|_{\overline{B}_{l}}(x)=T|_{\overline{B}_{l}}(y)$. Thus, we have proved that
$T$ is a diffeomorphism.

The competitiveness of $T$ will now follow once we have proved $(DT(x)_{\kappa(x)})^{-1}>0$ for all $x\in \mathbb{R}_+^n\setminus\{0\}$ by Proposition 4.1 in \cite{ruiz2011exclusion}. Recall that $\frac{\partial F_i(x)}{\partial x_j}<0$ for all $x\in \mathbb{R}^n_+$ and $i,j\in N$, so
the $(i,j)$-th entry of $DT(x)_{\kappa(x)}$ is negative for $i\neq j$. Then it follows from the proof of Theorem 3.1 in \cite{LG} that \eqref{condition-gamma-3} implies $(DT(x)_{\kappa(x)})^{-1}>0$ for all $x\in \mathbb{R}_+^n\setminus\{0\}$. This completes the proof.
\end{proof}

By \eqref{condition-gamma-3}, we know that  each map $T\in \mathrm{DCS}(n,f)$ satisfies the condition $\Upsilon$3) in Lemma \ref{simplex}. Since each map $T\in \mathrm{DCS}(n,f)$ also satisfies the conditions $\Upsilon$1) and $\Upsilon$2), it has a carrying simplex by Lemma \ref{simplex}, which had been proved in \cite[Corollary 3.3]{LG}.
\begin{lemma}\label{abstract-simplex}
Each map $T\in \mathrm{DCS}(n,f)$ admits a carrying simplex $\Sigma$.
\end{lemma}

\begin{remark}\label{remark:3-(ii)}
By Proposition \ref{prop-competitive} and Lemma \ref{abstract-simplex}, any map $T$ given by \eqref{equ:T2} is one-to-one and competitive on $\mathbb{R}_+^n$, and it has a carrying simplex unconditionally if each $f_i$ satisfies \eqref{cons:f} (i) and \eqref{cons:f} (ii). However, if \eqref{cons:f} (ii) does not hold for some $f_i$, then $T$ may not be one-to-one or competitive on $\mathbb{R}_+^n$. For example, the (Ricker) map $T$ with $f_i(z,r)=\exp(r-z)$, $i=1,\ldots,n$, is not one-to-one or competitive on $\mathbb{R}_+^n$ and in particular, it has a carrying simplex only under certain additional conditions; see \cite{GyllenbergRicker} for details.
\end{remark}

\begin{remark}\label{remark:DTp}
Let $T\in \mathrm{DCS}(n,f)$. If $T$ admits a unique positive fixed point $p=(p_1,\cdots,p_n)$, i.e.,
\begin{equation}\label{equ:p}
  (Ux^\tau)_i=1,\quad i=1,\cdots,n
\end{equation}
has a unique positive solution, then $1$ is not an eigenvalue of
\begin{equation}\label{equ:DTp}
DT(p)=I+\diag[p_i]\diag[\frac{\partial f_i}{\partial z}(r_i,r_i)]RU.
\end{equation}
Otherwise, $0$ is an eigenvalue of the matrix $DT(p)-I$, and hence $\det U=0$. Then \eqref{equ:p} has either no solution, or infinitely many solutions, a contradiction. Therefore, the index of $p$ which is given by $(-1)^{m}$ is either $1$ or $-1$, where $m$ is the sum of the multiplicities of all the eigenvalues of $DT(p)$ which are greater than one $($see {\rm\cite{Amann1976,Granas2003}}$)$. Let
$$
\mathcal{A}=-\diag[p_i]\diag[\frac{\partial f_i}{\partial z}(r_i,r_i)]RU, \quad \mathcal{B}=-\diag[\frac{\partial f_i}{\partial z}(r_i,r_i)]RU\diag[p_i].
$$
Then $DT(p)=I-\mathcal{A}$, and $\mathcal{A}$ is similar to $\mathcal{B}$. By the property \eqref{cons:f} {\rm(i)} of $f_i$, we know that $\mathcal{A},\mathcal{B}$ are positive matrices.
Note that $(Up^\tau)_i=1$, so the sum of the $i$-th row of $\mathcal{B}$ is $-r_i\frac{\partial f_i}{\partial z}(r_i,r_i)<1$ $($see \eqref{cons:f-2}$)$. It then follows from Perron-Frobenius theorem that $\rho(\mathcal{B})$, the spectral radius of $\mathcal{B}$, is an eigenvalue of $\mathcal{B}$ satisfying $0<\rho(\mathcal{B})<1$ and the magnitudes of the other eigenvalues of $\mathcal{B}$ are all less than $1$. Set $\lambda^*:=1-\rho(\mathcal{B})$. Since $\mathcal{A}$ and $\mathcal{B}$ have the same eigenvalues, $0<\lambda^*<1$ is a real eigenvalue of $DT(p)$ whose associated eigenvector is strictly positive and all the other eigenvalues possess real parts greater than $0$ and less than $2$. In particular, for the two-dimensional case, i.e. $n=2$, both of the two eigenvalues of $DT(p)$ are positive real numbers with one less than $1$, and $p$ is hyperbolic.
\end{remark}

In the remainder of this article, we will focus on analyzing the map $T\in\mathrm{DCS}(3,f)$ modeling three mutually competing species. We define an equivalence relation relative to local stability of fixed points on the boundary of $\Sigma$ for the set $\mathrm{DCS}(3,f)$ as that for all the
three dimensional Leslie-Gower maps
\cite{LG}
\begin{equation}\label{LG}
T: \mathbb{R}^3_+ \mapsto \mathbb{R}^3_+, \  T_i(x)=\frac{(1+r_i)x_i}{1+\sum_{j=1}^3a_{ij}x_j},\  r_i>0, a_{ij}=r_i\mu_{ij}>0, i,j=1,2,3.
\end{equation}
We show that the classification via this equivalence relation for three dimensional Leslie-Gower maps is valid for any $\mathrm{DCS}(3,f)$, and independent of the choice of generating function $f\in\mathcal{F}_3$. Furthermore, according to the equivalence classification, one can easily derive the permanence conditions in terms of simple inequalities on the parameters for $T\in\mathrm{DCS}(3,f)$.

\subsection{Classification via boundary dynamics}
In this subsection, we study the map $T \in \mathrm{DCS}(3,f)$:
\begin{equation}\label{map-T:3D}
T_i(x)=x_if_i(r_i\sum_{j=1}^3\mu_{ij}x_j,r_i)=x_if_i((RUx^\tau)_i,r_i), \quad i=1,2,3.
\end{equation}

It follows from Lemma \ref{abstract-simplex} that $T$ admits a 2-dimensional carrying simplex $\Sigma$ homeomorphic to $\Delta^{2}$. Each coordinate plane $\pi_i$ is positively invariant under $T$,
and the restriction of $T$ to $\pi_i$ is a 2-dimensional map $T|_{\pi_i} \in \mathrm{DCS}(2,f^{[i]})$, where $f^{[i]}=(f_j,f_k)\in \mathcal{F}_2$, $j<k$, so $\partial\Sigma$
is composed of the one-dimensional carrying simplices of $T|_{\pi_i}$. Therefore, before studying the three-dimensional map $T \in \mathrm{DCS}(3,f)$, we first study the two-dimensional case.
\subsubsection{The two-dimensional case}
Consider the map $T
\in \mathrm{DCS}(2,f)$:
\begin{equation}
T_i(x)=x_if_i(r_i\sum_{j=1}^2\mu_{ij}x_j,r_i)=x_if_i((RUx^\tau)_i,r_i), \quad i=1,2.
\end{equation}
By Lemma \ref{abstract-simplex}, $T$ admits a one-dimensional carrying simplex $\Sigma$ which is homeomorphic to the line segment joining the two points $(0,1)$ and $(1,0)$.
By Lemma \ref{abstract-simplex} and the arguments in Section \ref{subsection:permanence}, we conclude the following proposition.
\begin{proposition}\label{corollary:1}
Each map $T\in \mathrm{DCS}(2,f)$ has trivial dynamics, i.e., every nontrivial orbit converges to some fixed point on $\Sigma$.
\end{proposition}

Besides the trivial fixed point $0$ which is a hyperbolic repeller, $T$ admits two axial fixed points $q_{\{1\}}:(1/\mu_{11},0)$, $q_{\{2\}}:(0,1/\mu_{22})$. The fixed point $q_{\{i\}}$ is just the intersection of the line $\mathcal{S}_i=\{x\in\mathbb{R}^2_+: \mu_{i1}x_1+\mu_{i2}x_2=1\}$ and the $i$-th positive coordinate axis $\mathbb{H}^+_{i}$. If $\mathcal{S}_1$ and $\mathcal{S}_2$ intersect in $\dot{\mathbb{R}}_+^2$, then there also exists a positive fixed point $p$ at the intersection of $\mathcal{S}_1$ and $\mathcal{S}_2$.

Set $\mathbb{R}^2_+\setminus \mathcal{S}_i=\mathcal{U}_i \cup \mathcal{B}_i$, where $\mathcal{U}_i$ and $\mathcal{B}_i$ are the unbounded and bounded disjoint components of $\mathbb{R}^2_+\setminus \mathcal{S}_i$, respectively. Let $\gamma_{ij}:=\mu_{ii}-\mu_{ji}$ for $i,j=1,2$ and $i\neq j$. Then $q_{\{i\}}\in \mathcal{U}_j~($resp. $\mathcal{B}_j)~$ if and only if $\gamma_{ij}<0~($resp. $>0)$.

\begin{lemma}\label{lemma:3}
If $\gamma_{ij}>0$~$($resp. $<0)$, then $q_{\{i\}}$ is a saddle {\rm(}resp. an asymptotically stable node{\rm)}, and hence repels {\rm(}resp. attracts{\rm)} along $\Sigma$. Moreover, $q_{\{i\}}$ is hyperbolic if and only if $\gamma_{ij}\neq 0$.
\end{lemma}
\begin{proof}
Say $q_{\{1\}}$. The Jacobian matrix
$$
    DT(q_{\{1\}})=\left[
    \begin{array}{cc}
         1+r_1\frac{\partial f_1}{\partial z}(r_1,r_1) & \frac{r_1 \mu_{12}}{\mu_{11}}\frac{\partial f_1}{\partial z}(r_1,r_1) \\
        \noalign{\medskip}
        0 &  f_2(\frac{\mu_{21}r_2}{\mu_{11}},r_2)
    \end{array}
    \right],
$$
so $1+r_1\frac{\partial f_1}{\partial z}(r_1,r_1)$, $f_2(\frac{\mu_{21}r_2}{\mu_{11}},r_2)$ are its two positive eigenvalues. Note that $0<1+r_1\frac{\partial f_1}{\partial z}(r_1,r_1)<1$ and $\mathbb{H}^+_{1}$ is positively invariant, so every orbit emanating from $\mathbb{H}^+_{1}$ converges to $q_{\{1\}}$. If $f_2(\frac{\mu_{21}r_2}{\mu_{11}},r_2)>1~($resp. $<1)$, i.e., $\gamma_{12}>0~($resp. $<0)$, then $q_{\{1\}}$ is a saddle (resp. an asymptotically stable node), and hence repels (resp. attracts) along $\Sigma$. The last result is obvious.
\end{proof}


\begin{remark}\label{remark-2d-1}
In fact, the external eigenvalue at the axial fixed point $q_{\{i\}}$ in direction $j$ is $F_j(q_{\{i\}})=f_j((RUq^\tau_{\{i\}})_j,r_j)$ by  \eqref{jacobian-T}. Therefore, the sign of $r_j-(RUq^\tau_{\{i\}})_j$
is just the sign of $F_j(q_{\{i\}})-1$ (see the comments below \eqref{cons:f}), that is
\begin{equation}\label{eigenvalues:2D}
   \mathrm{sgn}(F_j(q_{\{i\}})-1)=\mathrm{sgn}(r_j-(RUq^\tau_{\{i\}})_j)=\mathrm{sgn}(\gamma_{ij}).
\end{equation}
Moreover, recall that $\gamma_{ij}>0~($resp.  $<0)$ if and only if $q_{\{i\}}\in \mathcal{B}_j~($resp. $\mathcal{U}_j)$. So the dynamics of the fixed point $q_{\{i\}}$ can be determined by the position of $q_{\{i\}}$ relative to the line $\mathcal{S}_j$, $i\neq j$. Moreover, if $\gamma_{12} \gamma_{21}>0~($resp. $<0)$, then $\mathcal{S}_1$ and $\mathcal{S}_2$ intersect {\rm(}resp. do not intersect{\rm)} in $\dot{\mathbb{R}}^2_+$, i.e., there exists {\rm(}resp. does not exist{\rm)} a positive fixed point $p$.
\end{remark}

Proposition \ref{prop:abstract-2} states that there are only four dynamical outcomes in $\mathrm{DCS}(2,f)$, which follows from Lemma \ref{lemma:3} and Remark \ref{remark-2d-1} directly. Just repeat the similar arguments in \cite[Theorem 4.1]{LG}.

\begin{proposition}\label{prop:abstract-2}
Let $T \in \mathrm{DCS}(2,f)$.
\begin{enumerate}[{\rm(a)}]
\item If $\gamma_{12}<0,\gamma_{21}>0$, then the positive fixed point $p$ does not exist and $q_{\{1\}}$ attracts all points not on the $x_2$-axis.
\item If $\gamma_{12}>0,\gamma_{21}<0$, then the positive fixed point $p$ does not exist and $q_{\{2\}}$ attracts all points not on the $x_1$-axis.
\item If $\gamma_{12},\gamma_{21}>0$, then $T$ has a hyperbolic positive fixed
point $p$ attracting all points in $\dot{\mathbb{R}}_+^2$.
\item If $\gamma_{12},\gamma_{21}<0$, then $T$ has a positive fixed point $p$ which is a hyperbolic saddle. Moreover, every nontrivial orbit tends to one of the asymptotically stable nodes $q_{\{1\}}$ or $q_{\{2\}}$ or to the saddle $p$.
\end{enumerate}
\end{proposition}

The following definition of equivalence appears to be unnecessarily pompous, but it prepares the way for the analogous definition in higher dimensions.

\begin{definition}\label{def-2d-1}
Two maps $T,\hat{T} \in \mathrm{DCS}(2,f)$ are said to be {\it equivalent relative to $\partial \Sigma$} if there exists a permutation $\sigma$ of $\{1,2\}$ such that $T$ has a
fixed point $q_{\{i\}}$ if and only if $\hat{T}$ has a fixed point
$\hat{q}_{\{\sigma(i)\}}$, and further 
$$
\mathrm{sgn}(F_j(q_{\{i\}})-1)=\mathrm{sgn}(\hat{F}_{\sigma(j)}(\hat{q}_{\{\sigma(i)\}})-1)
$$
for $j\neq i$, that is (see \eqref{eigenvalues:2D})
$$
\mathrm{sgn}(\gamma_{ij})=\mathrm{sgn}(\hat{\gamma}_{\sigma(i)\sigma(j)})
$$
for $j\neq i$.
\end{definition}

\begin{definition}
A map $T\in\mathrm{DCS}(2,f)$ is said to be {\it stable relative to $\partial \Sigma$} if
all the fixed points on $\partial \Sigma$ are hyperbolic. An equivalence class is said to be {\it stable} if each map in it is stable relative to $\partial \Sigma$.
\end{definition}
\begin{remark}\label{remark-2d-2}
Note that $\partial \Sigma \cap \mathcal{E}(T)=\{q_{\{1\}}, q_{\{2\}}\}$, so it follows from Lemma \ref{lemma:3} that $T$ is stable relative to $\partial \Sigma$ if and only if $\gamma_{12},\gamma_{21}\neq 0$, and hence an equivalence class is stable if there is a map in it which is stable relative to $\partial \Sigma$. 
Suppose that $T\in \mathrm{DCS}(2,f)$ is stable relative to $\partial \Sigma$ and possesses a positive fixed point $p$. Then the positive fixed point $p$ is  unique, and hence $\det U\neq 0$, where
\[
p=\left(\frac{\gamma_{21}}{\det U}, \frac{\gamma_{12}}{\det U}  \right).
\]
By the positivity of $p$, $\gamma_{12}$ and $\gamma_{21}$ have the same sign as $\det U$. Therefore, it follows from Proposition \ref{prop:abstract-2} (c) and (d) that $p$ attracts (resp. repels) on $\Sigma$ if and only if $\det U>0$ (resp. $\det U<0$). Moreover, it follows from Remark \ref{remark:DTp} that if $p$ attracts on $\Sigma$ then its two positive eigenvalues are less than $1$ while it has one eigenvalue greater than $1$ if it repels on $\Sigma$.
\end{remark}

By Proposition \ref{prop:abstract-2} we conclude the following result immediately.

\begin{figure}[ht]
 \begin{center}
    \includegraphics[width=0.92\textwidth]{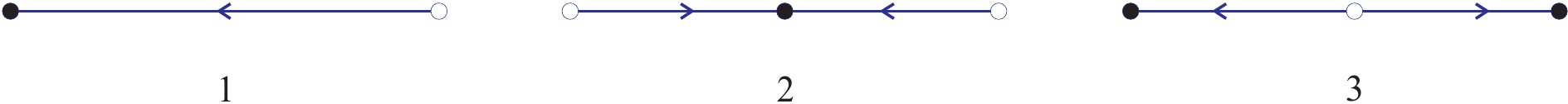}
 \end{center}
\caption{The phase portrait on $\Sigma$ replaced by $\Delta^1$. A closed dot
$\bullet$ denotes a fixed point which attracts on $\Sigma$, and an open dot $\circ$ denotes the one which repels on $\Sigma$. Each $\Sigma$ stands for an equivalence class. Class $1$ corresponds to Proposition \ref{prop:abstract-2} (a) and (b); class $2$ corresponds to Proposition \ref{prop:abstract-2} (c); class $3$  corresponds to Proposition \ref{prop:abstract-2} (d).} \label{fig:1}
\end{figure}

\begin{corollary}\label{coro-2d-1}
There are a total of $3$ stable equivalence classes in $\mathrm{DCS}(2,f)$. The three dynamical scenarios are presented in Fig. {\rm\ref{fig:1}}.
\end{corollary}

\begin{corollary}\label{coro-2d-2}
For a map $T\in \mathrm{DCS}(2,f)$ which is stable relative to $\partial \Sigma$, it is permanent if and only if it is in the stable class $2$, i.e. $\gamma_{12},\gamma_{21}>0$.
\end{corollary}
\begin{proof}
Note that $\gamma_{ij}>0~($resp. $<0)$ if and only if $f_j(\frac{\mu_{ji}r_j}{\mu_{ii}},r_j)>1~($resp. $<1)$, i.e., $F_j(q_{\{i\}})>1~($resp. $<1)$. So, it follows from the arguments in Section \ref{subsection:permanence} that $T$ is permanent if $\gamma_{12},\gamma_{21}>0$, and in this case there is a globally attracting positive fixed point. If some $\gamma_{ij}<0$, then $q_{\{i\}}$ is an attractor, so $T$ is impermanent in classes $1$ and $3$.
\end{proof}
\begin{remark}\label{remark-2d-3}
The statements of Proposition \ref{prop:abstract-2}, Corollaries \ref{coro-2d-1} and \ref{coro-2d-2} have clear biological interpretations.
\begin{enumerate}[{\rm (i)}]
\item If $\gamma_{ij}>0$, then species $j$ can invade species $i$ while it cannot invade if $\gamma_{ij} <0$.
\item If species $j$ can invade species $i$ but not vice versa, then species $i$ is driven to extinction, whilst species $j$ remains extant. In this case, the map is impermanent.
\item In the case of mutual invadability, that is, if both species can invade the other, then the map is permanent, and there will be coexistence in the form of an asymptotically stable positive fixed point.
\item If neither species can invade (mutual noninvadability), there is no coexistence: one of the species will oust the other. The surviving species depends on the initial conditions. (Convergence to the positive saddle happens only for initial conditions in a set of measure zero and is hence impossible in nature). In this case, the map is also impermanent.
\item When there is a positive fixed point, $\det U>0$  means in that both species can invade and the map is permanent, while $\det U<0$ means that none of them can  and the map is impermanent (Remark \ref{remark-2d-2}).
\end{enumerate}
The situations mentioned above are of particular interest when the two populations $1$ and $2$ are not different species, but different traits (resident and mutant) of the same species. To begin with, the resident ($i=1$) is at the fixed point $q_{\{1\}}$ and then the mutant $q_{\{2\}}$ is introduced in small quantities. Case (i) $\gamma_{12} >0$ gives the condition for successful invasion.  Case (ii) describes trait substitution.  Case (iii) is an example of protected dimorphism. For a discussion of these notions and their consequences for evolutionary dynamics we refer the reader to \cite{Geritz1997,Geritz1998,Geritz2002,Geritz2005}.
\end{remark}

\smallskip

\subsubsection{The three-dimensional case}
Now we analyze the $3$-dimensional map \eqref{map-T:3D}. We will define the equivalence relation on $\mathrm{DCS}(3,f)$ as Definition \ref{def-2d-1} and list the equivalence classification.

Besides the trivial fixed point $0$, $T$ has three axial
fixed points $q_{\{1\}}=(\frac{1}{\mu_{11}},0,0)$, $q_{\{2\}}=(0,\frac{1}{\mu_{22}},0)$, $q_{\{3\}}=(0,0,\frac{1}{\mu_{33}})$. In the interior of $\pi_k$, there may exist a planar fixed point $v_{\{k\}}$ satisfying
\begin{equation}\label{equ:v}
\mu_{ii}x_i+\mu_{ij}x_j+\mu_{ik}x_k=1, \  x_k=0, \ i\neq j  \neq k .
\end{equation}
In this case, $v_{\{k\}}$ is just the positive fixed point of the map $T|_{\pi_k}$.
$T$ may also admit a positive fixed point $p$ in
$\dot{\mathbb{R}}^3_+$ which satisfies
\begin{equation}\label{equ:equ2}
     \mu_{i1}x_1+\mu_{i2}x_2+\mu_{i3}x_3=1,\quad i=1,2,3.
\end{equation}
Hereafter, denote by
$$
\mathcal{S}_i=\{x\in\mathbb{R}^3_+:\mu_{i1}x_1+\mu_{i2}x_2+\mu_{i3}x_3=1\},~i=1,2,3.
$$
Let $\mathbb{R}^3_+\setminus \mathcal{S}_i=\mathcal{U}_i \cup \mathcal{B}_i$, where $\mathcal{U}_i$ and $\mathcal{B}_i$ are the unbounded and bounded disjoint components of $\mathbb{R}^3_+\setminus \mathcal{S}_i$, respectively. If $\mathcal{S}_i$ and $\mathcal{S}_j$ intersect in the interior of $\pi_k$, then $T$ has a fixed point $v_{\{k\}}$. There exists a positive fixed point $p$ if and only if $\mathcal{S}_1$, $\mathcal{S}_2$ and $\mathcal{S}_3$ intersect in $\dot{\mathbb{R}}^3_+$.


Let $\gamma_{ij}:=\mu_{ii}-\mu_{ji}$ for $i,j=1,2,3$ and $i\neq j$. By \eqref{jacobian-T}, we know that the external eigenvalue at the axial fixed point $q_{\{i\}}$ in direction $j$ is $F_j(q_{\{i\}})=f_j((RUq^\tau_{\{i\}})_j,r_j)$, and the external eigenvalue at the planar fixed point $v_{\{k\}}$ is $F_k(v_{\{k\}})=f_k((RUv_{\{k\}}^\tau)_k,r_k)$.
Therefore, the sign of $r_j-(RUq^\tau_{\{i\}})_j$
is just the sign of $F_j(q_{\{i\}})-1$, and that the sign of $r_k-(RUv_{\{k\}}^\tau)_k$ is just the sign of $F_k(v_{\{k\}})-1$ (see the comments below \eqref{cons:f}).  Specifically,
\begin{equation}\label{eigenvalues:3D}
\begin{array}{l}
      \mathrm{sgn}(F_j(q_{\{i\}})-1)=\mathrm{sgn}(r_j-(RUq^\tau_{\{i\}})_j)=\mathrm{sgn}(\gamma_{ij}), \\
      \noalign{\medskip}
      \mathrm{sgn}(F_k(v_{\{k\}})-1)=\mathrm{sgn}(r_k-(RUv_{\{k\}}^\tau)_k)=\mathrm{sgn}(1-(Uv_{\{k\}}^\tau)_k).
\end{array}
\end{equation}

By the positive invariance of $\pi_i$ and the analysis of the 2-dimensional case, the statements, proofs and classification program in \cite{LG} carry over to $\mathrm{DCS}(3,f)$ in a straightforward way, so we do not re-do it unless the need for special details and we only state the corresponding conclusions.

\begin{proposition}\label{lemma-3d-1}
If $\gamma_{ij}>0~($resp. $<0)$ then $q_{\{i\}}$ repels $($resp. attracts$)$ along $\partial \Sigma \cap \pi_k$, where $i,j,k$ are distinct. Furthermore, if $\gamma_{ij},\gamma_{ik}>0~($resp. $<0)$ then the fixed point $q_{\{i\}}$ is a repeller $($resp. an attractor$)$ on $\Sigma;$ if $\gamma_{ij}\gamma_{ik}<0$, then the fixed point $q_{\{i\}}$ is a saddle on $\Sigma;$ and $q_{\{i\}}$ is hyperbolic if and only if $\gamma_{ij}\gamma_{ik}\neq 0$.
\end{proposition}

\begin{proposition}\label{lemma-3d-2}
If $\gamma_{jk}\gamma_{kj}>0$~$($resp. $<0)$ then there is a unique $($resp. no$)$ fixed point $v_{\{i\}}$ in the interior of the coordinate plane $\pi_i$, where $i,j,k$ are distinct. Moreover, if $\gamma_{jk},\gamma_{kj}<0~($resp. $>0)$ then $v_{\{i\}}$ repels $($resp. attracts$)$ along $\partial \Sigma$.
\end{proposition}

The biological meaning of the condition $\gamma_{ij} > 0~({\rm resp.} <0)$ in Propositions \ref{lemma-3d-1} and \ref{lemma-3d-2} is that species $j$ can (resp. not) invade species $i$ in the absence of species $k$; here $i,j,k$ are distinct.

\begin{proposition}\label{lemma-3d-3}
Suppose the planar fixed point $v_{\{i\}}$ exists. Then  $(Uv_{\{i\}}^\tau)_i<1~($resp. $>1)$ implies that $v_{\{i\}}$ locally repels $($resp. attracts$)$ in $\dot{\Sigma}$. Moreover, $v_{\{i\}}$ is hyperbolic if and only if $(Uv_{\{i\}}^\tau)_i\neq 1$.
\end{proposition}

Propositions \ref{lemma-3d-1}--\ref{lemma-3d-3} imply that the local dynamics of $q_{\{i\}}$ and $v_{\{k\}}$ is generally determined by their external eigenvalues, i.e., $F_j(q_{\{i\}})$ ($j\neq i$) and $F_k(v_{\{k\}})$.
\begin{definition}
Two maps $T,\hat{T} \in \mathrm{DCS}(3,f)$ are said to be {\it equivalent relative to $\partial \Sigma$} if there exists a permutation $\sigma$ of $\{1,2,3\}$ such that 
\begin{enumerate}[{\rm (i)}]
\item $T$ has a fixed point $q_{\{i\}}$ if and only if $\hat{T}$ has a fixed
point $\hat{q}_{\{\sigma(i)\}}$, and further
$$
\mathrm{sgn}(F_j(q_{\{i\}})-1)=\mathrm{sgn}(\hat{F}_{\sigma(j)}(\hat{q}_{\{\sigma(i)\}})-1)
$$
for all $j\neq i$, that is (see \eqref{eigenvalues:3D})
$$
\mathrm{sgn}(\gamma_{ij})=\mathrm{sgn}(\hat{\gamma}_{\sigma(i)\sigma(j)})
$$
for all $j\neq i$;
\item $T$ has a fixed point $v_{\{k\}}$ if and only if $\hat{T}$ has a fixed
point $\hat{v}_{\{\sigma(k)\}}$, and further
$$
\mathrm{sgn}(F_k(v_{\{k\}})-1)=\mathrm{sgn}(\hat{F}_{\sigma(k)}(\hat{v}_{\{\sigma(k)\}})-1),
$$
that is (see \eqref{eigenvalues:3D})
$$
\mathrm{sgn}(1-(Uv_{\{k\}}^\tau)_k)=\mathrm{sgn}(1-(\hat{U}\hat{v}_{\{\sigma(k)\}}^\tau)_{\sigma(k)}).
$$
\end{enumerate}
\end{definition}

\begin{definition}
A map $T\in\mathrm{DCS}(3,f)$ is said to be {\it stable relative to $\partial \Sigma$} if
all the fixed points on $\partial \Sigma$ are hyperbolic. An
equivalence class is said to be {\it stable} if each map in
it is stable relative to $\partial \Sigma$.
\end{definition}

\begin{remark}\label{remark-3d-1}
 By Propositions \ref{lemma-3d-1} and \ref{lemma-3d-3}, a map $T \in \mathrm{DCS}(3,f)$ is stable relative to $\partial \Sigma$ if and only if $\gamma_{ij} \neq 0$ and $(Uv_{\{k\}}^\tau)_k\neq 1$ {\rm(}if $v_{\{k\}}$ exists{\rm)} for $i,j,k=1,2,3$ and $i\neq j$, and hence an equivalence class is stable if there is a map in it which is stable relative to $\partial \Sigma$. 

Suppose $\gamma_{ij},\gamma_{ji}\neq 0$ (here $i\neq j$). It follows from Proposition \ref{lemma-3d-2} that $v_{\{k\}}$ exists if and only if $\gamma_{ij}\gamma_{ji}>0$, which implies that $\det U_{\{i,j\}}\neq 0$ (i.e., $\mu_{ii}\mu_{jj}-\mu_{ij}\mu_{ji}\neq 0$) by noticing that  $v_{\{k\}}$ is the unique positive fixed point of $T|_{\pi_k}$ (see Remark \ref{remark-2d-2}), where 
$$U_{\{i,j\}}=\left[\begin{array}{cc}
\mu_{ii}&\mu_{ij}\\
\mu_{ji}&\mu_{jj}
\end{array}\right].$$
Therefore, for a map $T \in \mathrm{DCS}(3,f)$ which is stable relative to $\partial \Sigma$, if $v_{\{k\}}$ exists then $\mu_{ii}\mu_{jj}-\mu_{ij}\mu_{ji}\neq 0$ (here $i,j,k$ are distinct), and it
is easy to check that
\begin{equation}\label{equ:beta}
(Uv_{\{k\}}^\tau)_k<1~(>1) \Leftrightarrow \mu_{ki}\beta_{ij}+\mu_{kj}\beta_{ji}<1~(>1)\Leftrightarrow v_{\{k\}}\in \mathcal{B}_k~(\mathcal{U}_k),
\end{equation}
where 
$$\beta_{ij}:=\frac{\mu_{jj}-\mu_{ij}}{\mu_{ii}\mu_{jj}-\mu_{ij}\mu_{ji}}.$$
Thus a map $T \in \mathrm{DCS}(3,f)$ is stable relative to $\partial \Sigma$ if and only if $\gamma_{ij} \neq 0$ and $\mu_{ki}\beta_{ij}+\mu_{kj}\beta_{ji}\neq 1$, i.e., $(Uv_{\{k\}}^\tau)_k\neq 1$ {\rm(}if $v_{\{k\}}$ exists{\rm)}. Suppose that $T$ is stable relative to $\partial \Sigma$. It follows from Propositions \ref{lemma-3d-1}--\ref{lemma-3d-3} and \eqref{equ:beta} that
the existence and local dynamics of boundary fixed points on $\partial \Sigma$ for $T$ are completely determined by the parameters $\mu_{ij}$, i.e. the values $\gamma_{ij}$ and $\mu_{ki}\beta_{ij}+\mu_{kj}\beta_{ji}$, which are independent of the generating function $f$.

Moreover, if $T$ admits a positive fixed point
$p$ which satisfies \eqref{equ:equ2}, then $p$ is the unique positive fixed point. Otherwise, assume that $T$ has two different positive fixed points $p$ and $\tilde{p}$. Now $p_s:=sp+(1-s)\tilde{p}$ is a solution of \eqref{equ:equ2} for any $s\geq 0$. Let $\bar{s}:=\sup\{s > 0: p_s\in \Sigma\}$. Then $p_{\bar{s}}\in \partial \Sigma$ is a fixed point, which
is not hyperbolic, contradicting that $T$ is stable relative to $\partial \Sigma$. Thus, $1$ is not an eigenvalue of $DT(p)$ by Remark \ref{remark:DTp}.  Therefore, $T$ has only finitely many fixed points on $\Sigma$, i.e. three axial fixed points $q_{\{i\}}$, at most three planar fixed points $v_{\{i\}}$ and at most one positive fixed point $p$, and $1$ is not an eigenvalue of any of their Jacobian matrices.
\end{remark}

Let $Q=id-T$, where $id$ is the identity mapping.  Let $x$ be a fixed point of $T$, that is, a  zero of $Q$. The index of $T$ at $x$ is denoted by $\mathrm{Ind}(x,T)$ and the index of $Q$ at the zero $x$ is denoted by $\mathscr{I}(x,Q)$. The index $\mathscr{I}(x,Q)$ is defined as the sign of $\det DQ(x)$ if $\det DQ(x)\neq 0$, and the index $\mathrm{Ind}(x,T)$ as $\mathscr{I}(x,Q)$; for the general theory see \cite{Granas2003}.

\begin{lemma}[Index Formula on Carrying Simplex \cite{jiang2014}]\label{theory:index}
Suppose that $T:\mathbb{R}^3_+\to \mathbb{R}^3_+$ given by \eqref{map_T_dimn} satisfies $\partial F_i/\partial x_j< 0$ for all $x\in \mathbb{R}_+^3$. Assume that $T$ possesses a carrying simplex $\Sigma$ and the continuous-time system $\dot{x}=G(x)=T(x)-x$
is dissipative with the origin $0$ being a repeller. If $T$ has only finitely many fixed points on $\Sigma$
and $1$ is not an eigenvalue of any of their Jacobian matrices, then
$$
    \sum_{\hat{x}\in \mathcal{E}_v}\mathrm{Ind}(\hat{x},T)+2\sum_{\hat{x}\in \mathcal{E}_s}\mathrm{Ind}(\hat{x},T)+4\sum_{\hat{x}\in \mathcal{E}_p}\mathrm{Ind}(\hat{x},T)=1,
$$
where $\mathcal{E}_v$, $\mathcal{E}_s$, and $\mathcal{E}_p$ denote the set of all nontrivial axial, planar, and positive fixed points, respectively.
\end{lemma}

\begin{proposition} \label{prop:index}
Assume that $T \in \mathrm{DCS}(3,f)$ is stable relative to $\partial \Sigma$. Then we have the formula
\begin{equation}\label{equ:1}
    \sum^3_{i=1}(\mathrm{Ind}(q_{\{i\}},T)+2\mathrm{Ind}(v_{\{i\}},T))+4\mathrm{Ind}(p,T)=1.
\end{equation}
\end{proposition}
\begin{proof}
Let $G(x)=T(x)-x$. Consider the continuous-time system
\begin{equation}\label{continuous-discrete}
\dot{x}_i=G_i(x)=x_i(F_i(x)-1),\quad i=1,2,3.
\end{equation}
The origin $0$ is an equilibrium of system
\eqref{continuous-discrete}, and the eigenvalues of $DG(0)$ are
$F_i(0)-1>0$, that is, $0$ is a repeller. Note that
$$\frac{\partial
F_i}{\partial x_j}=\frac{\partial f_i}{\partial z}(r_i\sum_{j=1}^3\mu_{ij}x_j,r_i)r_i\mu_{ij}<0,$$
so system \eqref{continuous-discrete} is totally competitive. Since
$$G_i(x)=x_i(F_i(x)-1)=x_i(f_i(r_i\sum_{j=1}^3\mu_{ij}x_j,r_i)-1)<0,\quad i=1,2,3, $$
for $|x|$ sufficiently large, so system \eqref{continuous-discrete} is
dissipative. Recall that $T \in \mathrm{DCS}(3,f)$ admits a carrying simplex, so the result follows from Remark \ref{remark-3d-1} and Lemma \ref{theory:index}.
\end{proof}
\begin{lemma} \label{lemma:04}
Suppose that $T \in \mathrm{DCS}(3,f)$ is stable relative to $\partial \Sigma$. Then
\begin{enumerate}[{\rm (i)}]
\item $\mathrm{Ind}(q_{\{i\}},T)=1$ {\rm(}resp.\ $\mathrm{Ind}(v_{\{k\}},T)=1${\rm)} if $q_{\{i\}}$ {\rm(}resp.\ $v_{\{k\}}${\rm)} is a repeller or an attractor on $\Sigma$;
\item $\mathrm{Ind}(q_{\{i\}},T)=-1$ {\rm(}resp.\ $\mathrm{Ind}(v_{\{k\}},T)=-1${\rm)} if $q_{\{i\}}$ {\rm(}resp.\ $v_{\{k\}}${\rm)} is a saddle on $\Sigma$;
\item $\mathrm{Ind}(p,T)\neq 0$ if the positive fixed point $p$ exists.
\end{enumerate}
\end{lemma}
\begin{proof}
It follows from the analysis for the two-dimensional maps, \eqref{eigenvalues:3D} and Remark \ref{remark:DTp} that all the eigenvalues of $q_{\{i\}}$ and $v_{\{k\}}$ (if any) are positive real numbers and do not equal $1$. If $q_{\{i\}}$ (resp. $v_{\{k\}}$) is a repeller or an attractor on $\Sigma$ then the number of the eigenvalues of $DT(q_{\{i\}})$ (resp. $DT(v_{\{k\}})$) greater than $1$ is even, and hence $\mathrm{Ind}(q_{\{i\}},T)=1$ (resp. $\mathrm{Ind}(v_{\{k\}},T)=1$). If $q_{\{i\}}$ (resp. $v_{\{k\}}$) is a saddle on $\Sigma$ then the number of the eigenvalues of $DT(q_{\{i\}})$ (resp. $DT(v_{\{k\}})$) greater than $1$ is odd, and hence $\mathrm{Ind}(q_{\{i\}},T)=-1$ (resp. $\mathrm{Ind}(v_{\{k\}},T)=-1$). If there is a positive fixed point $p$, then it follows from Remark \ref{remark-3d-1} that it is unique and $1$ is not an eigenvalue of $DT(p)$. Thus, $\mathrm{Ind}(p,T)\neq 0$.
\end{proof}
\begin{remark}\label{remark:4.4}
For a map $T \in \mathrm{DCS}(3,f)$ which is stable relative to $\partial \Sigma$, it follows from Proposition \ref{prop:index} and Lemma \ref{lemma:04} that the existence of the positive fixed point $p$ and its index can be determined by the local dynamics of boundary fixed points.
\end{remark}

\begin{theorem}\label{theory:classification}
There are a total of $33$ stable equivalence classes in $\mathrm{DCS}(3,f)$, where the parameter conditions for each class with the corresponding phase portrait on the carrying simplex are listed in Table \ref{biao0}.
\end{theorem}

Recalling Remark \ref{remark-3d-1}, the existence and local dynamics of boundary fixed points on $\partial \Sigma$ for $T\in\mathrm{DCS}(3,f)$ are completely determined by the parameters $\mu_{ij}$, i.e. the values $\gamma_{ij}$ and $\mu_{ki}\beta_{ij}+\mu_{kj}\beta_{ji}$, which are independent of $f$, and the same as the Leslie-Gower map \eqref{LG}.
Therefore, the classifications are the same for them, which are independent of the choice of the generating function $f\in\mathcal{F}_3$. Any stable map in $\mathrm{DCS}(3,f)$ belongs to one of the $33$ classes in Table \ref{biao0} (modulo permutation of the indices). Moreover, there is no positive fixed point in classes $1-18$, which have trivial dynamics, i.e. every orbit converges to some fixed point. Each map from classes $19-25$ admits a unique positive fixed point with index $-1$, and every orbit also converges to some fixed point for these classes. Each map in classes $26-33$ has a unique positive fixed point with index $1$; and the positive fixed point is globally asymptotically stable in class 33; see Subsection \ref{stability} for details. Such a classification is also valid for the Ricker models admitting a carrying simplex \cite{GyllenbergRicker}, and we will discuss in Section \ref{sec:application}.
\subsection{Stability and permanence}\label{stability}
As befits the context, we shall consider the families of maps given in Table \ref{biao0} by permutation of the indices,
 i.e., we assume the parameters $\mu_{ij},r_i$ of the corresponding class satisfy the conditions listed in Table \ref{biao0}.

\begin{lemma}[Theorem 2.2 in \cite{ruiz2011exclusion} and Theorem 3.1 in \cite{nonlinearity2018}]\label{trivial-lemma}
Consider the three-dimensional map $T$ given by \eqref{map_T_dimn} which satisfies the conditions $\Upsilon 1)$, $\Upsilon 2)$ and $\Upsilon 3)$ in Lemma \ref{simplex}. Suppose that $T$ has only a finite number of fixed points. Then the following conclusions hold:

$\bullet$ If $T$ has no positive fixed point, then every nontrivial orbit converges to some fixed point on the boundary of the carrying simplex.

$\bullet$ If $T$ has a unique positive fixed point $p$ such that $\mathrm{Ind}(p,T)=-1$, then $p$ is a saddle on the carrying simplex, and moreover, every nontrivial orbit converges to some fixed point on the boundary of the carrying simplex, except those on the stable manifold of $p$.
\end{lemma}
\begin{remark}\label{remark:invariant-manifold}
For the three-dimensional map $T$ in Lemma \ref{trivial-lemma} which has a unique positive fixed point $p$ such that $\mathrm{Ind}(p,T)=-1$, it is proved in \cite{arXiv2019} that both the stable manifold and unstable manifold of the saddle $p$ are simple curves, and the phase portrait on the carrying simplex can be described clearly together with the dynamics on the boundary of the carrying simplex.
\end{remark}
Recall that each map $T\in\mathrm{DCS}(3,f)$ satisfies the conditions $\Upsilon 1)$, $\Upsilon 2)$ and $\Upsilon 3)$ in Lemma \ref{simplex}. Moreover, if $T$ is stable relative to $\partial \Sigma$, then there is at most one positive fixed point, say $p$, and $\mathrm{Ind}(p,T)\neq 0$ if $p$ exists. Thus, together with Lemma \ref{trivial-lemma}, Proposition \ref{prop:index} and Remark \ref{remark-3d-1} imply the following trivial dynamics via boundary fixed points.
\begin{proposition}\label{prop:trivial-dyns}
Assume that $T\in\mathrm{DCS}(3,f)$ is stable relative to $\partial \Sigma$. Suppose that
\begin{equation}\label{index-trivial}
\sum^3_{i=1}(\mathrm{Ind}(q_{\{i\}},T)+2\mathrm{Ind}(v_{\{i\}},T))=1
\end{equation}
or
\begin{equation}\label{index-trivial-2}
\sum^3_{i=1}(\mathrm{Ind}(q_{\{i\}},T)+2\mathrm{Ind}(v_{\{i\}},T))=5.
\end{equation}
Then $T$ has trivial dynamics, i.e. every nontrivial orbit converges to some fixed point on $\Sigma$.
\end{proposition}

Note that, each map $T$ in classes $1-18$ satisfies \eqref{index-trivial} in Proposition \ref{prop:trivial-dyns}, and hence $T$ has no positive fixed point. Therefore, such $T$ has trivial dynamics. That is, we have the following proposition.
\begin{proposition}\label{prop:3}
For each map $T$ in classes $1-18$, every nontrivial orbit converges to some fixed point on $\partial \Sigma$.
\end{proposition}

In the biological sense, Proposition \ref{prop:3} means that for three competing species modeled by $T$, if there is no coexistence state, then some of the species will be extinct.

For each map $T$ in classes $19-33$, there exits a unique positive fixed point $p$. Recall that $DT(p)=I-\mathcal{A}$, where
$$
\mathcal{A}=-\diag[p_i]\diag[\frac{\partial f_i}{\partial z}(r_i,r_i)]RU.
$$
\begin{lemma}\label{lemma:11}
For each map in classes $19-25$, we have $\mathrm{Ind}(p,T)=-1$ and $\det U<0;$
while for each map in classes $26-33$, we have $\mathrm{Ind}(p,T)=1$ and $\det U>0$.
\end{lemma}
\begin{proof}
For classes $19-25$ (resp. classes $26-33$), it follows from the local dynamics of fixed points on $\partial\Sigma$ in Table \ref{biao0}, Lemma \ref{lemma:04} and formula \eqref{equ:1} that $\mathrm{Ind}(p,T)=-1$ (resp. $\mathrm{Ind}(p,T)=1$). Moreover, if $\mathrm{Ind}(p,T)=-1$, then all the three eigenvalues of $DT(p)$ are positive real numbers with one eigenvalue greater than $1$ and the other two less than $1$ by Remark \ref{remark:DTp}. So, two eigenvalues of $\mathcal{A}$ are greater than $0$ and one is less than $0$, which implies that $\det \mathcal{A}<0$, and hence $\det U<0$. While $\mathrm{Ind}(p,T)=1$ ensures that there are
 zero or two eigenvalues of $DT(p)$ greater than $1$ by Remark \ref{remark:DTp}. For the former case, also by Remark \ref{remark:DTp} we have one eigenvalue  of $\mathcal{A}$ is greater than $0$ and the other two are either complex numbers or greater than $0$. For the latter case, two eigenvalues of $\mathcal{A}$ are less than $0$ and one is greater than $0$. Therefore, one always has $\det \mathcal{A}>0$, and hence $\det U>0$.
\qquad\end{proof}

\begin{proposition}\label{prop:5a}
 The positive fixed point $p$ is a saddle on $\Sigma$ in classes $19-25$, and every nontrivial orbit converges to some fixed point on the boundary of the carrying simplex, except those on the stable manifold of $p$.
\end{proposition}
\begin{proof}
Since $\mathrm{Ind}(p,T)=-1$ implies \eqref{index-trivial-2} holds, the result is immediate from Proposition \ref{prop:trivial-dyns} and Lemma \ref{lemma:11}.
\end{proof}

\begin{proposition}\label{prop:5}
 The positive fixed point $p$ is a repeller on $\Sigma$ in class $32$.
\end{proposition}
\begin{proof}
For each map $T$ in class $32$, there exists a planar fixed point $v_{\{k\}}$ in the interior of $\pi_k$ for each $k=1,2,3$, which is repelling along $\partial \Sigma \cap \pi_k$ (see Table \ref{biao0} (32)), so $v_{\{k\}}$ is a saddle for $T|_{\pi_k}$.
It then follows from Remark \ref{remark-2d-2} that
$\det U_{\{i,j\}}<0$ for any $i<j$, where
$U_{\{i,j\}}=\left[\begin{array}{cc}
\mu_{ii}&\mu_{ij}\\
\mu_{ji}&\mu_{jj}
\end{array}\right]$
is the principal $2\times 2$ submatrix of $U$. Therefore, $\det \mathcal{A}_{\{i,j\}}<0$ for any $i<j$. By $\det U>0$, one also has $\det \mathcal{A}>0$. It follows from Proposition 3.8 in \cite{Z993} and $\det \mathcal{A}>0$ that $\mathcal{A}$ has two eigenvalues with negative real parts. Therefore, $DT(p)=I-\mathcal{A}$ has two eigenvalues with real parts greater than $1$, i.e, $DT(p)$ has two eigenvalues with magnitudes greater than $1$ except $\lambda^*$, where $0<\lambda^*<1$ is defined in Remark \ref{remark:DTp}. So $p$ is a hyperbolic fixed point and it follows from \cite[Theorem 4.6]{arXiv2019} that the local dynamics of $p$ on $\Sigma$ is reflected by the other two eigenvalues except $\lambda^*$, which implies that $p$ is a repeller on $\Sigma$ and its two-dimensional unstable manifold is contained in $\Sigma$ (see \cite[Corollary 4.5]{arXiv2019}).
\end{proof}

The following lemma is the 3D specialization of Theorem 2.4 in \cite{Balreira2017} (see also Theorem 1.2 in \cite{Gyllenberg2019}), which can be used to establish our global stability for class $33$.
\begin{lemma}[\cite{Balreira2017}]\label{global-criterion}
Consider the three-dimensional map $T:\mathbb{R}^3_+\mapsto \mathbb{R}^3_+$ given by \eqref{map-T}, where $F_i$ are $C^1$ satisfying $F_i(x)>0$ for all $x\in \mathbb{R}_+^3$, $i=1,2,3$. Assume that
\begin{itemize}
  \item[{\rm(a)}] $\det DT(x)>0$ for all $x\in \mathbb{R}^3_+$;
  \item[{\rm(b)}] $DT(x)^{-1}>0$ for all $x\in \dot{\mathbb{R}}^3_+$;
  \item[{\rm(c)}] for each $i=1,2,3$, $T|_{\pi_i}$ has a unique interior fixed point $v_{\{i\}}$ that is globally asymptotically
stable in the interior of $\pi_i$, but a saddle for $T$;
  \item[{\rm(d)}] $T$ admits a carrying simplex;
  \item[{\rm(e)}] $T$ has a unique positive fixed point $p\in \dot{\mathbb{R}}^3_+$.
\end{itemize}
Then $p$ is globally asymptotically stable in $\dot{\mathbb{R}}^3_+$ for $T$.
\end{lemma}

\begin{figure}[h]
 \begin{center}
 \includegraphics[width=0.36\textwidth]{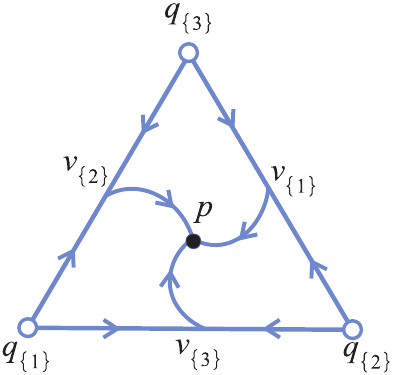}
\caption{The phase portrait on $\Sigma$ for class $33$. Every orbit in the interior of $\Sigma$ converges to $p$. The fixed point notation is as in Table \ref{biao0}.} \label{fig-global-33}
  \end{center}
\end{figure}
\begin{theorem}\label{global-stability}
The positive fixed point $p$ is globally asymptotically stable in $\dot{\mathbb{R}}^3_+$ for each map $T$ in class $33$, and the phase portrait on $\Sigma$ is as shown in Fig. \ref{fig-global-33}.
\end{theorem}
\begin{proof}
By Table \ref{biao0}, the map $T\in\mathrm{DCS}(3,f)$ is in class $33$ if the parameters satisfy the following inequalities
\begin{itemize}
  \item[{\rm(i)}] $\gamma_{12}>0, \gamma_{13}>0, \gamma_{21}>0,\gamma_{23}>0, \gamma_{31}>0, \gamma_{32}>0;$
  \item[{\rm(ii)}] $\mu_{12}\beta_{23}+\mu_{13}\beta_{32}<1;$
  \item[{\rm(iii)}] $\mu_{21}\beta_{13}+\mu_{23}\beta_{31}<1;$
  \item[{\rm(iv)}] $\mu_{31}\beta_{12}+\mu_{32}\beta_{21}<1$.
\end{itemize}
Besides the three axial fixed points $q_{\{1\}},q_{\{2\}}$ and $q_{\{3\}}$, which are all local repellers by (i), $T$ has three planar fixed points $v_{\{1\}},v_{\{2\}}$ and $v_{\{3\}}$ and a unique positive fixed point $p\in \dot{\mathbb{R}}^3_+$. By Proposition \ref{prop:abstract-2} (c) and Remark \ref{remark-2d-2}, each $v_{\{i\}}$ is globally asymptotically
stable for $T|_{\pi_i}$ in the interior of $\pi_i$ and the two internal eigenvalues of $DT(v_{\{i\}})$ are both positive and less than one.  Conditions (ii)-(iv) and \eqref{eigenvalues:3D} imply that the external eigenvalue of $DT(v_{\{i\}})$, i.e. $F_i(v_{\{i\}})$, is greater than one for each $v_{\{i\}}$, that is each $v_{\{i\}}$ is a saddle for $T$. Thus, the condition (c) in Lemma \ref{global-criterion} holds for $T$. By Proposition \ref{prop-competitive}, we know that $\det DT(x)>0$ for all $x\in \mathbb{R}^3_+$ and $DT(x)^{-1}>0$ for all $x\in \dot{\mathbb{R}}^3_+$, that is conditions (a) and (b) in Lemma \ref{global-criterion} hold for $T$. Therefore, the conclusion follows from Lemma \ref{global-criterion} immediately, because conditions (d) and (e) hold naturally for each map $T$ in class $33$.
\end{proof}

\begin{remark}
Propositions \ref{prop:3} and \ref{prop:5a} and Theorem \ref{global-stability} imply that nontrivial dynamics, e.g. bifurcations and invariant circles, can only occur in classes $26-32$. Proposition \ref{prop:5} implies that the positive fixed point $p$ in class $32$ is always hyperbolic, and has no eigenvalues of modulus $1$. So Neimark-Sacker bifurcations cannot occur in class $32$. However, within classes $26-31$, Neimark-Sacker bifurcations may occur for some specific $f\in\mathcal{F}_3$, such as the Atkinson-Allen model \cite{jiang2014,GyllenbergCGAA} and the Leslie-Gower model \cite{LG}; see Section \ref{sec:application} for details.
\end{remark}

For any map $T$ in class $27$, each axial fixed point $q_{\{i\}}$ is a saddle on $\Sigma$, and $\partial \Sigma\cap \pi_i$ is the heteroclinic connection between $q_{\{j\}}$ and $q_{\{k\}}$, where $i,j,k$ are distinct. Therefore, $\partial \Sigma$ is a heteroclinic cycle of May-Leonard type: $q_{\{1\}}\to q_{\{2\}} \to q_{\{3\}}\to q_{\{1\}}$ (or the arrows reversed), i.e., any map $T$ in class $27$ admits a heteroclinic cycle (see Table \ref{biao0} (27)).

Set $\mathcal{G}_{ij}=\ln F_j(q_{\{i\}})=\ln f_j((RUq^\tau_{\{i\}})_j,r_j)$, where $i\neq j$. Now the $\varrho$ which is defined in \eqref{heteroclinic-con} is written as
\begin{equation}\label{stability-rho}
\varrho=\mathcal{G}_{12}\mathcal{G}_{23} \mathcal{G}_{31}+\mathcal{G}_{21} \mathcal{G}_{13}\mathcal{G}_{32}.
\end{equation}

\begin{proposition}\label{stable_cycle}
Assume that $T\in\mathrm{DCS}(3,f)$ is in class $27$. If $\varrho>0~(resp. <0)$, then the heteroclinic cycle $\partial \Sigma$ of $T$ repels {\rm(}resp. attracts{\rm )}.
\end{proposition}
\begin{proof}
The conclusion follows from Lemma \ref{heteroclinic-cycle} immediately.
\end{proof}

From a biological point of view, these cycles in class $27$ may be seen to correspond to the biological environment where in purely pairwise competition species $2$ can invade species $1$ but not vice versa,
species $3$ can invade species $2$ but not vice versa, and species $1$ can invade species $3$ but not vice versa. It is this intransitivity in the pairwise competition, which underlies the cycle behavior. $\mathcal{G}_{ij}>0$ (resp. $<0$) means that species $j$ can (resp. not) invade species $i$; see Remark \ref{remark-2d-3}.

\begin{proposition}\label{permanence-3D}
Assume that $T\in\mathrm{DCS}(3,f)$ is stable relative to $\partial \Sigma$. Then
\begin{enumerate}[{\rm (i)}]
\item $T$ is permanent if it is in classes $29$, $31$, $33$ and class $27$ with $\varrho>0$;
\item $T$ is impermanent if it is in classes $1-26$, $28$, $30$, $32$ and class $27$ with $\varrho<0$.
\end{enumerate}
\end{proposition}
\begin{proof}
(i) Since the proofs for classes $29$, $31$ and $33$ are completely analogous, we only consider the class $29$; see Fig. \ref{fig:c29}.
\begin{figure}[h]
 \begin{center}
 \includegraphics[width=0.36\textwidth]{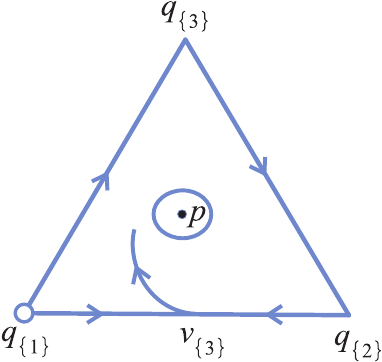}
\caption{The phase portrait on $\Sigma$ for class $29$. The fixed point notation is as in Table \ref{biao0}.} \label{fig:c29}
  \end{center}
\end{figure}\\
Note that $$\mathcal{E}(T)\cap \partial \Sigma=\{q_{\{1\}},q_{\{2\}},q_{\{3\}},v_{\{3\}}\}.$$
By Corollary \ref{coro:3-permanence}, it suffices to prove that
there are real numbers $\nu_1,\nu_2,\nu_3>0$ such that the following inequalities hold:
\begin{subequations}\label{inequalities}
\begin{align}
  \nu_1\ln F_1(q_{\{1\}})+\nu_2\ln F_2(q_{\{1\}})+\nu_3\ln F_3(q_{\{1\}})>0;\label{ineq:aa} \\[2pt]
  \nu_1\ln F_1(q_{\{2\}})+\nu_2\ln F_2(q_{\{2\}})+\nu_3\ln F_3(q_{\{2\}})>0;\label{ineq:bb}\\[2pt]
  \nu_1\ln F_1(q_{\{3\}})+\nu_2\ln F_2(q_{\{3\}})+\nu_3\ln F_3(q_{\{3\}})>0;\label{ineq:cc}\\[2pt]
  \nu_1\ln F_1(v_{\{3\}})+\nu_2\ln F_2(v_{\{3\}})+\nu_3\ln F_3(v_{\{3\}})>0.
  \label{ineq:dd}
\end{align}
\end{subequations}
Recall that for a fixed point $\hat{x}\in \mathcal{E}(T)$, one has $F_i(\hat{x})=1$ for all $i\in \kappa(\hat{x})$. Therefore,
$F_i(q_{\{i\}})=1$, $i=1,2,3$ and $F_1(v_{\{3\}})=F_2(v_{\{3\}})=1$. By \eqref{eigenvalues:3D}, Remark \ref{remark-3d-1} and the condition (ii) in Table \ref{biao0} (29), we have
$$ \mu_{31}\beta_{12}+\mu_{32}\beta_{21}<1 \Leftrightarrow  (Uv_{\{3\}}^\tau)_3<1\Leftrightarrow F_3(v_{\{3\}})>1.$$
So, \eqref{ineq:dd} holds for any $\nu_1,\nu_2,\nu_3>0$.
Since $\gamma_{12},\gamma_{13}>0$ by condition (i) in Table (29), one has $F_2(q_{\{1\}}),F_3(q_{\{1\}})>1$ (see \eqref{eigenvalues:3D}). Thus, \eqref{ineq:aa} holds for any $\nu_1,\nu_2,\nu_3>0$. The inequalities \eqref{ineq:bb} and \eqref{ineq:cc} can be written as
\begin{subequations}\label{inequalities-2}
\begin{align}
  \nu_1\ln F_1(q_{\{2\}})+\nu_3\ln F_3(q_{\{2\}})>0;\label{ineq:aa-2}\\[2pt]
  \nu_1\ln F_1(q_{\{3\}})+\nu_2\ln F_2(q_{\{3\}})>0.\label{ineq:bb-2}
\end{align}
\end{subequations}
We first fix a $\nu_2>0$. It follows from $\gamma_{32}>0$ and \eqref{eigenvalues:3D} that $\ln F_2(q_{\{3\}})>0$, and hence for sufficiently small $\nu_1>0$ one has
\eqref{ineq:bb-2} holds. Now fix some $\nu_1>0$ such that \eqref{ineq:bb-2} holds. Note that $\gamma_{21}>0$, so $\ln F_1(q_{\{2\}})>0$ (see \eqref{eigenvalues:3D}). Then we can choose some $\nu_3>0$ sufficiently small such that \eqref{ineq:aa-2} holds.  Such $\nu_1,\nu_2,\nu_3>0$ ensure that the inequalities \eqref{ineq:aa}--\eqref{ineq:dd} hold. This proves that each map $T$ in class $29$ is permanent. For the map $T$ in class $27$ such that $\varrho>0$, the conclusion follows from Corollary \ref{stability-permance}.

(ii) For each map $T$ in classes $1-26$, $28$, $30$ and $32$, there always exists a fixed point on $\partial \Sigma$ which is an attractor on $\Sigma$ (see Table \ref{biao0}), so it is impermanent. For the map $T$ in class $27$ such that $\varrho<0$, the conclusion follows from Corollary \ref{stability-permance}.
\end{proof}

\section{Applications to population models}
\label{sec:application}
In this section we apply the previous results in some concrete population models. Throughout this section, $A$ denotes the ${3\times 3}$ matrix with entries $a_{ij}>0$, $R=\diag[r_i]$ with $r_i>0$, and $U$ is the ${3\times 3}$ matrix with entries $\mu_{ij}>0$ such that $A=RU$, where $i,j=1,2,3$.
\subsection{Leslie-Gower model}
Consider the Leslie-Gower model \eqref{LG} due to Leslie and Gower \cite{leslie1958}.
The two-dimensional Leslie-Gower model is thoroughly
analyzed by Cushing \textit{et al.} \cite{cushing2004some}. The higher dimensional case was analyzed in \cite{hirsch2008existence,ruiz2011exclusion,jiang2015,LG}.

Denote the set of all Leslie-Gower maps \eqref{LG} by $\mathrm{CLG}(3)$. In symbols:
$$\mathrm{CLG}(3):=\{T\in \mathcal{T}(\mathbb{R}_+^3): T_i(x)=\frac{(1+r_i)x_i}{1+\sum_{j=1}^3a_{ij}x_j}, r_i>0, a_{ij}>0, i,j=1,2,3\}.$$
Set $f_i(z,r)=\frac{1+r}{1+z}$, $i=1,2,3$. Then the map $T=(T_1,T_2,T_3)$ with
$$
    T_i(x)=x_if_i((Ax^\tau)_i,r_i)=\frac{(1+r_i)x_i}{1+(RUx^\tau)_i}
$$
is just the Leslie-Gower model \eqref{LG}, i.e. $\mathrm{CLG}(3)$ is a special case of $\mathrm{DCS}(3,f)$.

Jiang and Niu \cite{LG} have listed the $33$ stable equivalence classes in $\mathrm{CLG}(3)$; see also Table \ref{biao0}. For $\mathrm{CLG}(3)$,
Proposition \ref{permanence-3D} is written in the following manner:
\begin{proposition}
The Leslie-Gower model $T\in\mathrm{CLG}(3)$ is permanent if it is in classes $29$, $31$, $33$ and class $27$ with $\varrho>0$ (defined by \eqref{stability-rho}), while $T$ is impermanent if it is in classes $1-26$, $28$, $30$, $32$ and class $27$ with $\varrho<0$.
\end{proposition}

In \cite{LG}, it was shown that  for $\mathrm{CLG}(3)$, Neimark-Sacker bifurcations can occur within each of classes $26-31$, so these classes can admit invariant closed curves. Here, we provide an example to show that the supercritical Neimark-Sacker bifurcation can occur in class {\rm 27} with $\varrho>0$ for $\mathrm{CLG}(3)$. We also provide a numerical example to show that the Chenciner (generalized Neimark-Sacker) bifurcation can occur in class {\rm 27} with $\varrho>0$, which implies that two isolated invariant closed curves can coexist on the carrying simplex in class {\rm 27} with a repelling heteroclinic cycle.
The Chenciner bifurcation is a two-parameter bifurcation phenomenon of a fixed point, which occurs when there is a pair of complex eigenvalues with modulus one and the first Lyapunov coefficient vanishes; see \cite{Kuznetsov,Govaerts2007Numerical} for more details.

\begin{figure}[h!]
    \centering
    \begin{tabular}{cc}
        \subfigure[The orbit with $x_0=(1,0.0667,0.0667)$]{
        \label{LG-27-1}
            \begin{minipage}[b]{0.46\textwidth}
                \centering                \includegraphics[width=\textwidth]{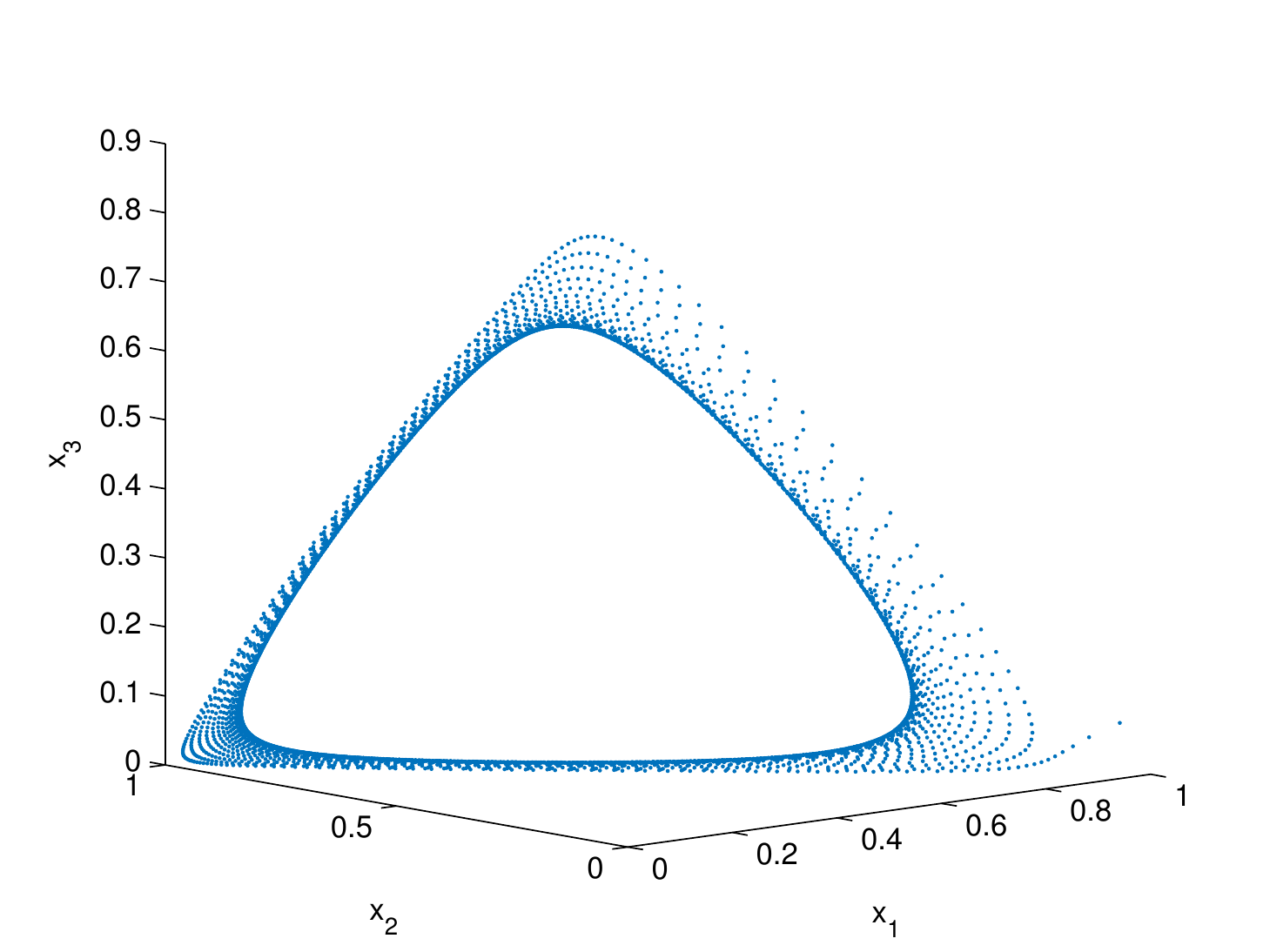}
            \end{minipage}
        } &
        \subfigure[The orbit with $x_0=(0.2151,0.746,0.0173)$]{\label{LG-27-2}
            \begin{minipage}[b]{0.46\textwidth}
                \centering
                \includegraphics[width=\textwidth]{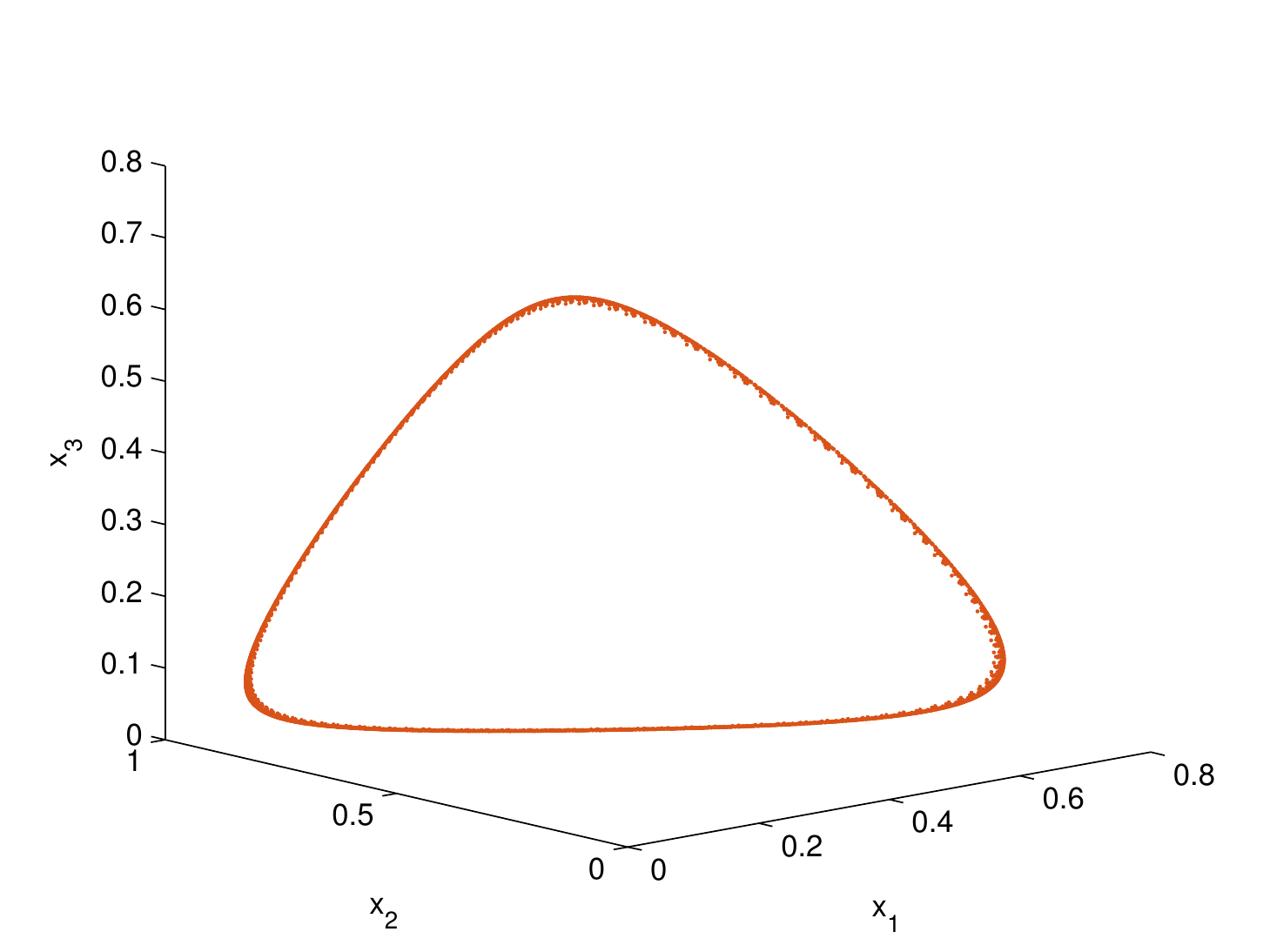}
            \end{minipage}
        }
    \end{tabular}
    \caption{The orbit emanating from $x_0=(1,0.0667,0.0667)$ for the map $T\in\mathrm{CLG}(3)$ with the parameter matrix $U$ given in Example \ref{exm:LG-1} and $r_1=1,r_2=0.2,r_3=1$ leads away from $\partial \Sigma$ and tends to an attracting invariant closed curve,  and the orbit emanating from $x_0=(0.2151,0.746,0.0173)$ also tends to an attracting invariant closed curve.} \label{fig:LG-27}
\end{figure}
\begin{example}\label{exm:LG-1}
Let $U=\left[ \begin {array}{ccc}
1&\frac{5}{4}&\frac{1}{2}\\
\noalign{\medskip}
\frac{1}{2}&1&\frac{3}{2}
\\ \noalign{\medskip}
\frac{3}{2}&\frac{3}{4}&1\end {array} \right]
$ and $r_1=1,r_2>0,r_3=1$. Consider the one-parameter family of maps $T^{[r_2]}\in\mathrm{CLG}(3)$ with the parameters $U$ and $r_i$. By Table \ref{biao0} {\rm(27)} we know that $T^{[r_2]}$ belongs to class {\rm 27} for all $r_2>0$. $T^{[r_2]}$ has a unique positive fixed point $p=(\frac{1}{4},\frac{1}{2},\frac{1}{4})$. When $r_2=r_2^*:=-{\frac{113}{194}}+{\frac {4\,\sqrt {295}}{97}}$, $DT^{[r_2]}(p)$ has a pair of complex conjugate eigenvalues of modulus $1$ which do not equal $\pm 1, \pm \mathrm{i}, (-1\pm \sqrt{3} \mathrm{i})/2$, where $\mathrm{i}$ stands for the imaginary unit. By calculating we obtain the first Lyapunov coefficient $l_1\approx -1.162\times 10^{-2}<0$. Since the Lyapunov coefficient is a rather lengthy expression, the approximate value was computed as a rational by using MATLAB \cite{Govaerts2007Numerical,Kuznetsov-Sacker,LG}. Therefore, a supercritical Neimark-Sacker bifurcation occurs at $r_2=r_2^*$, i.e., a stable invariant closed curve bifurcates from the fixed point $p$. On the other hand, it follows from \eqref{stability-rho} that $\varrho\approx  0.00078>0$ for $r_2=r_2^*$, so the heteroclinic cycle $\partial \Sigma$ of $T^{[r_2]}$ is repelling, i.e. $T^{[r_2]}$ is permanent, for any $r_2$ in a small neighborhood of $r_2^*$. Thus, a stable invariant closed curve can occur in class {\rm 27} with $\varrho>0$ for $\mathrm{CLG}(3)$. See Fig. \ref{fig:LG-27} for the orbit simulation.

Now let $r_1>0,r_2>0,r_3=1$, and consider the two-parameter family of maps $T^{[r_1,r_2]}\in\mathrm{CLG}(3)$ with the parameters $U$ and $r_i$. The map $T^{[r_1,r_2]}$ belongs to class {\rm 27} for all $r_1,r_2>0$ with a unique positive fixed point $p=(\frac{1}{4},\frac{1}{2},\frac{1}{4})$. By numerical calculation \cite{Govaerts2007Numerical,Govaerts2011A}, we find that $T^{[r_1,r_2]}$ admits a Chenciner bifurcation point at $p$ when $r_1\approx 0.248332$ and $r_2\approx 0.0633101$, where the second Lyapunov coefficient $l_2\approx -3.574\times 10^{-2}<0$.
Therefore, a stable fixed point and an attracting (large) invariant closed curve, separated by an unstable invariant closed curve can coexist in class $27$ for $\mathrm{CLG}(3)$ when the parameters $r_1$ and $r_2$ are properly disturbed near $0.248332$ and $0.0633101$ respectively; see \cite[Section 9.4]{Kuznetsov} or \cite[pp. 633--636]{GyllenbergCGAA} for details. Furthermore,  it follows from \eqref{stability-rho} that $\varrho\approx  0.00011>0$ for $r_1=0.248332$ and $r_2=0.0633101$, i.e. Chenciner bifurcation can also occur in class {\rm 27} with $\varrho>0$.
\end{example}

\subsection{Atkinson-Allen model}
Consider the generalized Atkinson-Allen model $T$ defined on $\mathbb{R}^3_+$ with
\begin{equation}\label{generalized Atkinson-Allen}
    T_i(x)=\frac{(1+r_i)(1-c_i)x_i}{1+\sum_{j=1}^3a_{ij}x_j}+c_ix_i, 0<c_i<1, a_{ij}, r_i>0, i,j=1,2,3.
\end{equation}
The model induced by the map \eqref{generalized Atkinson-Allen} is a discretized system of the competitive Lotka-Volterra equations, and see \cite{GyllenbergCGAA} for a mechanistic derivation of this model. A related two-dimensional discrete-time model for competition between populations of cyst-nematodes, due to Jones and Perry \cite{Jones1978}, was analyzed by Smith \cite{smith1998planar}.  When $r_i=1$ and $c_i=c$, the map \eqref{generalized Atkinson-Allen} reduces to the standard Atkinson-Allen map
\begin{equation}\label{model:AA}
    T: \mathbb{R}^3_+ \mapsto \mathbb{R}^3_+, \  T_i(x)=\frac{2(1-c)x_i}{1+\sum_{j=1}^3a_{ij}x_j}+cx_i,\  0<c<1, a_{ij}>0, i,j=1,2,3,
\end{equation}
which is a modified model derived from annual plants competition \cite{atkinson1997,Allen1998,roeger2004discrete}, and has been analyzed by Jiang and Niu in \cite{jiang2014}.

Since map \eqref{model:AA}
is a special case of the generalized Atkinson-Allen map \eqref{generalized Atkinson-Allen}, we apply the previous results to the map \eqref{generalized Atkinson-Allen}.
Denote the set of all generalized Atkinson-Allen maps \eqref{generalized Atkinson-Allen} by $\mathrm{CGAA}(3)$. In symbols:
$$
\mathrm{CGAA}(3):=\{T\in \mathcal{T}(\mathbb{R}_+^3): T_i(x)=\frac{(1+r_i)(1-c_i)x_i}{1+\sum_{j=1}^3a_{ij}x_j}+c_ix_i, 0<c_i<1, a_{ij}, r_i>0\}.
$$
Set $f_i(z,r)=\frac{(1+r)(1-c_i)}{1+z}+c_i$, $0<c_i<1$, $i=1,2,3$. Then the map $T=(T_1,T_2,T_3)$ with
$$
    T_i(x)=x_if_i((Ax^\tau)_i,r_i)=\frac{(1+r_i)(1-c_i)x_i}{1+(RUx^\tau)_i}+c_ix_i,\quad i=1,2,3
$$
is the generalized Atkinson-Allen model, i.e. $\mathrm{CGAA}(3)$ is also a special case of $\mathrm{DCS}(3,f)$.

Gyllenberg et al. \cite{GyllenbergRicker} have listed the $33$ stable equivalence classes in $\mathrm{CGAA}(3)$; see also Table \ref{biao0}. For $\mathrm{CGAA}(3)$,
Proposition \ref{permanence-3D} is written in the following manner:
\begin{proposition}
The generalized Atkinson-Allen model $T\in\mathrm{CGAA}(3)$ is permanent if it is in classes $29$, $31$, $33$ and class $27$ with $\varrho>0$ (defined by \eqref{stability-rho}), while $T$ is impermanent if it is in classes $1-26$, $28$, $30$, $32$ and class $27$ with $\varrho<0$.
\end{proposition}

It was shown in \cite{GyllenbergRicker} that for $\mathrm{CGAA}(3)$, classes $26-29$ and $31$ can admit supercritical Neimark-Sacker bifurcations, and class $30$ can admit subcritical Neimark-Sacker bifurcations. The authors also numerically show that Chenciner bifurcations can occur in classes $26-29$. Here, we give two examples to show that the supercritical Neimark-Sacker bifurcation can occur in class $27$ with $\varrho>0$ and can also occur in class $27$ with $\varrho<0$, respectively.

\begin{figure}[h!]
    \centering
    \begin{tabular}{cc}
        \subfigure[The orbit with $x_0=(1,0.0667,0.0667)$]{
        \label{GAA-27-1}
            \begin{minipage}[b]{0.45\textwidth}
                \centering                \includegraphics[width=\textwidth]{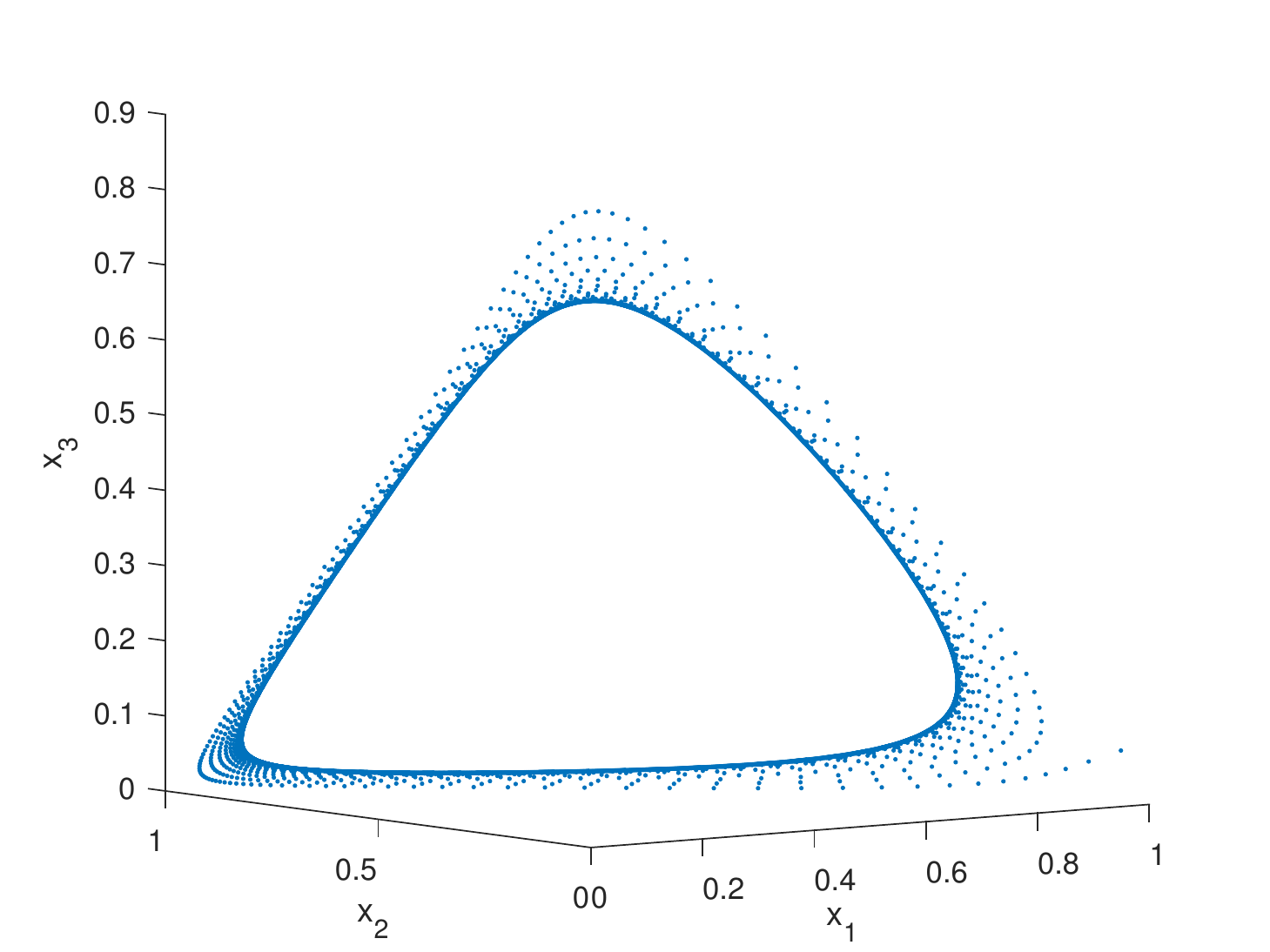}
            \end{minipage}
        } &
        \subfigure[The orbit with $x_0= (0.7,0.1642,0.1685)$]{\label{GAA-27-2}
            \begin{minipage}[b]{0.45\textwidth}
                \centering
                \includegraphics[width=\textwidth]{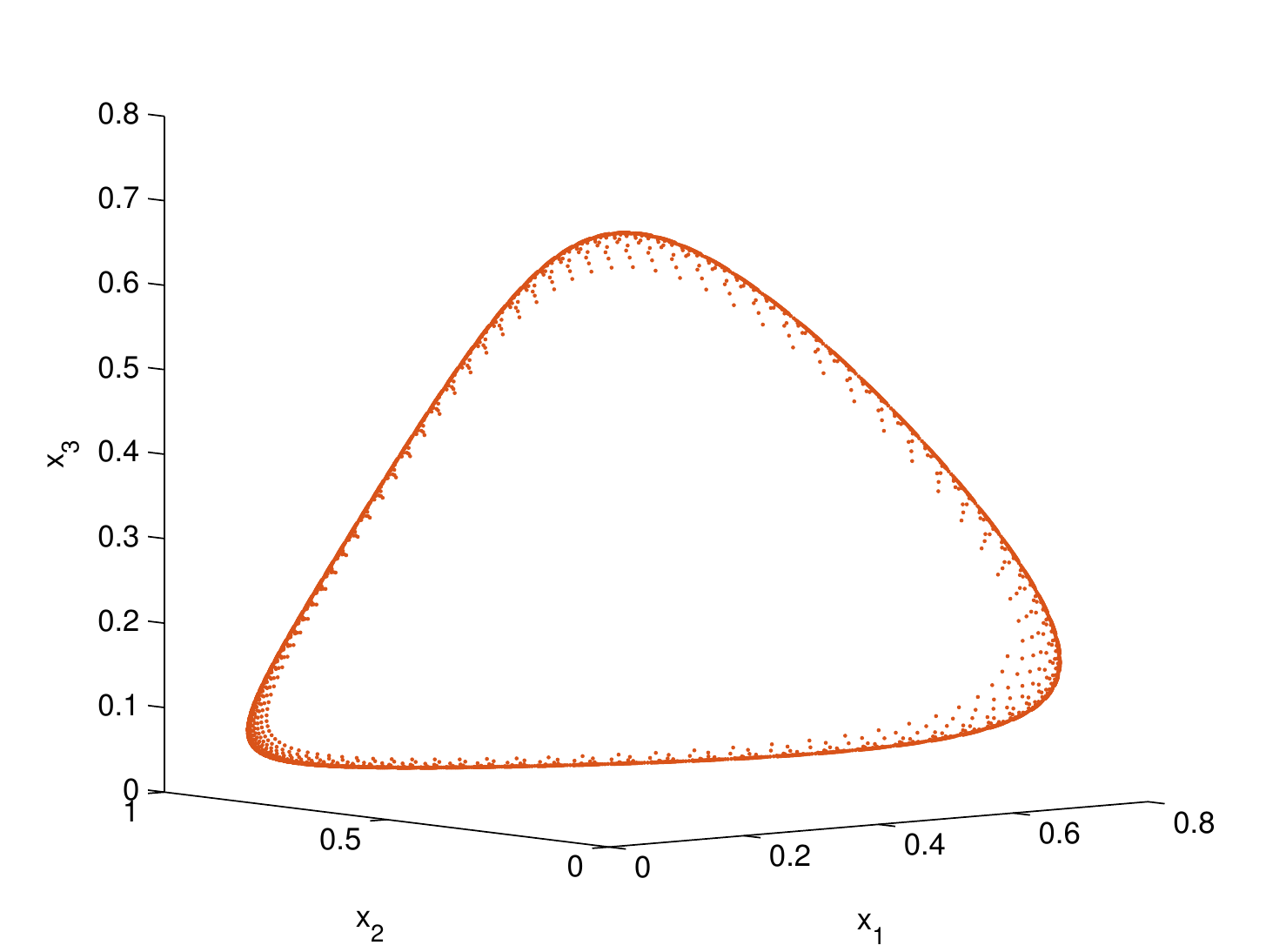}
            \end{minipage}
        }
    \end{tabular}
    \caption{The orbit emanating from $x_0=(1,0.0667,0.0667)$ for the map $T\in\mathrm{CGAA}(3)$ with the parameter matrix $U$ given in Example \ref{exm:LG-1} and $r_1=r_2=r_3=1$, $c_1=\frac{1}{10},c_2=\frac{1}{5}, c_3=\frac{1}{5}$ leads away from $\partial \Sigma$ and tends to an attracting invariant closed curve,  and the orbit emanating from $x_0=(0.7,0.1642,0.1685)$ also tends to an attracting invariant closed curve.} \label{fig:GAA-27-1}
\end{figure}
\begin{example}\label{exm:CGAA-1}
Let $r_1=r_2=r_3=1$, and $c_1=\frac{1}{10},c_2=\frac{1}{5}, 0<c_3<1$. Consider the one-parameter family of maps $T^{[c_3]}\in\mathrm{CGAA}(3)$ with the parameter matrix $U$ given in Example \ref{exm:LG-1} and the above $r_i$, $c_i$, $i=1,2,3$. By Table \ref{biao0} {\rm(27)} we know that $T^{[c_3]}$ belongs to class {\rm 27} for all $0<c_3<1$, whose unique positive fixed point $p=(\frac{1}{4},\frac{1}{2},\frac{1}{4})$. When $c_3=c_3^*:={\frac{432709}{80801}}-{\frac {80\,\sqrt {24656689}}{80801}}$, $DT^{[c_3]}(p)$ has a pair of complex conjugate eigenvalues of modulus $1$ which do not equal $\pm 1, \pm \mathrm{i}, (-1\pm \sqrt{3} \mathrm{i})/2$. By numerical calculation \cite{Govaerts2007Numerical,Kuznetsov-Sacker,GyllenbergCGAA}, we get the first Lyapunov coefficient $l_1\approx -1.814\times 10^{-2}<0$. Therefore, a supercritical Neimark-Sacker bifurcation occurs at $c_3=c_3^*$, i.e., a stable invariant closed curve bifurcates from the fixed point $p$. On the other hand, it follows from \eqref{stability-rho} that $\varrho\approx  0.0026>0$ for $c_3=c_3^*$, so the heteroclinic cycle $\partial \Sigma$ of $T^{[c_3]}$ is repelling, i.e. $T^{[c_3]}$ is permanent, for any $c_3$ in a small neighborhood of $c_3^*$. Thus, a stable invariant closed curve can occur in class {\rm 27} with $\varrho>0$ for $\mathrm{CGAA}(3)$. See Fig. \ref{fig:GAA-27-1} for the orbit simulation.
\end{example}

\begin{figure}[h!]
    \centering
    \begin{tabular}{cc}
        \subfigure[The orbit with $x_0=(0.04,0.12,0.36)$]{
        \label{GAA-27-3}
            \begin{minipage}[b]{0.45\textwidth}
                \centering                \includegraphics[width=\textwidth]{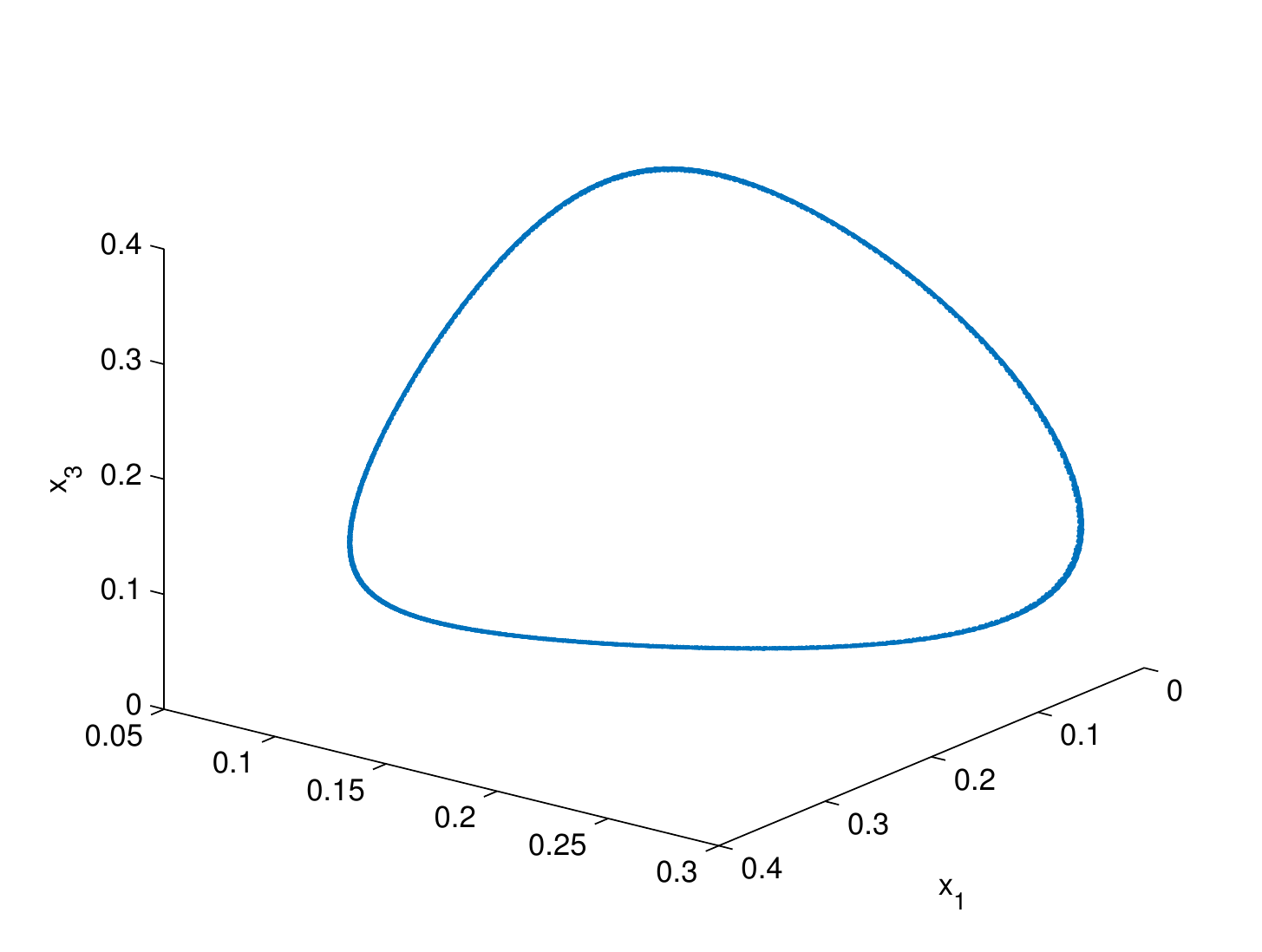}
            \end{minipage}
        } &
        \subfigure[The orbit with $x_0=(0.0002,0.023,0.486)$]{\label{GAA-27-4}
            \begin{minipage}[b]{0.45\textwidth}
                \centering
                \includegraphics[width=\textwidth]{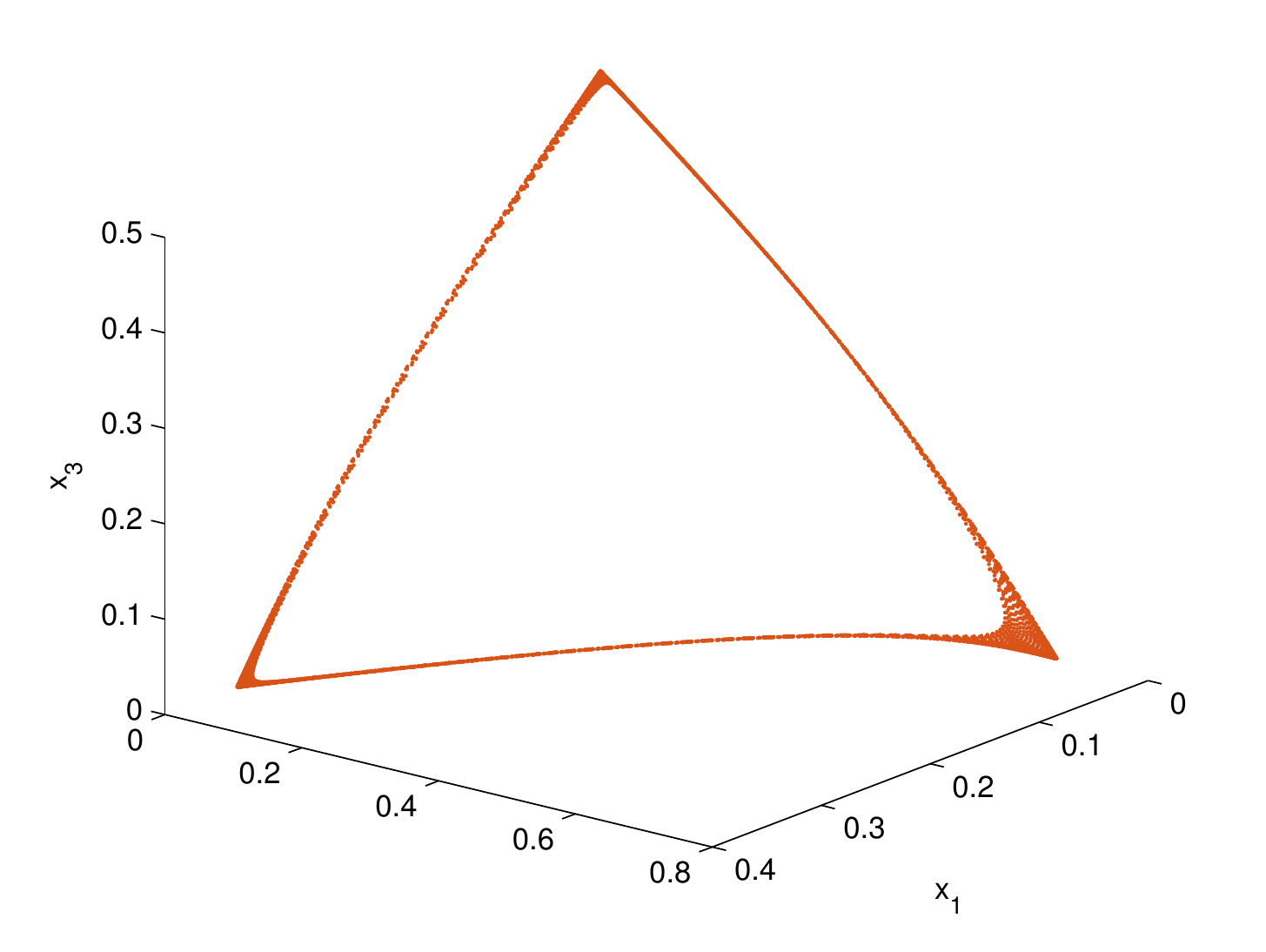}
            \end{minipage}
        }
    \end{tabular}
    \caption{The orbit emanating from $x_0=(0.04,0.12,0.36)$ for the map $T\in\mathrm{CGAA}(3)$ with the parameter matrix $U$ given in Example \ref{exm:CGAA-2} and $r_1=r_2=r_3=1$, $c_1=0.1,c_2=0.79, c_3=0.1$ tends to an attracting invariant closed curve, while the orbit emanating from $x_0=(0.0002,0.023,0.486)$ approaches the heteroclinic cycle $\partial \Sigma$.} \label{fig:GAA-27-2}
\end{figure}

\begin{example}\label{exm:CGAA-2}
Let $U= \left[
\begin {array}{ccc}
3&3&1\\
\noalign{\medskip}
\frac{3}{2}&\frac{3}{2}&4
\\
\noalign{\medskip}
4&1&2\end {array} \right]
$, and $r_1=r_2=r_3=1$, $c_1=c_3=\frac{1}{10},0<c_2<1$. Consider the one-parameter family of maps $T^{[c_2]}\in\mathrm{CGAA}(3)$ with the parameters $U$ and $r_i$, $c_i$, $i=1,2,3$. By Table \ref{biao0} {\rm(27)} we know that $T^{[c_2]}$ belongs to class {\rm 27} for all $0<c_2<1$, whose unique positive fixed point $p=(\frac{1}{7},\frac{1}{7},\frac{1}{7})$. When $c_2=c_2^*:={\frac{1822387}{382723}}-{\frac {840\,\sqrt {3257017}}{382723}}$, $DT^{[c_2]}(p)$ has a pair of complex conjugate eigenvalues of modulus $1$ which do not equal $\pm 1, \pm \mathrm{i}, (-1\pm \sqrt{3} \mathrm{i})/2$. By numerical calculation, we obtain the first Lyapunov coefficient $l_1\approx -5.039\times 10^{-2}<0$. Therefore, a supercritical Neimark-Sacker bifurcation occurs at $c_2=c_2^*$, and hence a stable invariant closed curve bifurcates from the fixed point $p$. On the other hand, it follows from \eqref{stability-rho} that $\varrho\approx -0.00058<0$ for $c_2=c_2^*$, so the heteroclinic cycle $\partial \Sigma$ of $T^{[c_2]}$ is attracting, i.e. $T^{[c_2]}$ is impermanent, for any $c_2$ in a small neighborhood of $c_2^*$. Thus, the supercritical Neimark-Sacker can occur in class {\rm 27} with $\varrho<0$ for $\mathrm{CGAA}(3)$. See Fig. \ref{fig:GAA-27-2} for the orbit simulation.
\end{example}

\subsection{Mixing growth functions}
Consider the following model $T=(T_1,T_2,T_3)$ on $\mathbb{R}_+^3$, in which the three competing species are assumed to have different types of growth functions:
\begin{equation}\label{map:mix}
\left\{
    \begin{array}{l}
    T_1(x)=\displaystyle\frac{(1+r_1)x_1}{1+a_{11}x_1+a_{12}x_2+a_{13}x_3},\\
    \noalign{\medskip}
    T_2(x)=\displaystyle \frac{(1+r_2)(1-c)x_2}{1+a_{21}x_1+a_{22}x_2+a_{23}x_3}+cx_2, \\
     \noalign{\medskip}
    T_3(x)=\displaystyle \frac{(1+\ln(1+r_3))x_3}{1+\ln (1+a_{31}x_1+a_{32}x_2+a_{33}x_3)}.
    \end{array}
    \right.
\end{equation}
Set $f_1(z,r)=\frac{1+r}{1+z}$, $ f_2(z,r)=\frac{(1+r)(1-c)}{1+z}+c$, $0<c<1$, $f_3(z,r)=\frac{1+\ln(1+r)}{1+\ln (1+z)}$. Then the map $T$ can be written as
$$
    T_i(x)=x_if_i((RUx^\tau)_i,r_i),\quad i=1,2,3.
$$
Note that each $f_i\in \mathscr{F}$, so $T$ admits a carrying simplex $\Sigma$.
Denote the set of all maps \eqref{map:mix} by
$$
\mathrm{MGF}(3):=\{T\in \mathcal{T}(\mathbb{R}_+^3): T_i(x)=x_if_i((RUx^\tau)_i,r_i),\, \mu_{ij}, r_i>0\}.
$$
Therefore, $\mathrm{MGF}(3)$ is a special case of $\mathrm{DCS}(3,f)$ with the generating function $f=(f_1,f_2,f_3)$. It follows from Theorem \ref{theory:classification} that there are $33$ stable equivalence classes in $\mathrm{MGF}(3)$, and furthermore, Proposition \ref{permanence-3D} is written in the following manner:
\begin{proposition}
The model $T\in\mathrm{MGF}(3)$ is permanent if it is in classes $29$, $31$, $33$ and class $27$ with $\varrho>0$ (defined by \eqref{stability-rho}), while $T$ is impermanent if it is in classes $1-26$, $28$, $30$, $32$ and class $27$ with $\varrho<0$.
\end{proposition}

We now provide two examples to show that Neimark-Sacker bifurcations can occur in the permanent classes $29$ and $31$ for $\mathrm{MGF}(3)$, respectively.
\begin{figure}[h!]
    \centering
    \begin{tabular}{cc}
        \subfigure[The orbit emanating from $x_0=(0.427,0.8574,0.014)$]{
            \begin{minipage}[b]{0.42\textwidth}
                \centering                \includegraphics[width=\textwidth]{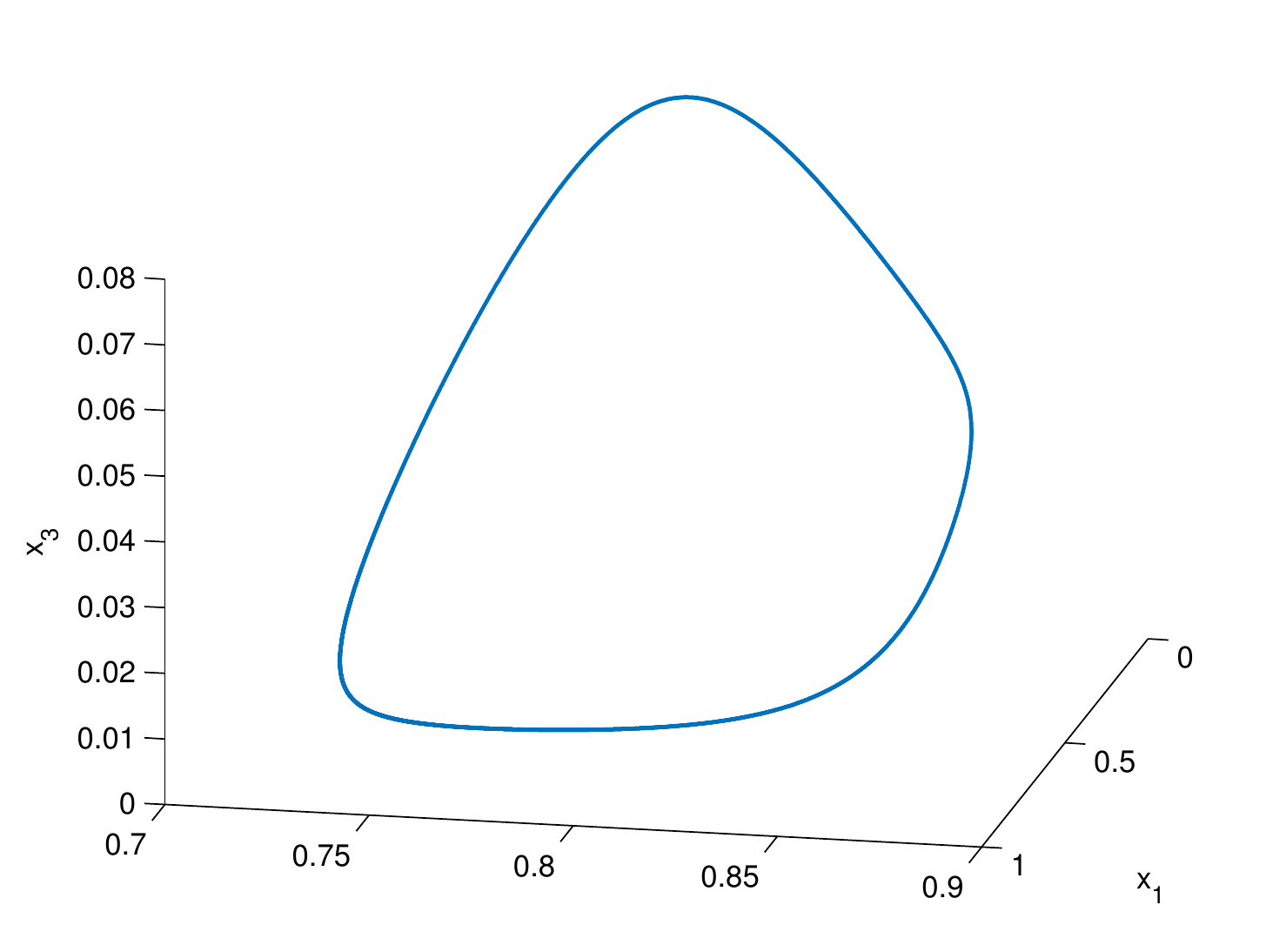}
            \end{minipage}
        } &
        \subfigure[The motion of components: $x_1$ (blue), $x_2$ (red) and $x_3$ (green)]{
            \begin{minipage}[b]{0.48\textwidth}
                \centering
                \includegraphics[width=\textwidth]{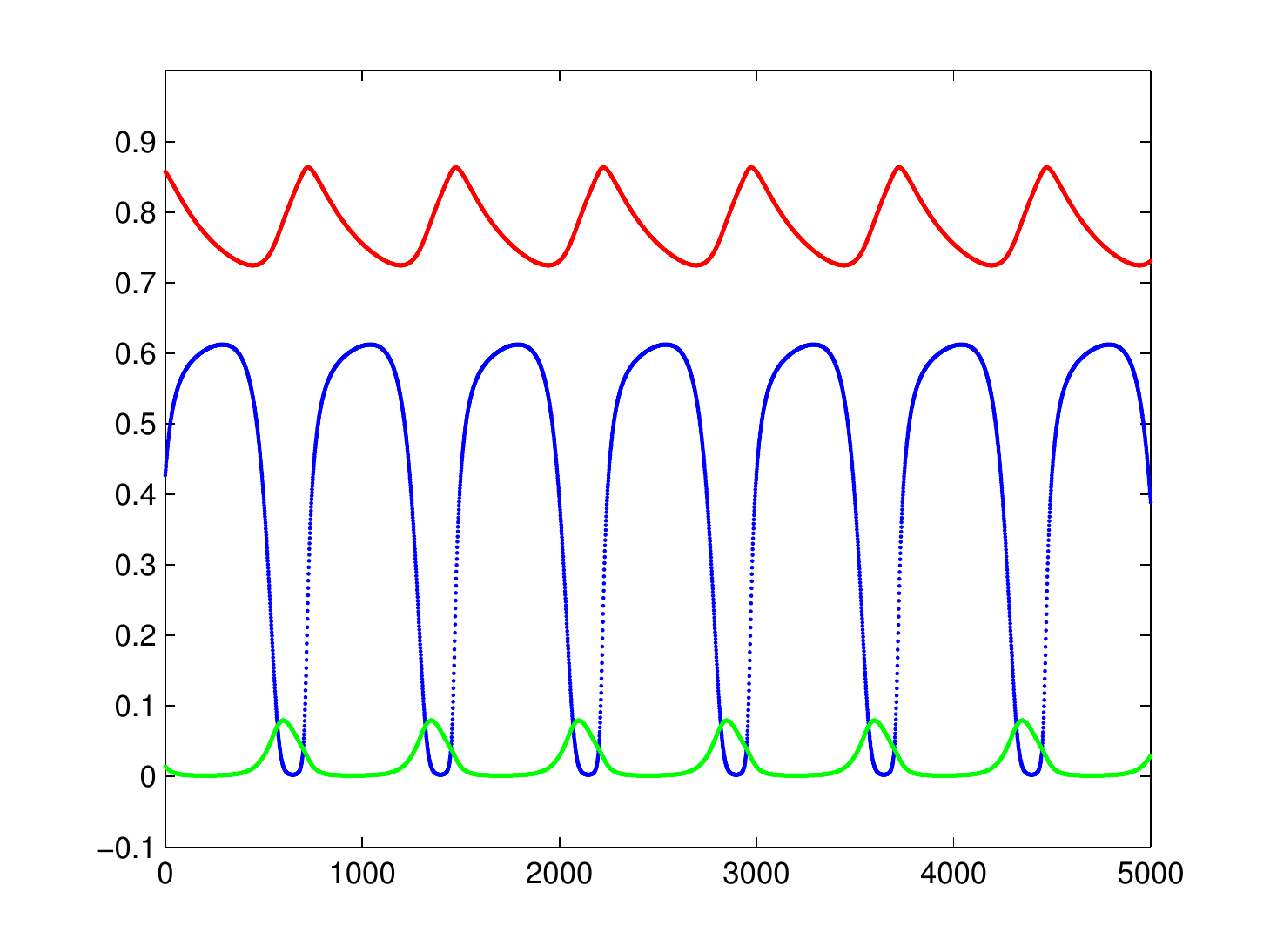}
            \end{minipage}
        }
    \end{tabular}
    \caption{The orbit emanating from $x_0=(0.427,0.8574,0.014)$ for the map $T\in\mathrm{MFC}(3)$ with the parameter matrix $U$ given in Example \ref{exm:Mix-1}, $c=\frac{4}{5}$ and $r_1=r_3=1,r_2=0.03$ tends to an attracting invariant closed curve.} \label{fig:Mix-29}
\end{figure}

\begin{example}\label{exm:Mix-1}
Let $U=\left[ \begin {array}{ccc}
1&\frac{1}{2}&9\\
\noalign{\medskip}
\frac{1}{2}&1&\frac{1}{2}\\
\noalign{\medskip}
\frac{1}{6}&\frac{7}{6}&1
\end {array} \right]
$ and $c=\frac{4}{5}$, $r_1=r_3=1,r_2>0$. Consider the one-parameter family of maps $T^{[r_2]}\in \mathrm{MGF}(3)$ with the parameters $U$, $c$ and $r_i$, $i=1,2,3$. By Table \ref{biao0} {\rm(29)} we known that $T^{[r_2]}$ belongs to class {\rm 29} with a unique positive fixed point $p=(\frac{8}{19},\frac{74}{95},\frac{2}{95})$ for all $r_2>0$. When $r_2\approx 0.032889$, $DT^{[r_2]}(p)$ has a pair of complex conjugate eigenvalues with modulus $1$ which do not equal $\pm 1, \pm \mathrm{i}, (-1\pm \sqrt{3} \mathrm{i})/2$. The first Lyapunov coefficient $l_1\approx -2.430\times 10^{-5}<0$.
Therefore, there is a supercritical Neimark-Sacker bifurcation in class {\rm 29} for $\mathrm{MGF}(3)$, i.e. a stable invariant closed curve bifurcates from the fixed point $p$. See Fig. \ref{fig:Mix-29} for the orbit simulation.
\end{example}

\begin{figure}[h!]
    \centering
    \begin{tabular}{cc}
        \subfigure[The orbit emanating from $x_0=(0.5962,0.4857,0.193)$]{
            \begin{minipage}[b]{0.42\textwidth}
                \centering                \includegraphics[width=\textwidth]{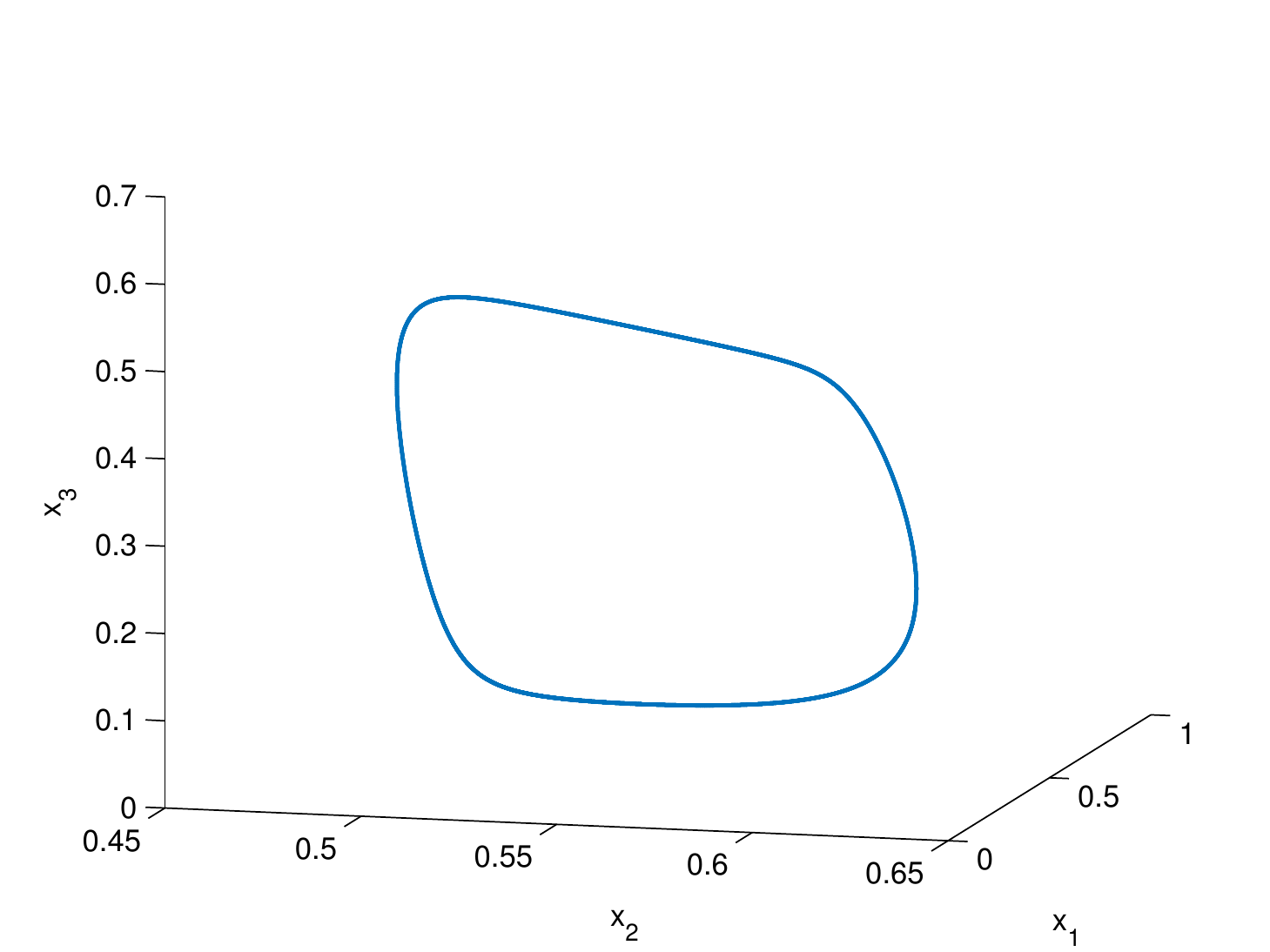}
            \end{minipage}
        } &
        \subfigure[The motion of components: $x_1$ (blue), $x_2$ (red) and $x_3$ (green)]{
            \begin{minipage}[b]{0.48\textwidth}
                \centering
                \includegraphics[width=\textwidth]{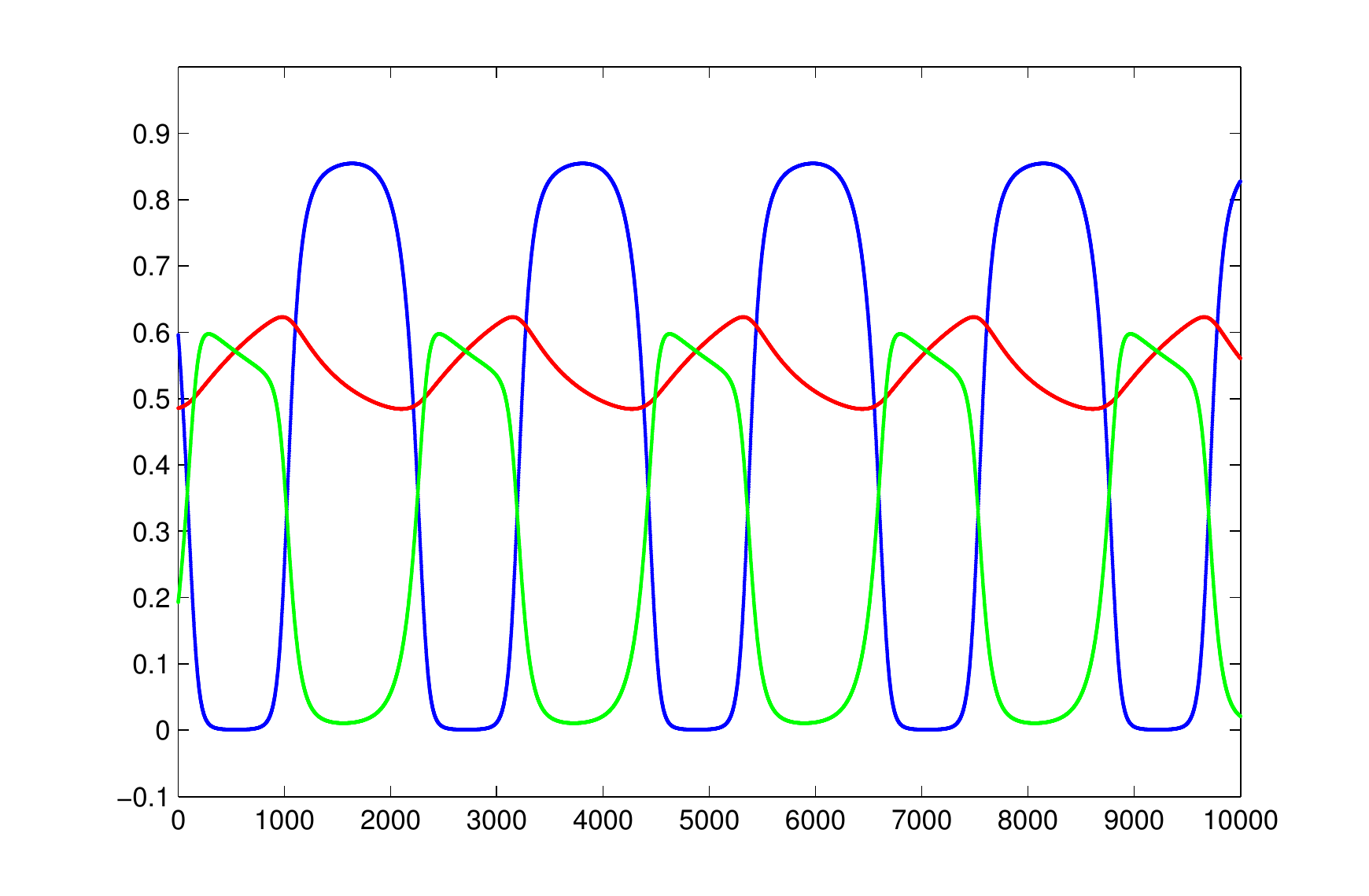}
            \end{minipage}
        }
    \end{tabular}
    \caption{The orbit emanating from $x_0=(0.5962,0.4857,0.193)$ for the map $T\in\mathrm{MFC}(3)$ with the parameter matrix $U$ given in Example \ref{exm:Mix-2}, $c=\frac{4}{5}$ and $r_1=r_3=1,r_2=0.02$ tends to an attracting invariant closed curve.} \label{fig:Mix-31}
\end{figure}

\begin{example}\label{exm:Mix-2}
Let $U=\left[ \begin {array}{ccc} 1&\frac{1}{4}&\frac{3}{2}\\
\noalign{\medskip}
\frac{5}{8}&1&\frac{5}{8}\\
\noalign{\medskip}
{\frac{7}{10}}&\frac{3}{4}&1\end {array}
\right]
$ and $c=\frac{4}{5}$, $r_1=r_3=1,r_2>0$. Consider the one-parameter family of maps $T^{[r_2]}\in \mathrm{MGF}(3)$ with the parameters $U$, $c$ and $r_i$, $i=1,2,3$. By Table \ref{biao0} {\rm(31)} we known that $T^{[r_2]}$ belongs to class {\rm 31} with a unique positive fixed point $p=(\frac{5}{11},\frac{6}{11},\frac{3}{11})$ for all $r_2>0$. When $r_2\approx 0.038917$, $DT^{[r_2]}(p)$ has a pair of complex conjugate eigenvalues with modulus $1$ which do not equal $\pm 1, \pm \mathrm{i}, (-1\pm \sqrt{3} \mathrm{i})/2$. The first Lyapunov coefficient $l_1\approx -3.968\times 10^{-3}<0$.
Therefore, there is a supercritical Neimark-Sacker bifurcation in class {\rm 31} for $\mathrm{MGF}(3)$, i.e. a stable invariant closed curve bifurcates from the fixed point $p$. See Fig. \ref{fig:Mix-31} for the orbit simulation.
\end{example}

\subsection{Ricker model}\label{subsection:Ricker}
Consider the Ricker map \cite{Ricker1954}
\begin{equation}\label{equ:Ricker}
T: \mathbb{R}^3_+ \mapsto \mathbb{R}^3_+,\ T_i(x)=x_i \exp(r_i-\sum_{j=1}^3a_{ij}x_j),~ r_i,a_{ij}>0, i,j=1,2,3.
\end{equation}
The one-dimensional map has been studied in detail by May and Oster \cite{May1976}, where they showed that every orbit converges to the positive fixed point for $r\leq 2$, and it will exhibit a scenario of chaotic behavior for large $r$. The two-dimensional map was analyzed in detail by Smith \cite{smith1998planar}, who showed that it has trivial dynamics provided $\nu_1,\nu_2<1$. Roeger \cite{Lih2005} studied the local dynamics of the positive fixed point and Neimark-Sacker bifurcations for the map \eqref{equ:Ricker} with $r_1=r_2=r_3$. Hofbauer et al. \cite{hofbauer1987coexistence} provided the criteria on permanence for map \eqref{equ:Ricker} and also the higher dimensional cases.

Set $f_i(z,r)=\exp(r-z)$, $i=1,2,3$. Then the Ricker map \eqref{equ:Ricker} can be written as
\begin{equation}\label{equ:Ricker2}
T_i(x)=x_if_i((Ax^\tau)_i,r_i)=x_i \exp(r_i(1-\sum_{j=1}^3\mu_{ij}x_j)),\quad i=1,2,3.
\end{equation}
Note that \eqref{cons:f} (ii) does not hold for $f_i$, that is $f_i\notin \mathscr{F}$, and unlike the maps in $\mathrm{DCS}(3,f)$ (such as the Leslie-Gower map or the Atkinson-Allen map discussed above), the Ricker map $T$ has a carrying simplex only under certain additional conditions (see \cite{GyllenbergRicker}).
Assume that the parameters satisfy
\begin{equation}\label{Ricker-con-1}
  r_i<1/(\sum_{j=1}^3\frac{\mu_{ij}}{\mu_{jj}}),\mathrm{~or~} r_i<\mu_{ii}/\sum_{j=1}^3\mu_{ij},~i=1,2,3.
\end{equation}
Then one can easily check that the Ricker map \eqref{equ:Ricker} satisfies the condition $\Upsilon$3) in Lemma \ref{simplex} and hence it admits a carrying simplex by Lemma \ref{simplex}.

Denote by
$$
\mathrm{CRC}(3):=\{T\in \mathcal{T}(\mathbb{R}_+^3): T_i(x)=x_i \exp(r_i(1-\sum_{j=1}^3\mu_{ij}x_j)), r_i,\mu_{ij}>0,  \eqref{Ricker-con-1} \mathrm{~holds}\}.
$$
Then each Ricker map \eqref{equ:Ricker2} in $\mathrm{CRC}(3)$ admits a carrying simplex. The classification program  via the dynamics on $\partial\Sigma$ and statements for the $3$-dimensional maps \eqref{equ:T1} are also applicable for $\mathrm{CRC}(3)$. Specifically, Gyllenberg et al. showed in \cite{GyllenbergRicker} that there are a total of $33$ stable equivalence classes in $\mathrm{CRC}(3)$ as shown in Table \ref{biao0}, where the parameters should satisfy the condition \eqref{Ricker-con-1} in addition to those listed in Table \ref{biao0} for each class. On the other hand, note that all the criteria on the permanence in Section \ref{subsection:permanence} do not depend on the condition \eqref{cons:f} (ii) and is applicable to any map $T$ given by \eqref{equ:T2} which has a carrying simplex. Moreover, by the proof of Proposition \ref{permanence-3D}, one can see that the existence of the carrying simplex and conditions (i) and (ii) in Table \ref{biao0} (29) imply the permanence of the class $29$, and similarly for classes $31$ and $33$, etc. Therefore, for the Ricker map \eqref{equ:Ricker2}, Proposition \ref{permanence-3D} is written in the following manner:
\begin{proposition}
The Ricker model $T\in\mathrm{CRC}(3)$ is permanent if it is in classes $29$, $31$, $33$ and class $27$ with $\varrho>0$ $($defined by \eqref{stability-rho}$)$, while $T$ is impermanent if it is in classes $1-26$, $28$, $30$, $32$ and class $27$ with $\varrho<0$.
\end{proposition}

It was shown in \cite{GyllenbergRicker} that for $\mathrm{CRC}(3)$, classes $26$ and $31$ can admit supercritical Neimark-Sacker bifurcations, while classes $27-30$ can admit subcritical Neimark-Sacker bifurcations. The authors also provided a numerical example to show that the Chenciner bifurcation can occur in class $29$. Here, we give an example to show that the supercritical Neimark-Sacker bifurcation can also occur in class $29$, and a numerical example to show that the Chenciner bifurcation can also occur in class $26$.
\begin{figure}[h!]
    \centering
    \begin{tabular}{cc}
        \subfigure[The orbit emanating from $x_0=(0.3128,0.8347,0.0199)$]{
            \begin{minipage}[b]{0.42\textwidth}
                \centering                \includegraphics[width=\textwidth]{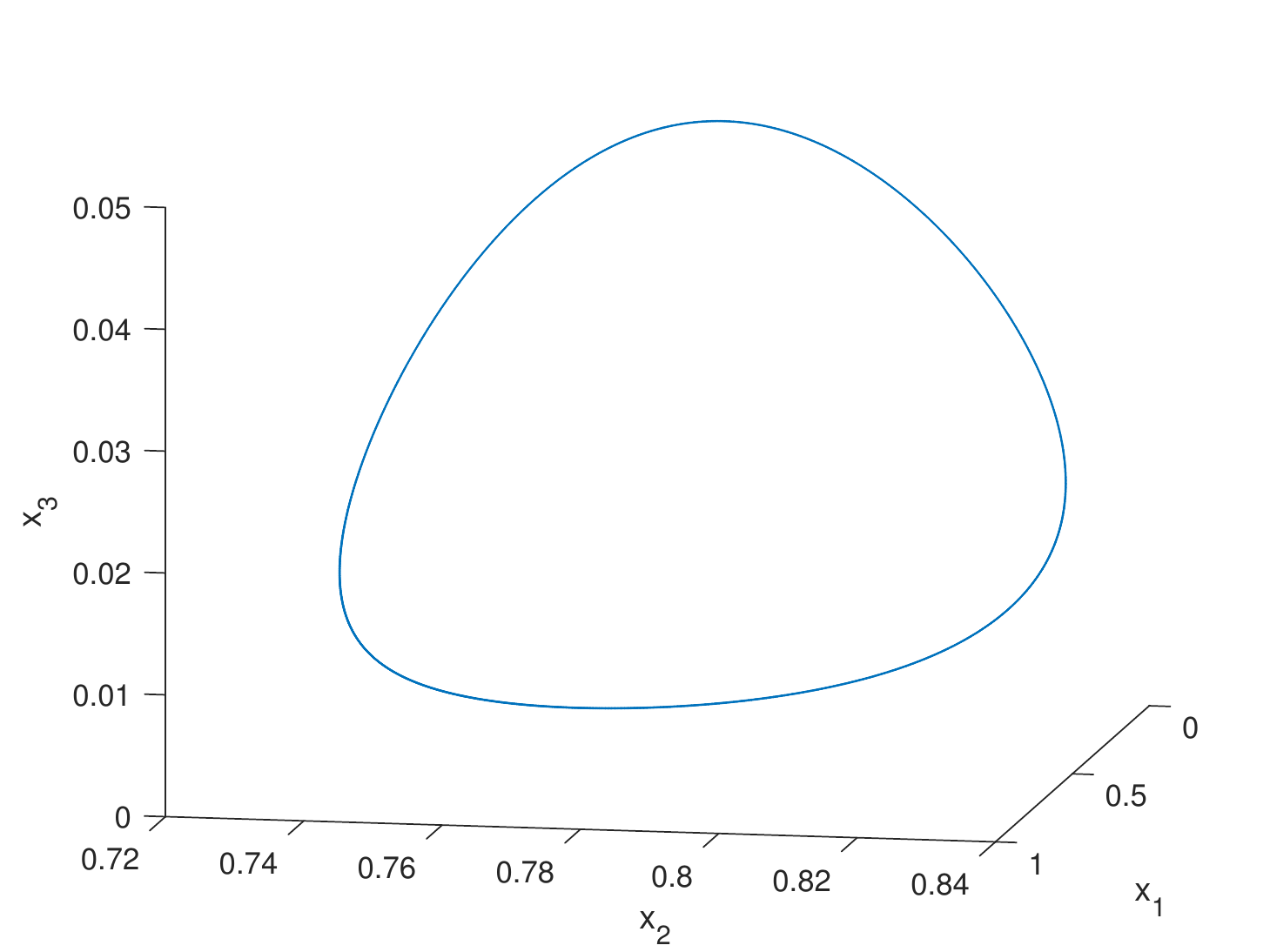}
            \end{minipage}
        } &
        \subfigure[The motion of components: $x_1$ (blue), $x_2$ (red) and $x_3$ (green)]{
            \begin{minipage}[b]{0.48\textwidth}
                \centering
                \includegraphics[width=\textwidth]{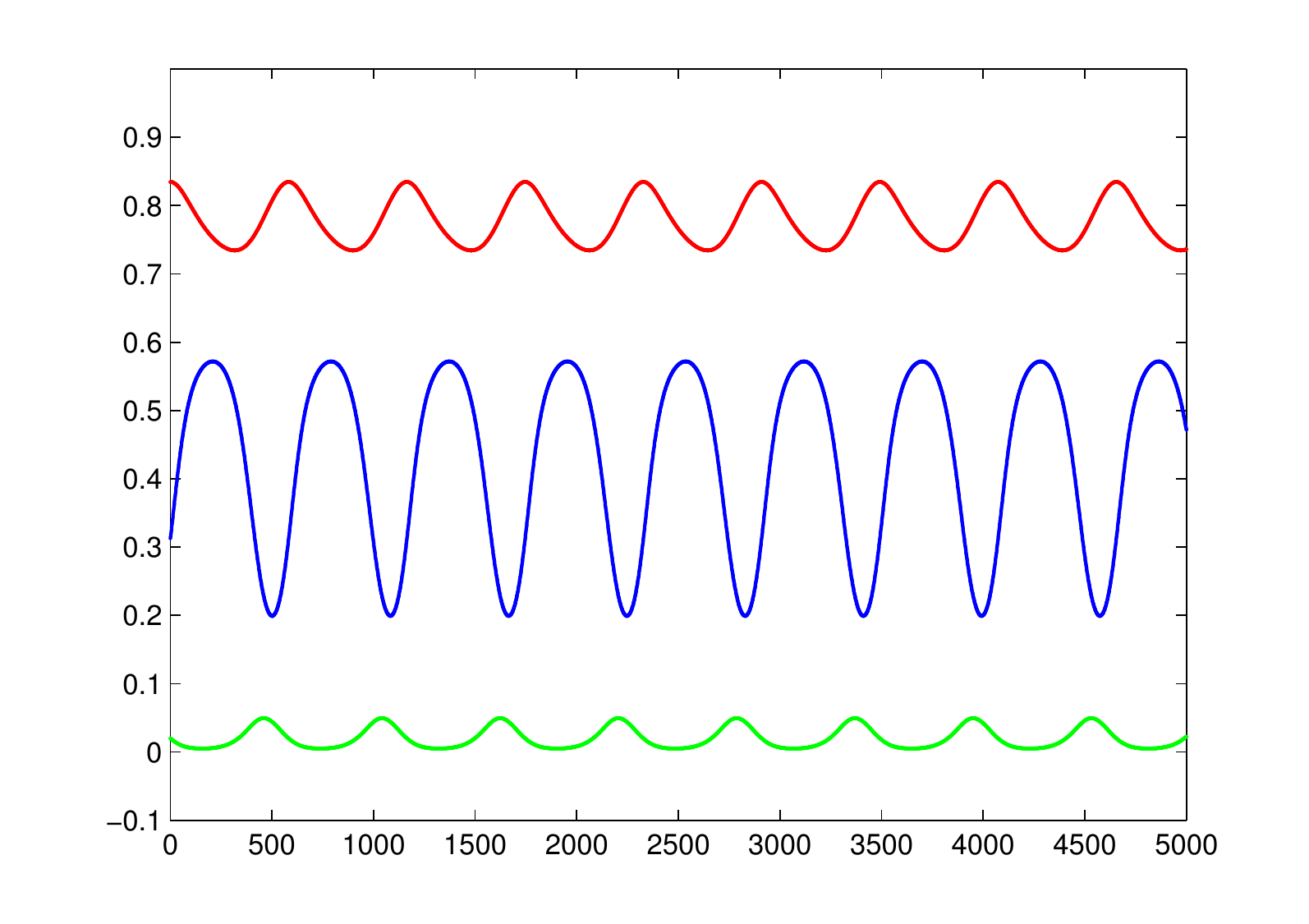}
            \end{minipage}
        }
    \end{tabular}
    \caption{The orbit emanating from $x_0=(0.3128,0.8347,0.0199)$ for the map $T\in\mathrm{CRC}(3)$ with the parameter matrix $U$ given in Example \ref{exm:Mix-1} and $r_1=\frac{1}{11},r_2=0.01,r_3=\frac{2}{7}$ tends to an attracting invariant closed curve.} \label{fig:Ricker-29}
\end{figure}

\begin{example}\label{exm:Ricker-2}
Let $r_1=\frac{1}{11},0<r_2<\frac{1}{2},r_3=\frac{2}{7}$. Consider the one-parameter family of maps $T^{[r_2]}$ given by \eqref{equ:Ricker2} with the parameter matrix $U$ given in Example \ref{exm:Mix-1} and the above $r_i$, $i=1,2,3$. It is easy to check that such $\mu_{ij},r_i$ satisfy \eqref{Ricker-con-1}, i.e. $T^{[r_2]}\in \mathrm{CRC}(3)$ for all $0<r_2<\frac{1}{2}$. It follows from Table \ref{biao0} {\rm(29)} that $T^{[r_2]}$ belongs to class {\rm 29}. $T^{[r_2]}$ possesses a unique positive fixed point $p=(\frac{8}{19},\frac{74}{95},\frac{2}{95})$. When $r_2=-{\frac{15}{2128}}+{\frac {3\,\sqrt {231729}}{78736}}$, $DT^{[r_2]}(p)$ has a pair of complex conjugate eigenvalues with modulus $1$ which do not equal $\pm 1, \pm \mathrm{i}, (-1\pm \sqrt{3} \mathrm{i})/2$. The first Lyapunov coefficient $l_1=-1.433\times 10^{-2}<0$.
Therefore, there is a supercritical Neimark-Sacker bifurcation in class {\rm 29} for $\mathrm{CRC}(3)$, i.e. a stable invariant closed curve bifurcates from the fixed point $p$. See Fig. \ref{fig:Ricker-29} for the orbit simulation.
\end{example}

\begin{example}
Let $U=\left[ \begin {array}{ccc} 1&4&\frac{3}{4}\\ \noalign{\medskip}\frac{1}{8}&1&\frac{5}{4}
\\ \noalign{\medskip}\frac{3}{4}&\frac{5}{4}&1\end {array} \right]
$ and $0<r_1<\frac{1}{6},r_2=\frac{1}{5},0<r_3<\frac{1}{4}$. Consider the two-parameter family of maps $T^{[r_1,r_3]}$ given by \eqref{equ:Ricker2} with the parameters $U$ and $r_i$, $i=1,2,3$. It is easy to check that such $\mu_{ij},r_i$ satisfy \eqref{Ricker-con-1}, i.e. $T^{[r_1,r_3]}\in \mathrm{CRC}(3)$ for all $0<r_1<\frac{1}{6},0<r_3<\frac{1}{4}$. It follows from Table \ref{biao0} {\rm(26)} that $T^{[r_1,r_3]}$ belongs to class {\rm 26}. $T^{[r_1,r_3]}$ has a unique positive fixed point $p=(\frac{80}{287},\frac{12}{287},\frac{212}{287})$. By numerical calculation, we find that $T^{[r_1,r_3]}$ admits a Chenciner bifurcation point at $p$ when $r_1\approx 0.026288$ and $r_3\approx 0.004706$, where the second Lyapunov coefficient $l_2=-0.1342<0$.
Therefore, a stable fixed point and an attracting (large) invariant closed curve, separated by an unstable invariant closed curve can coexist in class $26$ for $\mathrm{CRC}(3)$.
\end{example}

\section{Discussion}
\label{sec:Discussion}
This paper presents permanence and impermanence criteria for discrete-time dissipative Kolmogorov systems \eqref{map_T_dimn} (Theorem \ref{theorem:permanence}) and those admitting a carrying simplex $\Sigma$ (Theorem \ref{Sigma:permanence}), respectively. For three-dimensional maps admitting a carrying simplex, such criteria are finitely computable conditions which only depend on the nontrivial boundary fixed points (Corollary \ref{coro:3-permanence}).

The competitive systems induced by the maps \eqref{equ:T2} with linearly determined fixed points, i.e. all maps in the set $\mathrm{DCS}(n,f)$, always admit a carrying simplex. Particularly, we define an equivalence relation relative to local dynamics of nontrivial boundary fixed points for the set $\mathrm{DCS}(3,f)$ according to this linear structure. We say that two mappings in $\mathrm{DCS}(3,f)$ are equivalent if all their boundary
fixed points have the same local dynamics on the carrying simplices after a permutation of the indices $\{1, 2, 3\}$.
Via the index formula \eqref{equ:1}, which states that the sum of the indices of the fixed points on the carrying simplex is one, we list the stable equivalence classes for $\mathrm{DCS}(3,f)$ which are independent of generating functions $f\in \mathcal{F}_3$, and present the phase portraits on $\Sigma$. Specifically,

$\bullet$ there are always a total of $33$ stable equivalence classes, no matter what generating functions are, which are described in terms of inequalities on the parameters, and given in Table \ref{biao0};

$\bullet$  every nontrivial orbit converges to a fixed point on the boundary of $\Sigma$ in classes $1-18$;

 $\bullet$  each map in classes $19-25$ admits a unique positive fixed point $p$ which is a saddle, such that every nontrivial orbit converges to some fixed point on the boundary of the carrying simplex, except those on the stable manifold of $p$ which is a union of simple curves (see Remark \ref{remark:invariant-manifold});

 $\bullet$  {each map in classes $26-33$ has a unique positive fixed point $p$ with index $1$; $p$ is always a hyperbolic repeller in class $32$; and $p$ is globally asymptotically stable in class $33$; within classes $26-31$, Neimark-Sacker bifurcations might occur;}

 $\bullet$  there is a heteroclinic cycle in class $27$.

Applying our permanence and impermanence criteria to each class in $\mathrm{DCS}(3,f)$, we obtain that the systems in classes $29$, $31$, $33$ and class $27$ with a repelling heteroclinic cycle are permanent, while those in classes $1-26$, $28$, $30$, $32$ and class $27$ with an attracting heteroclinic cycle are impermanent; for systems in class $33$, the permanence can guarantee the global stability of the unique positive fixed point.

However, permanence {does not always imply} the global asymptotic stability of the unique positive fixed point $p$, and the local stability of $p$ depends on the generating function $f\in \mathcal{F}_3$ by \eqref{equ:DTp}.
Indeed, Neimark-Sacker bifurcations can happen in permanent classes $29$ and $31$ for the Leslie-Gower model, the generalized Atkinson-Allen model, the model with different types of growth functions, and the Ricker model. Neimark-Sacker bifurcations can also occur in class $27$ with repelling heteroclinic cycles for the Leslie-Gower model and the generalized Atkinson-Allen model.
So invariant cycles can occur in these classes, on which all orbits are periodic, or any orbit is dense. Numerical experiments show that Chenciner bifurcations can also happen in class $29$ for the generalized Atkinson-Allen model and the Ricker model, and in class $27$ with repelling heteroclinic cycles for the Leslie-Gower model, which means that two isolated invariant cycles can coexist on the carrying simplex for such systems.
In the impermanent classes, such as classes $26$, $28$, $30$ and class $27$ with attracting heteroclinic cycles, Neimark-Sacker bifurcations can also occur; see Section \ref{sec:application} and \cite{jiang2014,LG,GyllenbergCGAA,GyllenbergRicker} for more details. By the way, the dynamics in the same class which has a unique positive fixed point might be different for different kinds of generating functions $f\in \mathcal{F}_3$. For example, Neimark-Sacker bifurcations do not happen in classes $28$ and $30$ for the standard Atkinson-Allen model \cite{jiang2014}, while they can happen in these two classes for the Leslie-Gower model \cite{LG}.

Furthermore, the results imply that when all the boundary fixed points are unstable, the system may not be permanent, because impermanence can occur in class $27$ with attracting heteroclinic cycles, whose boundary fixed points are all unstable. When the system admits no heteroclinic cycle, i.e. it is not in class $27$, all the boundary fixed points being unstable implies the permanence for $T\in \mathrm{DCS}(3,f)$.

Biologically, the system is impermanent if one of the following conditions holds:

$\bullet$ there exists some species which cannot be invaded by any of the other two species (classes $1-3,7,8,13-23,26,28,30$ and $32$);

$\bullet$ there exists a two-species steady state which cannot be invaded by the third species (classes $4-6,9-12,24$ and $25$).\\
The system is permanent if the following conditions hold simultaneously (classes $29,31$ and $33$):

$\bullet$  each species can be invaded by at least one of the other two species;

 $\bullet$  there exists one species which can be invaded by both of the other two species;

 $\bullet$ any coexistence of two species can be invaded by the third species.

Such classification also presents a detailed classification for permanence and impermanence. Based on this, one can investigate the further long term dynamical properties within each of classes $26-32$. Finally, we propose some interesting open problems as follows.

 $\bullet$  {Give sufficient conditions to guarantee the global asymptotic stability of the positive fixed point for permanent systems in classes $29,31$ and class $27$ with repelling heteroclinic cycles.}

 $\bullet$ Investigate the nontrivial interesting dynamics, such as multiplicity of invariant cycles, in both permanent and impermanent systems.

\section*{Acknowledgments}
The authors are greatly indebted to two referees for the careful and patient reading of our original manuscript, many valuable comments and useful suggestions which led to much improvement in the presentation of our results.

\appendix
\section{Stable equivalence classes in $\mathrm{DCS}(3,f)$}
\begin{center}
  \begin{longtable}{c@{\extracolsep{\fill}}c@{\extracolsep{\fill}}c}
\caption{The $33$ equivalence classes in $\mathrm{DCS}(3,f)$, where $\gamma_{ij}=\mu_{ii}-\mu_{ji}$, $\beta_{ij}=\frac{\mu_{jj}-\mu_{ij}}{\mu_{ii}\mu_{jj}-\mu_{ij}\mu_{ji}}$ ($\beta_{ij}$ is well defined; see Remark \ref{remark-3d-1}), $i,j=1,2,3$ and $i\neq j$, and each $\Sigma$ is given by a representative map of that class. A fixed point is represented by a closed dot $\bullet$ if it attracts on $\Sigma$, by an open dot $\circ$ if it repels on $\Sigma$, and by the intersection of its stable and unstable manifolds if it is a saddle on $\Sigma$. For classes $1-25$ and $33$, every orbit converges to some fixed point; for classes $26-31$, Neimark-Sacker bifurcations might occur; for class $27$, $\partial \Sigma$ is a heteroclinic cycle; for class $32$, the unique positive fixed point is a repeller and Neimark-Sacker bifurcation cannot occur in this class.}\\[-2pt]
        \hline
         Class & Parameter conditions & Phase Portrait on $\Sigma$\\
        \hline
        \endfirsthead
        \caption[]{(continued)}\\
        \hline
        Class & Parameter conditions & Phase Portrait on $\Sigma$\\
        \hline
&&\\
        \endhead
        \hline
        \endfoot
        \endlastfoot
&&\\
1 &
\begin{tabular}{ll} &
$\gamma_{12}<0, \gamma_{13}<0, \gamma_{21}>0$,\\
& $\gamma_{23}>0, \gamma_{31}>0, \gamma_{32}<0$
\end{tabular}
&
    \parbox{2cm}{\includegraphics[width=2cm,height=1.6cm]{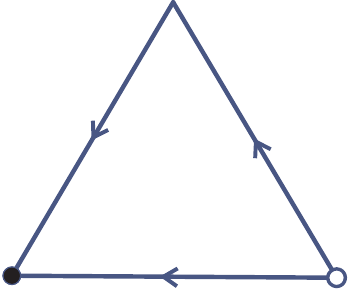}} \\

&&\\
    2 &
\begin{tabular}{@{}l@{~}l@{}} {
 (i)}& $\gamma_{12}<0, \gamma_{13}<0, \gamma_{21}<0$,\\
&$\gamma_{23}>0, \gamma_{31}>0, \gamma_{32}<0$\\
(ii)&$\mu_{31}\beta_{12}+\mu_{32}\beta_{21}<1$
 \end{tabular}
 &
    \parbox{2cm}{\includegraphics[width=2cm,height=1.6cm]{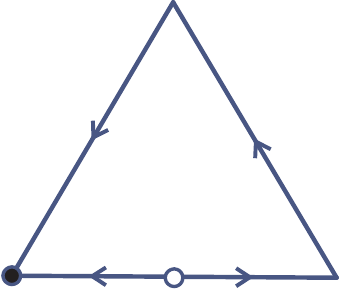}} \\
&&\\
    3 &
\begin{tabular}{ll}
 (i)&$\gamma_{12}<0, \gamma_{13}<0, \gamma_{21}>0$,\\
 &$\gamma_{23}<0, \gamma_{31}>0, \gamma_{32}<0$\\
(ii)&$\mu_{12}\beta_{23}+\mu_{13}\beta_{32}<1$
\end{tabular}
 &
    \parbox{2cm}{\includegraphics[width=2cm,height=1.6cm]{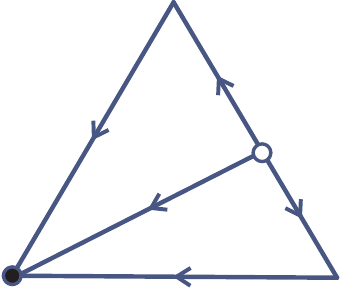}\vspace{2pt}}\\
&&\\
4 &
\begin{tabular}{ll}
 (i)& $\gamma_{12}>0, \gamma_{13}<0, \gamma_{21}>0$,\\
& $\gamma_{23}<0, \gamma_{31}>0, \gamma_{32}<0$\\
(ii)& $\mu_{12}\beta_{23}+\mu_{13}\beta_{32}<1$\\
(iii)& $\mu_{31}\beta_{12}+\mu_{32}\beta_{21}>1$
 \end{tabular}
 &
    \parbox{2cm}{\includegraphics[width=2cm,height=1.6cm]{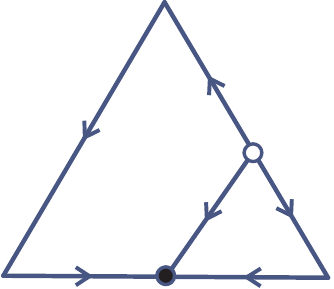}\vspace{2pt}} \\
&&\\
    5 &
\begin{tabular}{ll}
 (i)& $\gamma_{12}>0, \gamma_{13}>0, \gamma_{21}>0$,\\
& $\gamma_{23}<0, \gamma_{31}<0, \gamma_{32}>0$\\
(ii)& $\mu_{31}\beta_{12}+\mu_{32}\beta_{21}>1$
 \end{tabular}
 &
    \parbox{2cm}{\includegraphics[width=2cm,height=1.6cm]{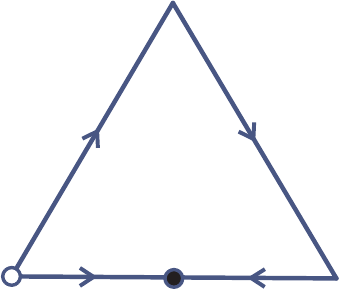}\vspace{2pt}} \\
    &&\\
6 &
\begin{tabular}{ll}
 (i)& $\gamma_{12}>0, \gamma_{13}>0, \gamma_{21}<0$,\\
& $\gamma_{23}>0, \gamma_{31}<0, \gamma_{32}>0$\\
 (ii)& $\mu_{12}\beta_{23}+\mu_{13}\beta_{32}>1$
 \end{tabular}
 &
    \parbox{2cm}{\includegraphics[width=2cm,height=1.6cm]{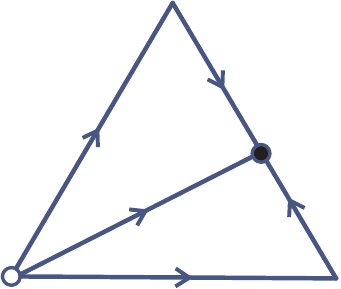}\vspace{2pt}} \\
&&\\
    7 &
\begin{tabular}{ll}
 (i)& $\gamma_{12}>0, \gamma_{13}>0, \gamma_{21}>0$,\\
& $\gamma_{23}>0, \gamma_{31}<0, \gamma_{32}<0$\\
 (ii)& $\mu_{31}\beta_{12}+\mu_{32}\beta_{21}<1$
 \end{tabular}
 &
    \parbox{2cm}{\includegraphics[width=2cm,height=1.6cm]{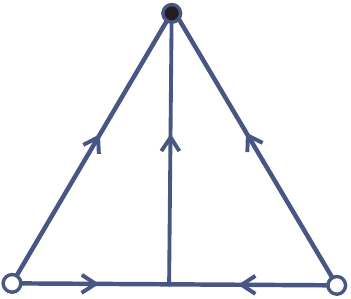}\vspace{2pt}} \\
&&\\
8 &
\begin{tabular}{ll}
 (i)& $\gamma_{12}>0, \gamma_{13}>0, \gamma_{21}>0$,\\
& $\gamma_{23}<0, \gamma_{31}<0, \gamma_{32}<0$\\
(ii)& $\mu_{12}\beta_{23}+\mu_{13}\beta_{32}<1$\\
(iii)& $\mu_{31}\beta_{12}+\mu_{32}\beta_{21}<1$
\end{tabular}
 &
    \parbox{2cm}{\includegraphics[width=2cm,height=1.6cm]{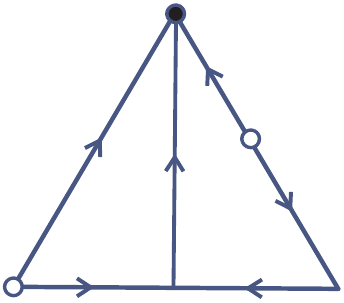}\vspace{2pt}} \\
&&\\
    9 &
\begin{tabular}{ll}
 (i)& $\gamma_{12}>0, \gamma_{13}>0, \gamma_{21}>0$,\\
& $\gamma_{23}>0, \gamma_{31}<0, \gamma_{32}>0$\\
 (ii)& $\mu_{12}\beta_{23}+\mu_{13}\beta_{32}>1$\\
(iii)& $\mu_{31}\beta_{12}+\mu_{32}\beta_{21}<1$
\end{tabular}
 &
    \parbox{2cm}{\includegraphics[width=2cm,height=1.6cm]{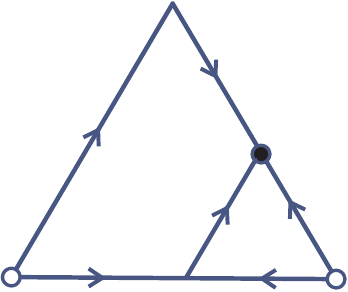}\vspace{2pt}} \\
&&\\
10 &
\begin{tabular}{ll}
 (i)& $\gamma_{12}>0, \gamma_{13}>0, \gamma_{21}>0$,\\
& $\gamma_{23}>0, \gamma_{31}<0, \gamma_{32}>0$\\
 (ii)& $\mu_{12}\beta_{23}+\mu_{13}\beta_{32}<1$\\
 (iii)& $\mu_{31}\beta_{12}+\mu_{32}\beta_{21}>1$
\end{tabular}
 &
    \parbox{2cm}{\includegraphics[width=2cm,height=1.6cm]{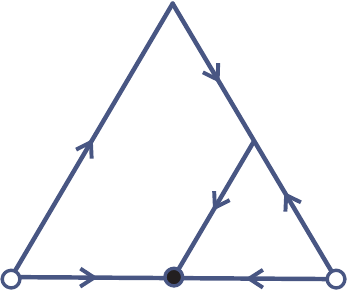}\vspace{2pt}} \\
&&\\
    11 &
\begin{tabular}{ll}
 (i)& $\gamma_{12}>0, \gamma_{13}>0, \gamma_{21}>0$,\\
& $\gamma_{23}<0, \gamma_{31}>0, \gamma_{32}<0$\\
 (ii)& $\mu_{12}\beta_{23}+\mu_{13}\beta_{32}<1$\\
 (iii)& $\mu_{21}\beta_{13}+\mu_{23}\beta_{31}<1$\\
 (iv)&  $\mu_{31}\beta_{12}+\mu_{32}\beta_{21}>1$
 \end{tabular}
 &
    \parbox{2cm}{\includegraphics[width=2cm,height=1.6cm]{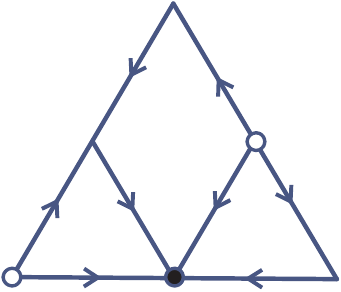}\vspace{2pt}} \\
&&\\
12 &
\begin{tabular}{ll}
 (i)& $\gamma_{12}>0, \gamma_{13}>0, \gamma_{21}>0$,\\
& $\gamma_{23}>0, \gamma_{31}>0, \gamma_{32}>0$\\
 (ii)& $\mu_{12}\beta_{23}+\mu_{13}\beta_{32}<1$\\
 (iii)& $\mu_{21}\beta_{13}+\mu_{23}\beta_{31}<1$\\
 (iv)& $\mu_{31}\beta_{12}+\mu_{32}\beta_{21}>1$
 \end{tabular}
 &
    \parbox{2cm}{\includegraphics[width=2cm,height=1.6cm]{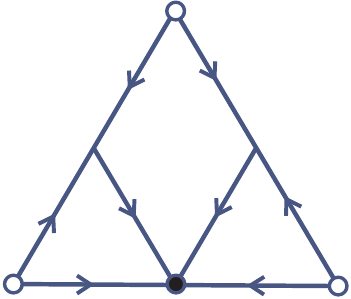}\vspace{2pt}} \\
&&\\
    13 &
\begin{tabular}{ll}
 (i)& $\gamma_{12}<0, \gamma_{13}<0, \gamma_{21}<0$,\\
& $\gamma_{23}<0, \gamma_{31}>0, \gamma_{32}>0$\\
 (ii)& $\mu_{31}\beta_{12}+\mu_{32}\beta_{21}>1$
 \end{tabular}
 &
    \parbox{2cm}{\includegraphics[width=2cm,height=1.6cm]{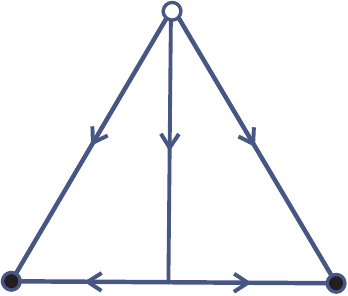}\vspace{2pt}} \\
&&\\
14 &
\begin{tabular}{ll}
 (i)& $\gamma_{12}<0, \gamma_{13}<0, \gamma_{21}<0$,\\
& $\gamma_{23}>0, \gamma_{31}>0, \gamma_{32}>0$\\
(ii)& $\mu_{12}\beta_{23}+\mu_{13}\beta_{32}>1$\\
 (iii)& $\mu_{31}\beta_{12}+\mu_{32}\beta_{21}>1$
 \end{tabular}
 &
    \parbox{2cm}{\includegraphics[width=2cm,height=1.6cm]{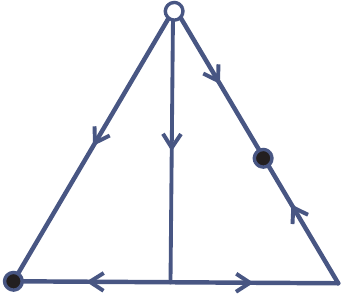}\vspace{2pt}} \\
&&\\
    15 &
\begin{tabular}{ll}
 (i)& $\gamma_{12}<0, \gamma_{13}<0, \gamma_{21}<0$,\\
 & $\gamma_{23}<0, \gamma_{31}>0, \gamma_{32}<0$\\
 (ii)& $\mu_{12}\beta_{23}+\mu_{13}\beta_{32}<1$\\
 (iii)& $\mu_{31}\beta_{12}+\mu_{32}\beta_{21}>1$
 \end{tabular}
 &
    \parbox{2cm}{\includegraphics[width=2cm,height=1.6cm]{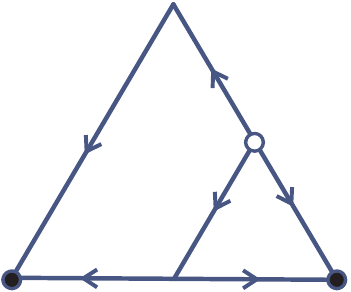}\vspace{2pt}} \\
&&\\
16 &
\begin{tabular}{ll}
 (i)& $\gamma_{12}<0, \gamma_{13}<0, \gamma_{21}<0$,\\
& $\gamma_{23}<0, \gamma_{31}>0, \gamma_{32}<0$\\
 (ii)& $\mu_{12}\beta_{23}+\mu_{13}\beta_{32}>1$\\
 (iii)& $\mu_{31}\beta_{12}+\mu_{32}\beta_{21}<1$
 \end{tabular}
 &
    \parbox{2cm}{\includegraphics[width=2cm,height=1.6cm]{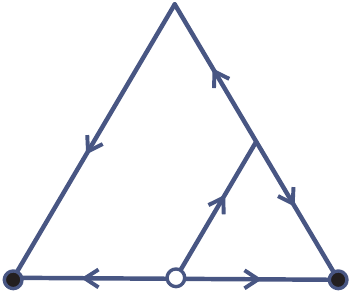}\vspace{2pt}} \\
&&\\
    17 &
\begin{tabular}{ll}
 (i)& $\gamma_{12}<0, \gamma_{13}<0, \gamma_{21}<0$,\\
& $\gamma_{23}>0, \gamma_{31}<0, \gamma_{32}>0$\\
 (ii)& $\mu_{12}\beta_{23}+\mu_{13}\beta_{32}>1$\\
 (iii)& $\mu_{21}\beta_{13}+\mu_{23}\beta_{31}>1$\\
 (iv)& $\mu_{31}\beta_{12}+\mu_{32}\beta_{21}<1$
 \end{tabular}
 &
    \parbox{2cm}{\includegraphics[width=2cm,height=1.6cm]{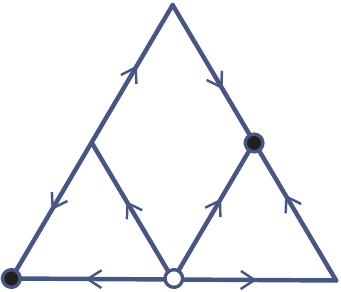}\vspace{2pt}} \\
&&\\
    18 &
\begin{tabular}{ll}
 (i)& $\gamma_{12}<0, \gamma_{13}<0, \gamma_{21}<0$,\\
& $\gamma_{23}<0, \gamma_{31}<0, \gamma_{32}<0$\\
 (ii)& $\mu_{12}\beta_{23}+\mu_{13}\beta_{32}>1$\\
 (iii)&  $\mu_{21}\beta_{13}+\mu_{23}\beta_{31}>1$\\
 (iv)&  $\mu_{31}\beta_{12}+\mu_{32}\beta_{21}<1$
 \end{tabular}
 &
    \parbox{2cm}{\includegraphics[width=2cm,height=1.6cm]{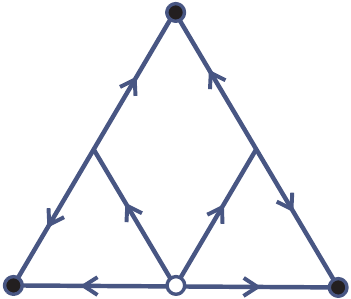}\vspace{2pt}}\\
&&\\
19 &
\begin{tabular}{ll}
 (i)& $\gamma_{12}>0, \gamma_{13}>0, \gamma_{21}<0$,\\
&  $\gamma_{23}<0, \gamma_{31}<0, \gamma_{32}<0$\\
 (ii)& $\mu_{12}\beta_{23}+\mu_{13}\beta_{32}<1$
 \end{tabular}
 &
    \parbox{2cm}{\includegraphics[width=2cm,height=1.6cm]{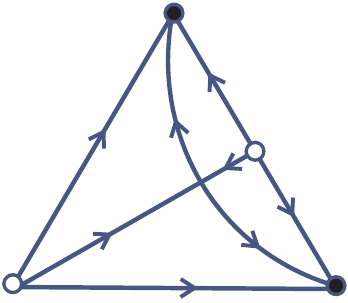}\vspace{2pt}} \\
&&\\
    20 &
\begin{tabular}{ll}
 (i)& $\gamma_{12}<0, \gamma_{13}<0, \gamma_{21}<0$,\\
& $\gamma_{23}<0, \gamma_{31}>0, \gamma_{32}<0$\\
 (ii)& $\mu_{12}\beta_{23}+\mu_{13}\beta_{32}<1$\\
 (iii)& $\mu_{31}\beta_{12}+\mu_{32}\beta_{21}<1$
 \end{tabular}
 &
    \parbox{2cm}{\includegraphics[width=2cm,height=1.6cm]{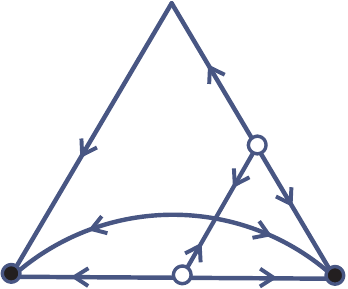}\vspace{2pt}} \\
&&\\
21 &
\begin{tabular}{ll}
 (i)& $\gamma_{12}<0, \gamma_{13}<0, \gamma_{21}<0$,\\
& $\gamma_{23}>0, \gamma_{31}<0, \gamma_{32}>0$\\
 (ii)& $\mu_{12}\beta_{23}+\mu_{13}\beta_{32}>1$\\
 (iii)& $\mu_{21}\beta_{13}+\mu_{23}\beta_{31}<1$\\
 (iv)&  $\mu_{31}\beta_{12}+\mu_{32}\beta_{21}<1$
 \end{tabular}
 &
    \parbox{2cm}{\includegraphics[width=2cm,height=1.6cm]{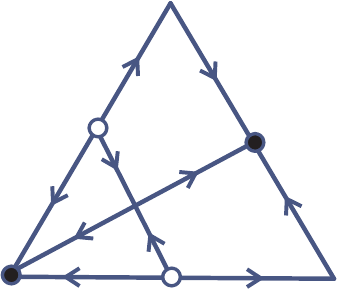}\vspace{2pt}} \\
&&\\
    22 &
\begin{tabular}{ll} {
 (i)}&{  $\gamma_{12}>0, \gamma_{13}>0, \gamma_{21}<0$,}\\
&{   $\gamma_{23}<0, \gamma_{31}>0, \gamma_{32}<0$}\\
{  (ii)}&{  $\mu_{12}\beta_{23}+\mu_{13}\beta_{32}<1$}\\
{  (iii)}& {  $\mu_{21}\beta_{13}+\mu_{23}\beta_{31}>1$
} \end{tabular}
 &
    \parbox{2cm}{\includegraphics[width=2cm,height=1.6cm]{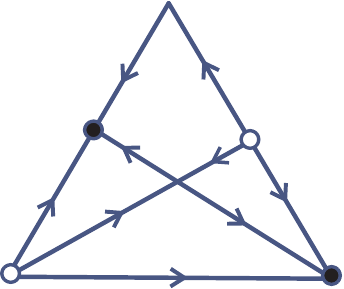}\vspace{2pt}} \\
&&\\
23 &
\begin{tabular}{ll} {
 (i)}&{  $\gamma_{12}>0, \gamma_{13}>0, \gamma_{21}>0$,}\\
&{   $\gamma_{23}>0, \gamma_{31}<0, \gamma_{32}<0$}\\
{  (ii)}&{  $\mu_{31}\beta_{12}+\mu_{32}\beta_{21}>1$
} \end{tabular}
 &
    \parbox{2cm}{\includegraphics[width=2cm,height=1.6cm]{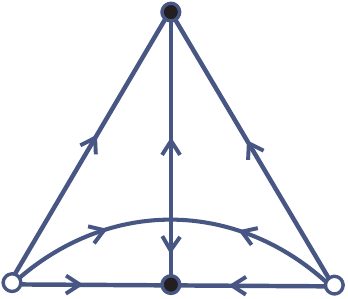}\vspace{2pt}} \\
&&\\
    24 &
\begin{tabular}{ll} {
 (i)}&{  $\gamma_{12}>0, \gamma_{13}>0, \gamma_{21}>0$,}\\
&{   $\gamma_{23}>0, \gamma_{31}<0, \gamma_{32}>0$}\\
{  (ii)}&{  $\mu_{12}\beta_{23}+\mu_{13}\beta_{32}>1$}\\
{  (iii)}&{  $\mu_{31}\beta_{12}+\mu_{32}\beta_{21}>1$
} \end{tabular}
&
    \parbox{2cm}{\includegraphics[width=2cm,height=1.6cm]{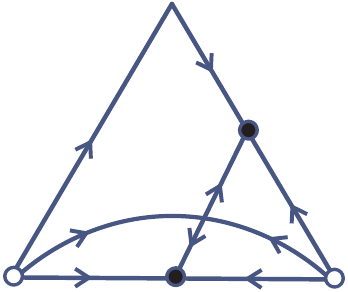}\vspace{2pt}} \\
&&\\
25 &
\begin{tabular}{ll} {
 (i)}&{  $\gamma_{12}>0, \gamma_{13}>0, \gamma_{21}>0$,}\\
&{   $\gamma_{23}<0, \gamma_{31}>0, \gamma_{32}<0$}\\
{  (ii)}&{  $\mu_{12}\beta_{23}+\mu_{13}\beta_{32}<1$}\\
{  (iii)}&{  $\mu_{21}\beta_{13}+\mu_{23}\beta_{31}>1$}\\
{  (iv)}&{  $\mu_{31}\beta_{12}+\mu_{32}\beta_{21}>1$
} \end{tabular}
 &
    \parbox{2cm}{\includegraphics[width=2cm,height=1.6cm]{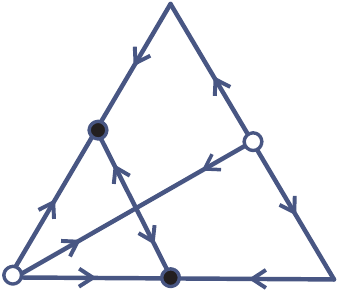}\vspace{2pt}} \\
&&\\
    26 &
\begin{tabular}{ll} {
 (i)}&{  $\gamma_{12}>0, \gamma_{13}>0, \gamma_{21}<0$,}\\
&{   $\gamma_{23}<0, \gamma_{31}>0, \gamma_{32}<0$}\\
{  (ii)}&{  $\mu_{12}\beta_{23}+\mu_{13}\beta_{32}>1$}\\
{  (iii)}&{  $\mu_{21}\beta_{13}+\mu_{23}\beta_{31}<1$
} \end{tabular}
 &
    \parbox{2cm}{\includegraphics[width=2cm,height=1.6cm]{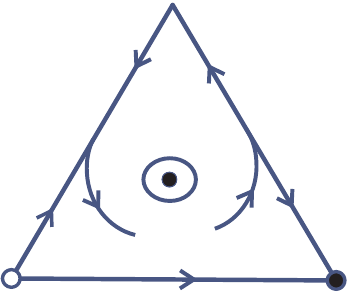}\vspace{2pt}} \\
&&\\
27 &
\begin{tabular}{ll} &{
 $\gamma_{12}>0, \gamma_{13}<0, \gamma_{21}<0$,}\\
&{   $\gamma_{23}>0, \gamma_{31}>0, \gamma_{32}<0$
} \end{tabular}
 &
    \parbox{2cm}{\includegraphics[width=2cm,height=1.6cm]{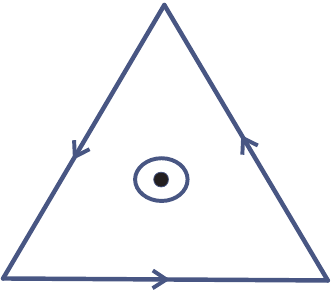}\vspace{2pt}} \\
&&\\
    28 &
\begin{tabular}{ll} {
 (i)}&{  $\gamma_{12}<0, \gamma_{13}<0, \gamma_{21}<0$,}\\
 &{  $\gamma_{23}>0, \gamma_{31}>0, \gamma_{32}<0$}\\
{  (ii)}&{   $\mu_{31}\beta_{12}+\mu_{32}\beta_{21}>1$
} \end{tabular}
&
    \parbox{2cm}{\includegraphics[width=2cm,height=1.6cm]{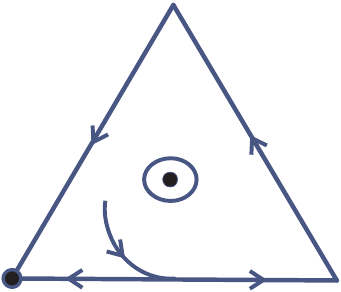}\vspace{2pt}} \\
&&\\
29 &
\begin{tabular}{ll} {
 (i)}&{  $\gamma_{12}>0, \gamma_{13}>0, \gamma_{21}>0$,}\\
&{   $\gamma_{23}<0, \gamma_{31}<0, \gamma_{32}>0$}\\
{  (ii)}&{  $\mu_{31}\beta_{12}+\mu_{32}\beta_{21}<1$
} \end{tabular}
 &
    \parbox{2cm}{\includegraphics[width=2cm,height=1.6cm]{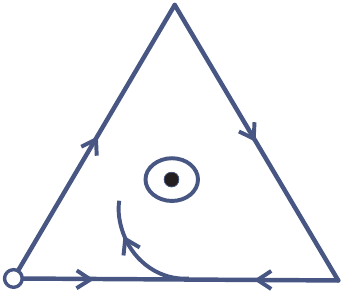}\vspace{2pt}} \\
&&\\
    30 &
\begin{tabular}{ll} {
(i)}&{  $\gamma_{12}<0, \gamma_{13}<0, \gamma_{21}<0$,}\\
&{  $\gamma_{23}<0, \gamma_{31}>0, \gamma_{32}<0$}\\
{  (ii)}&{  $\mu_{12}\beta_{23}+\mu_{13}\beta_{32}>1$}\\
{  (iii)}&{  $\mu_{31}\beta_{12}+\mu_{32}\beta_{21}>1$
} \end{tabular}
&
    \parbox{2cm}{\includegraphics[width=2cm,height=1.6cm]{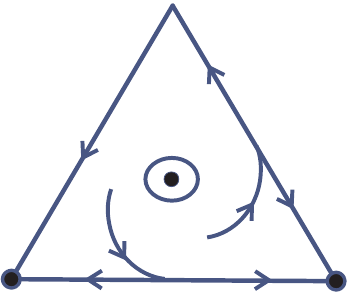}\vspace{2pt}} \\
&&\\
    31 &
\begin{tabular}{ll} {
(i)}&{  $\gamma_{12}>0, \gamma_{13}>0, \gamma_{21}>0$,}\\
&{  $\gamma_{23}>0, \gamma_{31}<0, \gamma_{32}>0$}\\
{  (ii)}&{  $\mu_{12}\beta_{23}+\mu_{13}\beta_{32}<1$}\\
{  (iii)}&{  $\mu_{31}\beta_{12}+\mu_{32}\beta_{21}<1$
} \end{tabular}
&
    \parbox{2cm}{\includegraphics[width=2cm,height=1.6cm]{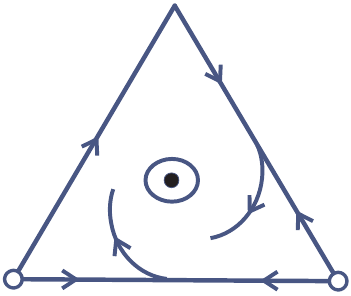}\vspace{2pt}} \\
&&\\
32 &
\begin{tabular}{ll} {
 (i)}&{  $\gamma_{12}<0, \gamma_{13}<0, \gamma_{21}<0$,}\\
&{   $\gamma_{23}<0, \gamma_{31}<0, \gamma_{32}<0$}\\
{  (ii)}&{  $\mu_{12}\beta_{23}+\mu_{13}\beta_{32}>1$}\\
{  (iii)}&{  $\mu_{21}\beta_{13}+\mu_{23}\beta_{31}>1$}\\
{  (iv)}&{  $\mu_{31}\beta_{12}+\mu_{32}\beta_{21}>1$
} \end{tabular}
&
    \parbox{2cm}{\includegraphics[width=2cm,height=1.6cm]{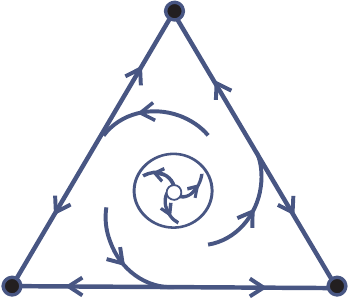}\vspace{2pt}} \\
&&\\
    33 &
\begin{tabular}{ll} {
 (i)}&{  $\gamma_{12}>0, \gamma_{13}>0, \gamma_{21}>0$,}\\
&{   $\gamma_{23}>0, \gamma_{31}>0, \gamma_{32}>0$}\\
{  (ii)}&{  $\mu_{12}\beta_{23}+\mu_{13}\beta_{32}<1$}\\
{  (iii)}&{  $\mu_{21}\beta_{13}+\mu_{23}\beta_{31}<1$}\\
{  (iv)}&{  $\mu_{31}\beta_{12}+\mu_{32}\beta_{21}<1$
} \end{tabular}
&
    \parbox{2cm}{\includegraphics[width=2cm,height=1.6cm]{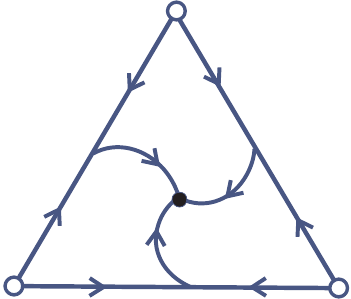}\vspace{2pt}} \\
[-8pt]\label{biao0}
\end{longtable}
\end{center}

\bibliographystyle{AIMS}

\medskip
Received xxxx 20xx; revised xxxx 20xx.
\medskip

\end{document}